\documentclass[11pt]{article}
\usepackage[authoryear,round]{natbib}
\usepackage{graphicx}
\usepackage{latexsym}
\usepackage[left=1.5in,top=1.5in,right=1.5in,bottom=1.5in]{geometry}\usepackage{amsmath}
\usepackage{amssymb}
\usepackage{booktabs}
\usepackage{color}
\usepackage{tikz}
\usetikzlibrary{shapes}
\usetikzlibrary{arrows}
\definecolor{darkblue}{rgb}{0,0.4,0.9}
\definecolor{gray10}{rgb}{0.1,0.1,0.1}
\definecolor{gray20}{rgb}{0.2,0.2,0.2}
\definecolor{gray30}{rgb}{0.3,0.3,0.3}
\definecolor{gray40}{rgb}{0.4,0.4,0.4}
\definecolor{gray60}{rgb}{0.6,0.6,0.6}
\definecolor{gray80}{rgb}{0.8,0.8,0.8}
\definecolor{gray90}{rgb}{0.9,0.9,.9}
\definecolor{gray95}{rgb}{0.95,0.95,.95}
\definecolor{gray96}{rgb}{0.96,0.96,.96}
\definecolor{lgreen} {RGB}{180,210,100}
\definecolor{dblue}  {RGB}{20,66,129}
\definecolor{ddblue} {RGB}{11,36,69}
\definecolor{lred}   {RGB}{220,0,0}
\definecolor{nred}   {RGB}{224,0,0}
\definecolor{norange}{RGB}{230,120,20}
\definecolor{nyellow}{RGB}{255,221,0}
\definecolor{ngreen} {RGB}{98,158,31}
\definecolor{dgreen} {RGB}{78,138,21}
\definecolor{nblue}  {RGB}{28,130,185}
\definecolor{jblue}  {RGB}{20,50,100}
\definecolor{nnyellow}{RGB}{235,200,0}
\definecolor{purple}{RGB}{150, 0, 120}
\definecolor{sgGreen} {RGB}{20, 180, 50}
\definecolor{revised}{rgb}{0,0,0.9}
\newtheorem{definition}{Definition}

\newtheorem{theorem}{Theorem}
\newtheorem{corollary}{Corollary}
\newtheorem{lemma}{Lemma}

\newcommand{\nl}{\newline}
\newcommand{\pl}{\parallel}
\newcommand{\openr}{\hbox{${\rm I\kern-.2em R}$}}
\newcommand{\openn}{\hbox{${\rm I\kern-.2em N}$}}

\title{Finite Sample Inference for Targeted Learning}

\author{Mark van der Laan \\ Division of Biostatistics, University of California, Berkeley\\ {\tt laan@berkeley.edu}
\\
} \date{\today}
 
\begin{document}
\maketitle
 
\begin{abstract}
The Highly-Adaptive-Lasso(HAL)-TMLE is an efficient estimator of a pathwise differentiable parameter in a statistical model that at minimal (and possibly only) assumes that the sectional variation norm of the true nuisance parameters are finite. It relies on an initial estimator (HAL-MLE) of the nuisance parameters by minimizing the empirical risk over the parameter space under the constraint that sectional variation norm is bounded by a constant, where this constant can be selected with cross-validation. In the formulation of the HAL-MLE this sectional variation norm corresponds with the sum of absolute value of coefficients for an indicator basis.   Due to its reliance on machine learning, statistical inference for the TMLE has been based on its normal limit distribution, thereby potentially ignoring a large second order remainder in finite samples.

In this article we present four methods for construction of a finite sample 0.95-confidence interval that use the nonparametric bootstrap to estimate the
finite sample distribution of the HAL-TMLE or a conservative distribution dominating the true finite sample distribution. We prove that it consistently estimates the optimal normal limit distribution, while its approximation error is driven by the performance of the bootstrap for a well behaved empirical process. We demonstrate our general
inferential methods for 1) nonparametric estimation of the average treatment
effect based on observing on each unit a covariate vector, binary treatment,
and outcome, and for 2) nonparametric estimation of the integral of the square of the multivariate density of the data distribution.

\end{abstract}


\section{Introduction}

We consider estimation of a pathwise differentiable real valued target parameter based on observing $n$ independent and identically distributed observations $O_1,\ldots,O_n$ with a data distribution $P_0$ known to belong in a highly nonparametric statistical model ${\cal M}$. 
A target parameter $\Psi:{\cal M}\rightarrow \openr$ is a  mapping that maps a possible data distribution $P\in {\cal M}$ into real number, while $\psi_0=\Psi(P_0)$ represents the answer to the question of interest about the data experiment. 
The canonical gradient $D^*(P)$ of the pathwise derivative of the target parameter at a $P$ defines an asymptotically efficient estimator among the class of regular estimators \citep{Bickeletal97}: An estimator $\psi_n$ is asymptotically efficient at $P_0$ if and only if it is asymptotically linear at $P_0$ with influence curve $D^*(P_0)$:
\[
\psi_n-\psi_0=\frac{1}{n}\sum_{i=1}^n D^*(P_0)(O_i)+o_P(n^{-1/2}).\]
The target parameter depends on the data distribution $P$ through a parameter $Q=Q(P)$, while the canonical gradient $D^*(P)$ possibly also depends on another nuisance parameter $G(P)$: $D^*(P)=D^*(Q(P),G(P))$. Both of these nuisance parameters are chosen so that they can be defined as a minimizer of the expectation of a specific loss function: $P L_1(Q(P))=\min_{P_1\in {\cal M}}PL_1(Q(P_1))$ and $PL_2(G(P))=\min_{P_1\in {\cal M}}P L_2(G(P_1))$, where we used the notation $Pf\equiv \int f(o)dP(o)$. We assume that the parameter spaces $Q({\cal M})=\{Q(P):P\in {\cal M}\}$ and $G({\cal M})=\{G(P):P\in {\cal M}\}$  for these nuisance parameters $Q$ and $G$ are contained in the set of multivariate cadlag functions with sectional variation norm $\pl \cdot\pl_v^*$ \citep{Gill&vanderLaan&Wellner95} bounded by a constant (this norm will be defined in the next section). 

We consider a targeted minimum loss-based (substitution) estimator $\Psi(Q_n^*)$  \citep{vanderLaan&Rubin06,vanderLaan08, vanderLaan&Rose11}
of the target parameter that uses as initial estimator   of these nuisance parameters $(Q_0,G_0)$ the highly adaptive lasso minimum loss-based estimators (HAL-MLE)  $(Q_n,G_n)$ defined by minimizing the empirical mean of the loss over the parameter space \citep{Benkeser&vanderLaan16}.  Since the HAL-MLEs converge at a rate faster than $n^{-1/2}$ w.r.t.\ the loss-based quadratic dissimilarities (which corresponds with  a rate faster than $n^{-1/4}$ for estimation of $Q_0$ and $G_0$), this HAL-TMLE has been shown to be asymptotically efficient under weak regularity conditions \citep{vanderLaan15}. Statistical inference could therefore be based on the normal limit distribution in which the asymptotic variance is estimated with an estimator of the variance of the canonical gradient. 
In that case, inference is ignoring the potentially very large contributions of the higher order remainder which could in finite samples easily dominate the first order empirical mean of the efficient influence curve term when the size of the nuisance parameter spaces is large (e.g., dimension of data is large and model is nonparametric). 

In this article we present four methods for inference that use the nonparametric bootstrap to estimate the  finite sample distribution of the HAL-TMLE or a  conservative distribution  dominating its true finite sample distribution. 
\subsection{Organization}

Firstly, in Section \ref{sectformulation} we formulate the estimation problem and motivate the challenge for statistical inference. We also provide  an easy to implement finite sample highly conservative confidence interval whose width converges to zero at the usual square-root sample size rate, but is not asymptotically sharp. We use this result to demonstrate the potential impact of the dimension of the data and sectional variation norm bound on the  width of a finite sample confidence interval. 

In Section \ref{sectnp} we  present the nonparametric bootstrap estimator of the actual sampling distribution of the HAL-TMLE which  thus incorporates estimation of its higher order stochastic behavior, and can thereby be expected to outperform the Wald-type confidence intervals.   We prove that this nonparametric bootstrap is asymptotically consistent for the optimal normal limit distribution. Our results also prove that the nonparametric bootstrap preserves the asymptotic behavior of the HAL-MLEs of our nuisance parameters $Q$ and $G$, providing further evidence for good performance of the  nonparametric bootstrap.
In the second subsection of Section \ref{sectnp} we propose to bootstrap the exact second-order expansion of the  HAL-TMLE. This results in a very direct estimator of the exact sampling distribution of the HAL-TMLE, although it comes at a cost of not respecting that the HAL-TMLE is a substitution estimator. 
 Importantly, our results demonstrate that the approximation error of the two nonparametric bootstrap estimates of the true finite sample distribution of the HAL-TMLE is mainly driven by the approximation  error of the nonparametric bootstrap for estimating the finite sample distribution of a well behaved empirical process. 
We suggest that these two  nonparametric bootstrap methods are the preferred methods for {\em accurate} inference, among our proposals, by not being aimed to be conservative.

In Section \ref{sectupperb} we  upper-bound  the absolute value of the exact remainder for the second-order expansion of the  HAL-TMLE in terms of a specified function of the loss-based dissimilarities for the HAL-MLEs of the nuisance parameters $Q$ and $G$. The resulting conservative finite sample second-order expansion is highly conservative  but is still asymptotically sharp by converging to the actual normal limit distribution of the HAL-TMLE (but from above). We then  propose to use the nonparametric bootstrap to estimate  this conservative finite  sample distribution. 
In the Appendix Section \ref{sectsupupperb} we further upper bound the  previously obtained conservative finite sample expansion by taking a supremum over a set of possible realizations of the HAL-MLEs  that will contain the true $Q_0$ and $G_0$ with probability tending to 1, where this probability is  controlled/set by the user.
We also propose a simplified conservative approximation of this supremum which is easy to implement. Even though these two sampling distributions are even more conservative they are  still asymptotically sharp, so that also the corresponding nonparametric bootstrap method is asymptotically converging to the optimal normal limit distribution.

In Section \ref{sectexample} we demonstrate our methods for two examples involving a nonparametric model and a specified target parameter (average treatment effect and integral of the square of the data density). 
We conclude with a discussion in Section \ref{sectdisc}. 
Some of the technical results and proofs  have been deferred to the Appendix, while the overall proofs are presented in the main part of the article.

\subsection{Why does it work, and how it applies to adaptive TMLE}

The key behind the validity of the nonparametric bootstrap for estimation of the sampling distribution of the HAL-MLE and HAL-TMLE is that the HAL-MLE is an actual MLE thereby avoiding data adaptive trade-off of bias and variance as naturally achieved with cross-validation. 
However, even though the inference is based on such a non-adaptive HAL-TMLE, one can still use an highly adaptive HAL-TMLE as point estimate in  our reported confidence intervals. Specifically, one can  use our confidence intervals with the point estimate defined as  a TMLE using a super-learner \citep{vanderLaan&Dudoit03,vanderVaart&Dudoit&vanderLaan06,vanderLaan&Dudoit&vanderVaart06,vanderLaan&Polley&Hubbard07,Chpt3} that includes the HAL-MLE as one of the candidate estimators in its library. By the oracle inequality for the cross-validation selector, such a super-learner will improve on the HAL-MLE so that the proposed inference based on the non-adaptive HAL-TMLE will be more conservative. In addition, our confidence intervals can be used with the point estimate defined by adaptive TMLEs incorporating additional refinements such as collaborative TMLE \citep{vanderLaan:Gruber10,Gruber:vanderLaan10a,Stitelman:vanderLaan10,vanderLaan:Rose11,Wang:Rose:vanderLaan11,gruber2012ijb}; cross-validated TMLE 
\citep{Zheng&vanderLaan12,vanderLaan&Rose11}; higher order TMLE \citep{Carone2014techreport,Caroneetal17,diaz2016ijb};  and double robust inference TMLE \citep{vanderLaan14a,Benkeseretal17}. Again, such refinements generally improve the finite sample accuracy of the estimator, so that  it will improve the coverage of the confidence intervals based on the non-adaptive HAL-TMLE.

Our confidence intervals can also be used if the statistical model ${\cal M}$ has  no known bound on the sectional variation norm of the nuisance parameters. In that case, we recommend  to select such a  bound with cross-validation (just as one selects the $L_1$-norm penalty in Lasso regression with cross-validation), which, by the oracle inequality for the cross-validation selector \citep{vanderLaan&Dudoit03,vanderVaart&Dudoit&vanderLaan06,vanderLaan&Dudoit&vanderVaart06}
is guaranteed to be larger that the  sectional variation norm of the true nuisance parameters $(Q_0,G_0)$  with probability tending to 1. In that case, the confidence intervals will still be asymptotically correct, incorporate most of the higher order variability, but ignores the potential finite sample underestimation of the true sectional variation norm.    In addition, in that case the inference adapts to the underlying unknown  sectional variation norm of the true nuisance parameters $(Q_0,G_0)$. We plan to evaluate the practical performance of our  methods in the near future.

\subsection{Relation to literature on higher order influence functions}
J. Pfanzagl \citep{pfanzagl1985} introduced the notion of higher order pathwise differentiability of finite dimensional target parameters and corresponding higher order gradients. He used these higher order expansions of the target parameter to define higher order one-step estimators that might result in asymptotically linear estimators where regular one-step estimators \citep{levit1975,ibragimov1981,pfanzagl1982,bickel1982annals} might fail to behave well due to a too large second-order remainder. 
 This is the perspective that inspired the seminal contributions of J. Robins, L. Li, E. Tchetgen \& A. van der Vaart (e.g., \citealp{robins2008imsnotes,robins2009metrika,li2011statsprobletters,vandervaart2014statscience}).  They develop a rigorous theory for (e.g.)  second-order one-step estimators, including the typical case that the parameter is not second-order pathwise differentiable. They allow the case that the second-order remainder asymptotically dominates the first order term, resulting in estimators and confidence intervals that converge to zero at a slower rate than $n^{-1/2}$. Their second-order expansion uses  approximations of ''would be'' second-order gradients, where the approximation results in a bias term they termed the representation error. Unfortunately, this representation error, due to the lack of second order pathwise differentiability, obstructs the construction of estimators with a third order remainder (and thereby asymptotic linearity under the condition that a third order term is $o_P(n^{-1/2})$) These second-order one-step estimators involve careful selection of tuning/smoothing parameters for approximating the ''would be'' second-order gradient in order to obtain an optimal bias-variance trade-off.  These authors applied their theory to nonparametric estimation of a mean with missing data and the  integral of the square of the density. 
 The higher-order expansions that come with the construction of higher order one-step estimators  can be  directly incorporated in
the construction of confidence intervals, thereby possibly leading to improved finite sample  coverage. These higher order expansions  rely on hard to estimate objects such as a multivariate density in a denominator, giving rise to enormous practical challenges to construct robust higher order confidence intervals, as noted in the above articles. 
 
 \cite{pfanzagl1982} already pointed out that the one-step estimators and till a larger degree higher order one-step estimators fail to respect global known bounds implied by the model and target parameter mapping, by adding to an initial estimator  an empirical mean of a first order influence function and higher order U-statistics (i.e. higher orders empirical averages) of higher order influence functions. He suggested that to circumvent this problem one would have to carry out the updating process in the model space instead of in the parameter space. This is precisely what is carried out by the general TMLE framework \citep{vanderLaan&Rubin06,vanderLaan&Rose11}, and higher order TMLE based on approximate higher order influence functions were developed in \citep{Carone2014techreport,diaz2016ijb,Caroneetal17}. These higher order TMLE represents the TMLE-analogue of higher order one-step estimators, just as the regular TMLE is an analogue of the regular one-step estimator. 
 These TMLEs automatically satisfy the known bounds and thus never produce non-sensical output such a negative number for a probability.  The higher order TMLE is just another TMLE but using  a least favorable submodel with an extra parameter, thereby providing a crucial safeguard against erratic behavior due to estimation of the higher order influence functions, while also being able to utilize the C-TMLE framework to select the tuning parameters for approximating these higher order influence functions.   
 

The approach in this article for construction of higher order confidence intervals is quite different from the construction of higher order one-step estimators or higher order TMLE and using the corresponding higher order expansion for inference. To start with, we use an asymptotically efficient HAL-TMLE so that we preserve the $n^{-1/2}$-rate of convergence, asymptotic normality and efficiency, even in nonparametric models that only assume that the true nuisance parameters have finite sectional variation norm. As point estimate we can still use an adaptive HAL-TMLE which can, for example, include the  higher-order HAL-TMLE refinement, beyond refinements mentioned above. However, for inference,  we avoid the  delicate higher order expansions based on approximate higher order gradients, but instead use the exact second-order expansion $\Psi(Q_n^*)-\Psi(Q_0)=(P_n-P_0)D^*(Q_n^*,G_n)+R_{20}(Q_n^*,G_n,Q_0,G_0)$  implied by the definition of the exact second-order remainder $R_{20}()$ (\ref{exactremainder}), which thus incorporates any higher order term. In addition, by using the robust HAL-MLE as estimators of $Q_0,G_0$,  the HAL-TMLE is not only efficient but  one can also use nonparametric bootstrap to estimate its sampling distribution. We then  use the nonparametric bootstrap to estimate the sampling distribution of HAL-TMLE itself, or its exact expansion, or an exact conservative expansion in which $R_{20}()$ is replaced by a robust upper bound which only depends on well behaved empirical processes for which the nonparametric bootstrap works  (again, due to using the HAL-MLE). Our confidence intervals have width of order $n^{-1/2}$ and are asymptotically sharp by converging to the optimal normal distribution based confidence interval as sample size increases. In addition, they are easy to implement as a by product of the computation of the HAL-TMLE itself.

\section{General formulation of statistical estimation problem and motivation for finite sample inference}\label{sectformulation}

\subsection{Statistical model and target parameter}
Let $O_1,\ldots,O_n$ be $n$ i.i.d. copies of a random variable $O\sim P_0\in {\cal M}$. 
Let $P_n$ be the empirical probability measure of $O_1,\ldots,O_n$.
Let $\Psi:{\cal M}\rightarrow\openr$ be a real valued parameter that is pathwise differentiable at each $P\in {\cal M}$ with 
canonical gradient $D^*(P)$.  That is, given a collection of one dimensional  submodels $\{P_{\epsilon}^S:\epsilon\}\subset {\cal M}$ through $P$ at $\epsilon =0$ with score $S$, for each of these submodels the derivative $\left . \frac{d}{d\epsilon}\Psi(P_{\epsilon}^S)\right |_{\epsilon =0}$ can be represented as $ E_P D(P)(O)S(O)$.
The latter is  an inner product of a gradient $D(P)\in L^2_0(P)$ with the score $S$ in the Hilbert space $L^2_0(P)$ of functions of $O$ with mean zero (under $P$) endowed with inner product $\langle S_1,S_2\rangle_P=P S_2S_2$. Let $\pl f\pl_P\equiv \sqrt{\int f(o)^2dP(o)}$ be the Hilbert space norm.
Such an element $D(P)\in L^2_0(P)$ is called a gradient of the pathwise derivative of $\Psi$ at $P$.
The canonical gradient  $D^*(P)$ is the unique gradient that is an element of the tangent space defined as the closure of the linear span of the collection of scores generated by this family of submodels. 

Define the exact second-order  remainder
\begin{equation}\label{exactremainder}
R_2(P,P_0)=\Psi(P)-\Psi(P_0)+(P-P_0)D^*(P),\end{equation}
where $(P-P_0)D^*(P)=-P_0D^*(P)$ since $D^*(P)$ has mean zero under $P$.

Let $Q:{\cal M}\rightarrow Q({\cal M})$ be a function valued parameter so that $\Psi(P)=\Psi_1(Q(P))$ for some $\Psi_1$.
For notational convenience, we will abuse  notation  by referring to the target parameter with $\Psi(Q)$ and $\Psi(P)$ interchangeably. 
Let $G:{\cal M}\rightarrow G({\cal M})$ be a function valued parameter so that $D^*(P)=D^*_1(Q(P),G(P))$ for some $D^*_1$. Again, we will use the notation $D^*(P)$ and $D^*(Q,G)$ interchangeably. 

Suppose that $O\in [0,\tau]\subset\openr^d_{\geq 0}$ is a $d$-variate random variable with support contained in a d-dimensional cube $[0,\tau]$.
Let $D_d[0,\tau]$ be the Banach space of $d$-variate real valued cadlag functions endowed with a supremum norm $\pl\cdot\pl_{\infty}$ \citep{Neuhaus71}.
Let $L_1:Q({\cal M})\rightarrow D_d[0,\tau]$ and $L_2:G({\cal M})\rightarrow D_d[0,\tau]$ be loss functions that identify the true $Q_0$ and $G_0$ in the sense that $P_0L_1(Q_0)=\min_{Q\in Q({\cal M})}P_0L_1(Q)$ and
$P_0L_2(G_0)=\min_{G\in G({\cal M})} P_0L_2(G)$. Let $d_{01}(Q,Q_0)=P_0L_1(Q)-P_0L_1(Q_0)$ and $d_{02}(G,G_0)=P_0L_2(G)-P_0L_2(G_0)$ be the loss-based dissimilarities for these two nuisance parameters.

{\bf Loss functions and canonical gradient have a uniformly bounded sectional variation norm:}
We assume that these loss functions and the canonical gradient map into functions in $D_d[0,\tau]$ with a sectional variation norm bounded by some universal finite constant:
\begin{eqnarray}
M_1\equiv \sup_{P\in {\cal M}}\pl L_1(Q(P))\pl_v^* &<&\infty\nonumber \\
 M_2\equiv \sup_{P\in {\cal M}}\pl L_2(G(P))\pl_v^*&<& \infty\nonumber \\
 M_3\equiv \sup_{P\in {\cal M}}\pl D^*(P)\pl_v^*&<& \infty.
  \label{sectionalvarbound}
 \end{eqnarray}
 

 For a given function $F\in D_d[0,\tau]$, we define the sectional variation norm as follows. For a given subset $s\subset\{1,\ldots,d\}$,  let $F_s(x_s)=F(x_s,0_{-s})$ be the $s$-specific section of $F$ that sets the coordinates outside the subset $s$ equal to 0, where we used the notation $(x_s,0_{-s})$ for the vector whose $j$-th component equals $x_j$ if $j\in s$ and $0$ otherwise. 
 The sectional variation norm is now defined by
 \[
 \pl F\pl_v^*=\mid F(0)\mid +\sum_{s\subset\{1,\ldots,d\}}\int_{(0_s,\tau_s]}\mid dF_s(u_s)\mid ,\]
 where the sum is over all subsets $s$ of $\{1,\ldots,d\}$.
 Note that $\int_{(0_s,\tau_s]}\mid dF_s(u_s)\mid $ is the standard variation norm of the measure $dF_s$ generated by its $s$-specific section $F_s$ on the $\mid s\mid$-dimensional edge $(0_s,\tau_s]\times \{0_{-s}\}$ of the $d$-dimensional cube $[0,\tau]$. Thus, the sectional variation norm is the sum of the variation of $F$ itself and of all its $s$-specific sections, plus $F(0)$.
 We also note that any function $F\in D_d[0,\tau]$ with finite sectional variation norm (i.e., $\pl F\pl_v^*<\infty$) can be represented as follows \citep{Gill&vanderLaan&Wellner95}:
 \begin{equation}\label{Frepresentation}
 F(x)=F(0)+\sum_{s\subset\{1,\ldots,d\}}\int_{(0_s,x_s]} dF_s(u_s) .\end{equation}
As utilized in \citep{vanderLaan15} to define the HAL-MLE,  since $\int_{(0_s,x_s]}dF_s(u_s)=\int I_{u_s\leq x_s}dF_s(u_s)$, this representation shows that $F$ can be written as an infinite linear combination of
 $s$-specific indicator basis functions $x\rightarrow I_{u_s\leq x_s}$ indexed by  a cut-off $u_s$, across all subsets $s$, where the coefficients in front of the indicators are equal to the infinitesimal increments $dF_s(u_s)$ of $F_s$ at $u_s$. For discrete measures $F_s$ this integral becomes a finite linear combination of such $\mid s\mid$-way indicators. One could think of this representation as a saturated model of a function $F$ in terms of single way indicators, two-way indicators, etc, till the final $d$-way indicator basis functions. 
For a function $f\in D_d[0,\tau]$,  we also define the supremum norm $\pl f\pl_{\infty}=\sup_{x\in [0,\tau]}\mid f(x)\mid$.

 
 {\bf Assuming that parameter spaces for $Q$ and $G$ are cartesian products of sets of cadlag functions with bounds on sectional variation norm:}
Although the above bounds $M_1,M_2,M_3$ are the only relevant bounds for the asymptotic performance of the HAL-MLE and HAL-TMLE, for practical formulation of a model ${\cal M}$ one might prefer to state the sectional variation norm restrictions on the parameters $Q$ and $G$ themselves.
For that purpose, let's assume that $Q=(Q_1,\ldots,Q_{K_1})$ for variation independent parameters $Q_k$  that are themselves $m_{1k}$-dimensional cadlag functions on $[0,\tau_{1k}]\subset \openr^{m_{1k}}_{\geq 0}$ with sectional variation norm bounded by some upper-bound $C_{1k}^u$ and lower bound $C_{1k}^l$, $k=1,\ldots,K_1$, and similarly  for $G=(G_1,\ldots,G_{K_2})$ with sectional variation norm bounds
$C_{2k}^u$ and $C_{1k}^l$, $k=1,\ldots,K_2$.
Typically, we have $C_{1k}^l=0$.
Specifically, let
\begin{eqnarray*}
{\cal F}_{1k}\equiv Q_k({\cal M})\\
{\cal F}_{2k}\equiv G_k({\cal M}),
\end{eqnarray*}
denote the parameter spaces for $Q_k$ and $G_k$,
and assume that these parameter spaces ${\cal F}_{jk}$ are contained in the class ${\cal F}_{jk}^{np}$ of $m_{jk}$-variate cadlag functions with sectional variation norm bounded from above by $C_{jk}^u$ and from below by $C_{jk}^l$, $k=1,\ldots,K_j$, $j=1,2$.
These bounds $C_1^u=(C_{1k}^u:k)$ and $C_2^u=(C_{2k}^u:k)$ will then imply bounds $M_1,M_2,M_3$. 
In such a setting, $L_1(Q)$ would  be defined as a sum loss function $L_1(Q)=\sum_{k=1}^{K_1}L_{1k}(Q_k)$
and $L_2(G)=\sum_{k=1}^{K_2}L_{2k}(G_k)$. We also define the vector losses ${\bf L}_1(Q)=(L_{1k}(Q_k):k=1,\ldots,K_1)$, ${\bf L}_2(G)=(L_{2k}(G_k):k=1,\ldots,K_2)$,
and corresponding vector dissimilarities ${\bf d}_{01}(Q,Q_0)=(d_{01,k}(Q_k,Q_{k0}):k=1,\ldots,K_1)$ and ${\bf d}_{02}(G,G_0)=(d_{02,k}(G_k,G_{k0}):k=1,\ldots,K_2)$.

In a typical case we would have that the parameter space ${\cal F}_{jk}$ of $Q_k$ ($j=1$) or $G_k$ ($j=2$) would be equal to
\begin{equation}\label{calFmodel}
{\cal F}_{jk,A_{jk}}^{np}\equiv \{F\in {\cal F}_{jk}^{np}: dF_s(u_s)=I_{(s,u_s)\in A_{jk}}dF_s(u_s)\},
\end{equation}  for some set $A_{jk}$ of possible values for $(s,u_s)$, $k=1,\ldots,K_j$, $j=1,2$,
where one evaluates this restriction on  $F$ in terms of the  representation (\ref{Frepresentation}).
Note that we used short-hand notation $g(x)=I_{x\in A} g(x)$ for $g$ being zero for $x\not \in A$.
We will make the convention that if $A$ excludes $0$, then it corresponds with assuming $F(0)=0$.

This subset  ${\cal F}_{1k,A_{1k}}^{np}$ of all cadlag functions ${\cal F}_{1k}^{np}$ with sectional variation norm smaller than $C_{1k}^u$ further  restricts the support of these functions to a set $A_{1k}$.
For example, $A_{1k}$ might set $dF_s=0 $ for subsets $s$ of size larger than $3$ for all values $u_s\in (0_s,\tau_s]$, in which case one assumes that the nuisance parameter $Q_k$ can be represented as a sum over all subsets $s$ of size $1,2$ and $3$ of a  function of the variables indicated by $s$.

In order to allow modeling of monotonicity (e..g, nuisance parameter $Q_k$ is an actual cumulative distribution function), we also allow that this set restricts $dF_s(u_s)\geq 0$ for all $(s,u_s)\in A_{jk}$. We will denote the latter parameter space with 
\begin{equation}\label{calFmodelplus}
{\cal F}_{jk,A_{jk}}^{np,+}=\{F\in {\cal F}_{jk}^{np}: dF_s(u_s)=I_{(s,u_s)\in A_{jk}}dF_s(u_s), dF_s\geq 0, F(0)\geq 0\}.
\end{equation}
For the parameter space (\ref{calFmodelplus}) of monotone functions we allow that the sectional variation norm is known by setting $C_{jk}^u=C_{jk}^l$ (e.g, for the class of cumulative distribution functions we would have
$C_{jk}^u=C_{jk}^l=1$),  while for the parameter space (\ref{calFmodel}) of cadlag functions with sectional variation norm between $C_{jk}^l$ and $C_{jk}^u$ we assume 
$C_{jk}^l<C_{jk}^u$.


Although not necessary at all, for the analysis of our proposed  nonparametric bootstrap sampling distributions  {\em we  assume this extra structure} that ${\cal F}_{jk}={\cal F}_{jk,A_{jk}}^{np}$ or ${\cal F}_{jk}={\cal F}_{jk,A_{jk}}^{np,+}$ for some set $A_{jk}$, $k=1,\ldots,K_j$, $j=1,2$.
This extra structure allows us to obtain concrete results for the validity of the nonparametric bootstrap for the HAL-MLEs $Q_n$ and $G_n$ defined below, and thereby the HAL-TMLE (see Appendix \ref{AppendixB}).
In addition, the implementation of the HAL-MLE for such a parameter space ${\cal F}_{jk,A_{jk}}^{np}$ still corresponds with fitting a linear combination of indicator basis functions
$I_{u_s\leq x_s}$
under the sole constraint that the sum of the absolute value of the coefficients is bounded by $C_{jk}^u$ (and possibly from below by $C_{jk}^l$), and possibly that the coefficients are non-negative, where the set $A_{jk}$ implies the set of indicator basis functions that are included.
Specifically, in the case that the nuisance parameter is a conditional mean we can compute the  HAL-MLE with standard lasso regression software \citep{Benkeser&vanderLaan16}. Therefore, this restriction on our set of models allows straightforward computation of its HAL-MLEs and  corresponding HAL-TMLE.

Thus, a typical statistical model would be of the form ${\cal M}=\{P: Q_{k_1}(P)\in {\cal F}_{1k_1,A_{1k_1}}^{np},G_{k_2}(P)\in {\cal F}_{2k_2,A_{2k_2}}^{np},k_1,k_2\}$ for sets $A_{1k_1},A_{2k_2}$, but the model might  include additional restrictions on $P$ beyond restricting  the variation independent components of $Q(P)$ and $G(P)$ to be elements of these sets ${\cal F}_{jk_j,A_{jk_j}}^{np}$, as long as their parameter spaces equal these sets ${\cal F}_{jk_j,A_{jk_j}}^{np}$ or ${\cal F}_{jk_j,A_{jk_j}}^{np,+}$.
\paragraph{Remark regarding creating  nuisance parameters with parameter space of type (\ref{calFmodel}) or (\ref{calFmodelplus}):}
In our first example we have a nuisance parameter ${G}(W)=E_P(A\mid W)$ that is not just assumed to be cadlag and have bounded sectional variation norm but is also bounded between $\delta$ and $1-\delta$ for some $\delta>0$. This means that the parameter space for this ${G}$ is not exactly of type (\ref{calFmodel}). 
This is easily resolved by reparameterizing ${G}=\delta +(1-2\delta)\mbox{expit}(f(W))$ where $f$ can be any cadlag function with sectional variation norm bounded by some constant. One now defines the nuisance parameter as $f({G})$ instead of ${G}$ itself. Similarly, in our second example, $Q$ is the data density $p$ itself, which is assumed to be bounded from below by a $\delta\geq 0$ and from above by an $M<\infty$, beyond being cadlag and having a bound on the sectional variation norm. 
In this case, we could parameterize $p$ as $p(o)=c(f)\{ \delta+(M-\delta)\mbox{expit}(f(o))\}$, where $c(f)$ is the normalizing constant guaranteeing that $\int p(o)d\mu(o)=1$. One now defines the nuisance parameter as $f(Q)$  instead of $Q$ itself.  These just represent a few examples showcasing that one can reparametrize the natural nuisance parameters $Q$ and $G$ in terms of nuisance parameters that have a parameter space of the form (\ref{calFmodel}) or (\ref{calFmodelplus}). These representations are actually natural steps for the implementation of the HAL-MLE since they allow us now to minimize the empirical risk over a linear model with the sole constraint that the sum of absolute value of coefficients is bounded  (and possibly coefficients are non-negative).

{\bf Bounding the exact second-order remainder in terms of loss-based dissimilarities:}
Let \[
R_2(P,P_0)=R_{20}(Q,G,Q_0,G_0)\]
for some mapping $R_{20}()=R_{2P_0}()$ possibly indexed by $P_0$. We often have that 
$R_{20}(Q,G,Q_0,G_0)$ is a sum of second-order terms of the types $\int (H_1(Q)-H_1(Q_0))^2 f(P,P_0)dP_0$, $\int (H_2(G)-H_2(G_0))^2 f(P,P_0)dP_0$ and
$\int (H_1(Q)-H_1(Q_0))(H_2(G)-H_2(G_0))f(P,P_0) dP_0$ for certain specifications of $H_1, H_2 $ and $f()$.
Specifically, in all our applications it has the form $\int R_2(Q,G,Q_0,G_0) dP_0$ for some quadratic function $R_2(Q,G,Q_0,G_0)$.
If it only involves terms of the third type, then $R_2(P,P_0)$ has a double robust structure allowing the construction of double robust estimators whose consistency relies on consistent estimation  of either $Q$ or $G$. In particular, in that case the HAL-TMLE is double robust as well.

We assume the following upper bound:
\begin{equation}\label{boundingR2}\mid R_2(P,P_0)\mid =\mid R_{20}(Q,G,Q_0,G_0)\mid \leq f({\bf d}_{01}^{1/2}(Q,Q_0),{\bf d}_{02}^{1/2}(G,G_0))\end{equation}
for some function $f:\openr^K_{\geq 0}\rightarrow\openr_{\geq 0}$, $K=K_1+K_2$, of the form $f(x)=\sum_{i,j} a_{ij} x_ix_j$, a quadratic polynomial with positive coefficients $a_{ij}\geq 0$.  In all our examples, one simply uses the Cauchy-Schwarz inequality to bound $R_{20}(P,P_0)$ in terms of $L^2(P_0)$-norms of $Q_{k_1}-Q_{k_10}$ and $G_{k_2}-G_{k_20}$, and subsequently one relates these  $L^2(P_0)$-norms to its loss-based dissimilarities $d_{01,k_1}(Q_{k_1},Q_{k_10})$ and $d_{02,k_2}(G_{k_2},G_{k_20})$, respectively. This bounding step will also rely on a positivity assumption so that denominators in $R_{20}(P,P_0)$ are uniformly bounded away from zero.

{\bf Continuity of efficient influence curve as function of $P$:}
We also assume a basic uniform continuity condition on the efficient influence curve:
\begin{equation}\label{contDstar}
\sup_{P\in {\cal M} }\frac{P_0\{D^*(P)-D^*(P_0)\}^2}{ d_{01}(Q(P),Q_0)+d_{02}(G(P),G_0)}<\infty .\end{equation}
The above two uniform bounds (\ref{boundingR2}) and (\ref{contDstar}) on the model ${\cal M}$  will generally hold under a strong positivity assumption that guarantees that there are no nuisance parameters (e.g., a parameter of $G$) in the denominator of $D^*(P)$ and 
$R_2(P,P_0)$ that can be arbitrarily close to 0 on the support of $P_0$. 

\subsection{HAL-MLEs of nuisance parameters}
We estimate $Q_0,G_0$ with HAL-MLEs $Q_n,G_n$ satisfying
\begin{eqnarray*}
P_n L_1(Q_n)&=&\min_{Q\in Q({\cal M})} P_n L_1(Q)\\
P_n L_2(G_n)&=& \min_{G\in G({\cal M})}P_n L_2(G).
\end{eqnarray*}
Due to the sum-loss and variation independence of the components of $Q$ and $G$, these HAL-MLEs correspond with separate HAL-MLEs for each component. 
We have the following previously established result \citep{vanderLaan15} for these HAL-MLEs. We represent estimators as mappings on the nonparametric model ${\cal M}_{np}$ containing all possible realizations of the empirical measure $P_n$.
\begin{lemma}\label{lemma2}
Let $O\sim P_0\in {\cal M}$. Let $Q:{\cal M}\rightarrow Q({\cal M})$ be a function valued parameter and let $L:Q({\cal M})\rightarrow D_d[0,\tau]$ be a loss function so that $Q_0\equiv Q(P_0)=\arg\min_{Q\in Q({\cal M})}P_0L(Q)=\arg\min_{Q\in Q({\cal M})} \int L(Q)(o)dP_0(o)$. 
 Let $\hat{Q}:{\cal M}_{np}\rightarrow Q({\cal M})$ be an estimator  $Q_n\equiv \hat{Q}(P_n)$ so that $P_n L_1(Q_n)=\min_{Q\in Q({\cal M})}P_n L(Q)$. 
 Let $d_0(Q,Q_0)=P_0L(Q)-P_0L(Q_0)$ be the loss-based dissimilarity. 
 Then,
 \[
 d_0(Q_n,Q_0)\leq -(P_n-P_0)\{L(Q_n)-L(Q_0)\}. 
 \]
 If $\sup_{Q\in Q({\cal M})}\pl L(Q)\pl_v^*<\infty$, then 
 \[
E_0 d_0(Q_n,Q_0)=O(n^{-1/2-\alpha(d)}),\]
where $\alpha(d)=1/(2d+4)$.
 \end{lemma}
 Application of this general lemma proves that $d_{01}(Q_n,Q_0)=O_P(n^{-1/2-\alpha(d)})$ and $d_{02}(G_n,G_0)=O_P(n^{-1/2-\alpha(d)})$. 
It also shows that we have the following actual empirical process upper-bounds:\begin{eqnarray*}
d_{01}(Q_n,Q_0)&\leq&-(P_n-P_0)L_1(Q_n,Q_0)\\
d_{02}(G_n,G_0)&\leq&-(P_n-P_0)L_2(G_n,G_0),
\end{eqnarray*}
where we defined $L_1(Q,Q_0)\equiv L_1(Q)-L_1(Q_0)$ and $L_2(G,G_0)\equiv L_2(G)-L_2(G_0)$. These upper bounds will be utilized in our proposed conservative sampling distributions of the HAL-TMLE in Appendix \ref{sectsupupperb}.

\paragraph{Super learner including HAL-MLE outperforms HAL-MLE}
Suppose that we estimate $Q_0$ and $G_0$ instead with super-learners $\tilde{Q}_n,\tilde{G}_n$  in which the library of the super-learners contains this HAL-MLE $Q_n$ and $G_n$. Then, by the oracle inequality for the super-learner, we know that 
$d_{01}(\tilde{Q}_n,Q_0)$ and $d_{02}(\tilde{G}_n,G_0)$ will be asymptotically equivalent with the oracle selected estimator, so that $d_{01}(Q_n,Q_0)$ and $d_{02}(G_n,G_0)$ represent asymptotic upper bounds for $d_{01}(\tilde{Q}_n,Q_0)$ and $d_{02}(\tilde{G}_n,G_0)$ \citep{vanderLaan15}.
In addition, practical experience has demonstrated that the super-learner outperforms its library candidates in finite samples. 
Therefore, assuming that each estimator in the library of the super-learners for $Q_0$ and $G_0$  falls in the parameter spaces ${\cal F}_1$ and ${\cal F}_2$ of $Q$ and $G$, respectively, our proposed  estimators of the sampling distribution of the HAL-TMLE can also be used to construct a confidence interval around the super-learner based TMLE. These width of these confidence intervals are not adapting to possible superior performance of the  super-learner and could thus be overly conservative in case the super-learner outperforms the HAL-MLE.

\subsection{HAL-TMLE}
Consider a finite dimensional local  least favorable model $\{Q_{n,\epsilon}:\epsilon\}\subset Q({\cal M})$ through $Q_n$ at $\epsilon =0$ so that  the linear span of the components of $\frac{d}{d\epsilon}L_1(Q_{n,\epsilon})$ at $\epsilon =0$ includes $D^*(Q_n,G_n)$. 
Let $Q_n^*=Q_{n,\epsilon_n}$ for $\epsilon_n=\arg\min_{\epsilon}P_n L_1(Q_{n,\epsilon})$. We assume that this one-step TMLE  $Q_n^*$ already satisfies
\begin{equation}\label{efficeqn}
r_n\equiv \mid P_n D^*(Q_n^*,G_n)\mid=o_P(n^{-1/2}).\end{equation}
 As shown in \citep{vanderLaan15} this holds for the one-step HAL-TMLE under regularity conditions. Alternatively,  one could use the one-dimensional canonical universal least favorable model satisfying $\frac{d}{d\epsilon}L_1(Q_{n,\epsilon})=D^*(Q_{n,\epsilon},G_n)$ at each $\epsilon$ (see our second example in Section \ref{sectexample}). In that case, the efficient influence curve equation (\ref{efficeqn}) is solved exactly with the one-step TMLE: i.e.,  $r_n=0$ \citep{vanderLaan&Gruber15}. 
 The HAL-TMLE of $\psi_0$ is now the plug-in estimator $\psi_n^*=\Psi(Q_n^*)$. Sometimes, we will refer to this estimator as the HAL-TMLE$(C^u)$ to indicate its dependence on the specification of $C^u=(C_1^u,C_2^u)$.
 
In the Appendix \ref{AppendixA} we show that under  smoothness condition on the least favorable submodel (as function of $\epsilon$)   $d_{01}(Q_{n,\epsilon_n},Q_0)$ converges at the same rate as $d_{01}(Q_n,Q_0)=O_P(n^{-1/2-\alpha(d)})$ (see (\ref{thalmle})). This also implies this result for any $K$-th step TMLE with $K$ fixed.
The advantage of  a one-step or $K$-th step TMLE is that it is always well defined, and it easily follows that it converges at the same rate as the initial $Q_n$ to $Q_0$. Even though we derive some more explicit results for the one-step TMLE (and thereby $K$-th step TMLE), our results are presented so that they can   be applied  to any TMLE $Q_n^*$, including iterative TMLE, but we then simply assume that it has been shown that $d_{01}(Q_n^*,Q_0)$ converges at same rate to zero as $d_{01}(Q_n,Q_0)$.

It is assumed that for any $Q$ in its parameter space $\sup_{\epsilon}\pl Q_{\epsilon}\pl_v^*< C \pl Q\pl_v^*$ for some $C<\infty$ so that the least favorable model preserves the bound on the sectional variation norm. Since the HAL-MLE $Q_n$ has the maximal allowed uniform sectional variation norm $C_1^u$, it is likely that $Q_n^*$ has a slightly larger
variation norm than this bound. 

\subsection{Asymptotic efficiency theorem for HAL-TMLE and CV-HAL-TMLE}
The bound  ${\bf d}_{01}(Q_n^*,Q_0)=O_P({\bf d}_{01}(Q_n,Q_0))$, the rate results for ${\bf d}_{01}(Q_n,Q_0)$ and ${\bf d}_{02}(G_n,G_0)$ implied by Lemma \ref{lemma2},  combined with (\ref{boundingR2}), now shows that the second-order term $R_{20}(Q_n^*,G_n,Q_0,G_0)=O_P(n^{-1/2-\alpha(d)})$.

We have the following identity for the HAL-TMLE:
\begin{eqnarray}\nonumber
\Psi(Q_n^*)-\Psi(Q_0)&=&(P_n-P_0)D^*(Q_n^*,G_n)+R_{20}(Q_n^*,G_n,Q_0,G_0)+r_n\\
&=&(P_n-P_0)D^*(Q_0,G_0)+(P_n-P_0)\{D(Q_n^*,G_n)-D^*(Q_0,G_0)\}\nonumber \\
&&+R_{20}(Q_n^*,G_n,Q_0,G_0)+r_n.\label{exactexpansiontmle}
\end{eqnarray}
The second term on the right-hand side is $O_P(n^{-1/2-\alpha(d)})$ by  empirical process theory and the continuity condition (\ref{contDstar}) on $D^*$.
Thus, this proves the following asymptotic efficiency theorem.

\begin{theorem}\label{thefftmle}
Consider the statistical model ${\cal M}$, target parameter $\Psi:{\cal M}\rightarrow\openr$ and the model assumptions (\ref{sectionalvarbound}), (\ref{boundingR2}), (\ref{contDstar}). In addition, assume that the HAL-TMLE $Q_n^*$ is such that it solves the efficient influence curve equation  (\ref{efficeqn}) up till $r_n=o_P(n^{-1/2})$; 
it preserves the sectional variation norm in the sense that $\pl Q_n^*\pl_v^*<C \pl Q_n\pl_v^*$ for some $C<\infty$; and
$d_{01}(Q_n^*,Q_0)=O_P(d_{01}(Q_n,Q_0))$. 

Then the HAL-TMLE $\Psi(Q_n^*)$ of $\psi_0$ is asymptotically efficient:
\[
\Psi(Q_n^*)-\Psi(Q_0)=(P_n-P_0)D^*(Q_0,G_0)+O_P(n^{-1/2-\alpha(d)}).\]
\end{theorem}

{\bf Wald type confidence interval:}
A first order asymptotic 0.95-confidence interval is given by $\psi_n^*\pm 1.96 \sigma_n/n^{1/2}$ where $\sigma_n^2=P_n \{D^*(Q_n^*,G_n)\}^2$ is a consistent estimator of $\sigma^2_0=P_0\{D^*(Q_0,G_0)\}^2$.
Clearly, this first order confidence interval ignores the exact remainder $\tilde{R}_{2n}$ in the exact expansion $\Psi(Q_n^*)-\Psi(Q_0)=(P_n-P_0)D^*(Q_0,G_0)+\tilde{R}_{2n}$ as presented in (\ref{exactexpansiontmle}):
\begin{equation}\label{exactremainder}
\tilde{R}_{2n}\equiv R_{20}(Q_n^*,G_n,Q_0,G_0)+(P_n-P_0)\{D^*(Q_n^*,G_n)-D^*(Q_0,G_0)\}+r_n.\end{equation}
 
 The asymptotic efficiency proof  above of the HAL-TMLE$(C^u)$   relies on the HAL-MLEs $(Q_{n,C_1^u},G_{n,C_2^u})$ to converge to the true $(Q_0,G_0)$ at rate faster than $n^{-1/4}$, and that their sectional variation norm is uniformly bounded from above by $C^u=(C_1^u,C_2^u)$. Both of these conditions are still known to hold for the CV-HAL-MLE $(Q_{n,C_{1n}},G_{n,C_{2n}})$ in which the constants $(C_1,C_2)$ are selected with the cross-validation selector $(C_{1n},C_{2n})$ \citep{vanderLaan15}. This follows since the cross-validation selector is asymptotically equivalent with the oracle selector, thereby guaranteeing that $C_n$ will exceed the sectional variation norm of the true $(Q_0,G_0)$ with probability tending to 1.
Therefore, we  have that this CV-HAL-TMLE is also asymptotically efficient. Of course, this CV-HAL-TMLE is more practical and powerful than the HAL-TMLE at an apriori specified $(C_1^u,C_2^u)$  since it adapts the choice of bounds $(C_1,C_2)$ to  the true sectional variation norms $C_0=(C_{10},C_{20})$ for $(Q_0,G_0)$.
 
 \begin{theorem}\label{thefftmlecv}
 Let $C_{10}=\pl Q_0\pl_v^*$, $C_{20}=\pl G_0\pl_v^*$.
 Suppose that $C_1^u$ and $C_2^u$  that define the HAL-MLEs $Q_n=Q_{n,C_1^u}$ and $G_n=G_{n,C_2^u}$ are replaced by data adaptive selectors $C_{1n}$ and $C_{2n}$ for which 
 \begin{equation}
 P(C_{10}\leq C_{1n}<C_1^u,C_{20}\leq C_{2n}<C_2^u)\rightarrow 1,\mbox{ as $n\rightarrow\infty$.}\label{Cn}
 \end{equation}
  Then, under the same assumptions as in Theorem \ref{thefftmle}, the TMLE
 $\Psi(Q_n^*)$,  using $Q_n=Q_{n,C_{1n}}$ and $G_n=G_{n,C_{2n}}$ as initial estimators,  is asymptotically efficient. 
 \end{theorem}
 In general, when the model is defined by global  constraints, then one should use cross-validation to select these constraints, which will only improve the performance of the initial estimators and corresponding TMLE, due to its asymptotic equivalence with the oracle selector. So our model might have more global constraints beyond $(C_1^u,C_2^u)$ and these could then also be selected with cross-validation resulting in a CV-HAL-MLE and corresponding HAL-TMLE (see also our two examples).
 
\subsection{Motivation for finite sample inference}

In order to understand how large the exact remainder $\tilde{R}_{2n}$ could be relative to the leading first order term, we need to understand the size of $d_{01}(Q_n,Q_0)$ and $d_{02}(G_n,G_0)$. This will then motivate us to propose methods that estimate the finite sample distribution of the HAL-TMLE or conservative versions thereof.

To establish this behavior of $d_{01}(Q_n,Q_0)$ we will use the following general integration by parts formula, and a resulting bound.

\begin{lemma}\label{lemmaip}
Let $F,Z\in D_d[0,\tau]$. 
For a given function $Z\in D_d[0,\tau]$ we define
\[
\bar{Z}(u)=Z([u,\tau])=\int_{[u,\tau]} dZ(s),\]
the measure $dZ$ assigns to the cube $[u,\tau]$, which is a generalized difference across the $2^d$-corners of $[u,\tau]$ \citep{Gill&vanderLaan&Wellner95}.
For any two functions $F,Z\in D[0,\tau]$ with $\pl F\pl_v^*<\infty$ and $\pl Z\pl_v^*<\infty$, we have the following integration by parts formula:
\[
\int_{[0,\tau]} F(x)dZ(x)=F(0)\bar{Z}(0) +\sum_s \int_{u_s} \bar{Z}(u_s,0_{-s})dF_s(u_s).\]
This implies 
\[
\int_{[0,\tau]} F(x)dZ(x)\leq \pl \bar{Z}\pl_{\infty}\pl F\pl_v^* .\]
\end{lemma}
{\bf Proof:}
The representation of $F(x)$ is presented in \citep{Gill&vanderLaan&Wellner95,vanderLaan15}. Using this representation yields the presented integration by parts formula as follows:
\begin{eqnarray*}
\int F dZ&=&\int \{F(0)+\sum_s \int_{(0_s,x_s]} dF_s(u)\} dZ(x)\\
&=&F(0) Z([0,\tau])+\sum_s \int_x\int_{u_s} I_{x_s\geq u_s} dF_s(u_s)dZ(x)\\
&=& F(0)Z([0,\tau])+\sum_s \int_{u_s} Z([u_s,\tau_s]\times [0_{-s},\tau_{-s}]) dF_s(u_s)\\
&\leq & \max_s \sup_{u_s\in [0_s,\tau_s]}\mid Z ([u_s,\tau_s]\times [0_{-s},\tau_{-s}])\mid \pl F\pl_v^*\\
&=& \pl Z\pl_{\infty} \pl F\pl_v^*.\Box
\end{eqnarray*}

For a $P\in {\cal M}$ and $P=P_n$ we define
\[
\bar{P}(u)=P([u,\tau])=\int_{[u,\tau]} dP(s).\]

By Lemma \ref{lemmaip}, we have
\begin{eqnarray*}
d_{01}(Q_n,Q_0)&\leq & \pl \bar{P}_n-\bar{P}_0\pl_{\infty}\pl L_1(Q_n,Q_0)\pl_v^* \\
d_{02}(G_n,G_0)&\leq &\pl \bar{P}_n-\bar{P}_0\pl_{\infty}\pl L_2(G_n,G_0)\pl_v^* .
\end{eqnarray*}
By \citep{vanderVaart&Wellner11}, we can bound the expectation of the  supremum norm of an empirical process over a  class of functions with uniformly bounded envelope  by the entropy integral:
\[
E\pl \sqrt{n}(\bar{P}_n-\bar{P}_0)\pl_{\infty}\lesssim J(1,{\cal F}_I)\equiv \sup_Q \int_{0}^1 \sqrt{\log N(\epsilon,L^2(Q),{\cal F}_I)} d\epsilon .\]
The covering number $N(\epsilon,L^2(Q),{\cal F}_I)$ for the class of indicators ${\cal F}_I=\{I_{[u,\tau]}:u\in [0,\tau]\}$ behaves as $\epsilon^{-d}$.   
This proves that
$E\pl \sqrt{n}(\bar{P}_n-\bar{P}_0)\pl_{\infty}=O(d^{1/2})$,
and thus
\begin{eqnarray*}
Ed_{01}(Q_n,Q_0)&=& O(d^{1/2}n^{-1/2}M_1)\\
Ed_{02}(G_n,G_0)&=& O(d^{1/2}n^{-1/2}M_2).
\end{eqnarray*}
In particular, this shows that the exact second-order remainder (\ref{exactremainder}) can be bounded in expectation as follows:
\[
E\mid \tilde{R}_{2n} \mid =O(n^{-1/2}d^{1/2} \sqrt{M_1M_2}).\]
Even though these bounds are overly conservative, these bounds provide a clear indication how the size of $d_{01}(Q_n,Q_0)$ and $d_{02}(G_n,G_0)$, and thereby the second-order remainder is  potentially affected by the dimension $d$ (i.e., for nonparametric models) and the allowed complexity of the model as measured by the bounds $M_1,M_2$.

One can thus conclude that there are many settings in which the exact second-order remainder $\tilde{R}_{2n}$ will dominate the leading linear term $(P_n-P_0)D^*(P_0)$ in finite samples.
Therefore, for the sake of accurate inference we will need methods that estimate the actual finite sample sampling distribution of the HAL-TMLE.

\subsection*{A very conservative finite sample confidence interval}
Consider the case that $r_n=0$.
Let ${\bf M}_1^*=\sup_{P\in {\cal M}}\max_{\epsilon}\pl {\bf L}_1(Q(P)_{\epsilon})\pl_v^*$ and ${\bf M}_2=\sup_{P\in {\cal M}}\pl {\bf L}_2(G(P))\pl_v^*$ be the  deterministic upper bound on the sectional variation norms of  ${\bf L}_1(Q_n^*)$ and ${\bf L}_2(G_n)$. 
 Let $\bar{Z}_n=n^{1/2}(\bar{P}_n-\bar{P}_0)$. The integration by parts bound applied to (\ref{exactexpansiontmle}) yields the following bound:
 \begin{eqnarray*}
\mid n^{1/2}(\Psi(Q_n^*)-\Psi(Q_0))\mid&\leq& \pl D^*(Q_n^*,G_n)\pl_v^*\pl \bar{Z}_n\pl_{\infty}\\
&&\hspace*{-2cm}
+f(\pl {\bf L}_1(Q_n^*,Q_0)\pl_v^{*1/2}\pl \bar{Z}_n\pl_{\infty}^{1/2},\pl {\bf L}_2(G_n,G_0)\pl_v^{*1/2}\pl\bar{Z}_n\pl_{\infty}^{1/2}) .\end{eqnarray*}
Let $M_3^*\equiv \sup_{P\in {\cal M}}\max_{\epsilon}\pl D^*(Q(P)_{\epsilon},G)\pl_v^*$. Then, we obtain the following bound:
\begin{equation}\mid n^{1/2}(\Psi(Q_n^*)-\Psi(Q_0))\mid\leq \left\{ M_3^*+f({\bf M}_1^{*1/2},{\bf M}_2^{1/2})\right\} \pl \bar{Z}_n\pl_{\infty}.\label{consbound}
\end{equation}
 Let $q_{n,0.95}$ be the $0.95$-quantile of $\pl \bar{Z}_n\pl_{\infty}$. A conservative finite sample $0.95$-confidence interval is then given by:
 \[
 \Psi(Q_n^*)\pm C({\bf M}_1^*,{\bf M}_2,M_3^*)q_{n,0..95}/n^{1/2},\]
where $C({\bf M}_1,{\bf M}_2,M_3)= M_3+f({\bf M}_1^{1/2},{\bf M}_2^{1/2})$.
One could estimate the distribution of $\bar{Z}_n$ with the nonparametric bootstrap and thereby obtain an bootstrap-estimate $q_{n,0.95}^{\#}$ of $q_{n,0.95}$.
One could push the conservative nature of this confidence interval further by using theoretical bounds for the tail-probability $P(\pl \bar{Z}_n\pl_{\infty}>x)$
and define the quantile $q_{n,0.95}$ in terms of this theoretical upper bound (such exponential bounds are available in (e.g.) \citep{vanderVaart&Wellner96}, but the constants in these exponential bounds appear to not  be concretely specified).

 The bound (\ref{consbound}) simplifies if we focus on the sampling distribution of the one-step estimator $\psi_n^1=\Psi(Q_n)+P_n D^*(Q_n,G_n)$ by being able to replace the targeted version $Q_n^*$ by $Q_n$.
 For the one-step estimator we have
 \[
 \psi_n^1-\psi_0=(P_n-P_0)D^*(Q_n,G_n)+R_{20}(Q_n,G_n,Q_0,G_0) .\]
 Let ${\bf M_1}=\sup_{P\in {\cal M}}\pl {\bf L}_1(Q(P))\pl_v^*$ and ${\bf M}_2=\sup_{P\in {\cal M}}\pl {\bf L}_2(G(P))\pl_v^*$ be the upper bound on the sectional variation norms of ${\bf L}_1(Q_n)$ and ${\bf L}_2(G_n)$.
Analogue to above, we obtain
 \begin{eqnarray*}
\mid n^{1/2}(\psi_n^1-\Psi(Q_0))\mid&\leq& \pl D^*(Q_n,G_n)\pl_v^*\pl \bar{Z}_n\pl_{\infty}\\
&&\hspace*{-2cm}
+f(\pl {\bf L}_1(Q_n,Q_0)\pl_v^{*1/2}\pl \bar{Z}_n\pl_{\infty}^{1/2},\pl {\bf L}_2(G_n,G_0)\pl_v^{*1/2}\pl\bar{Z}_n\pl_{\infty}^{1/2}) .\end{eqnarray*}
Recall $M_3\equiv \sup_{P\in {\cal M}}\pl D^*(P)\pl_v^*$. Then, we obtain the following conservative sampling distribution:
\begin{equation}
Z_n^+\equiv \mid n^{1/2}(\psi_n^1-\Psi(Q_0))\mid\leq \left\{ M_3+f({\bf M}_1^{1/2},{\bf M}_2^{1/2})\right\} \pl \bar{Z}_n\pl_{\infty},\label{consbound1}
\end{equation}
and  conservative finite sample $0.95$-confidence interval 
 \[
 \psi_n^1\pm C({\bf M}_1,{\bf M}_2,M_3)q_{n,0..95}/n^{1/2},\]
where $C({\bf M}_1,{\bf M}_2,M_3)= M_3+f({\bf M}_1^{1/2},{\bf M}_2^{1/2})$. Clearly, this same confidence interval can be applied to the TMLE  since the TMLE is asymptotically equivalent with the one-step estimator and generally performs better in finite samples by being a substitution estimator. 
 
 Above, we pointed out that  $\pl \bar{Z}_n\pl_{\infty}=O_P((n/d)^{-1/2})$, which shows that this confidence interval has a width of order $(n/d)^{-1/2}$.
 This confidence interval is not only finite sample conservative but is also not  asymptotically sharp. 
Nonetheless, this formula appears to demonstrate  that the dimension $d$ of the data $O$ enters directly into the rate of convergence as $(n/d)^{1/2}$. In addition, it shows  and that the actual bounds $(C_1^u,C_2^u)$ on the sectional  variation norms of $Q$ and $G$ are directly affecting the width of the confidence interval (essentially linearly).
 In addition, the dimension $d$ itself naturally affects the chosen upper bounds $(C_1^u,C_2^u)$ and thereby ${\bf M}_1,{\bf M}_2,M_3$, so that the dimension $d$ may also affect the width of the finite sample confidence interval through the constant $C({\bf M}_1,{\bf M}_2,M_3)$. 
 
 Though interesting,  we suggest that in most applications this bound is much too conservative for practical use. This motivates us to construct much more accurate estimators of the actual sampling distribution of the HAL-TMLE.
 
 The above finite sample bound could also be applied to choices $M_{1n},M_{2n},M_{3n}$ implied by the cross-validation selector $(C_{1n},C_{2n})$ of $(C_1,C_2)$. We suggest that this would make the resulting confidence interval more reasonable, by not being so conservative (by being forced to select conservative upper bounds $C^u=(C_1^u,C_2^u)$). 



\section{The nonparametric bootstrap for the HAL-TMLE}\label{sectnp}
Let $O_1^{\#},\ldots,O_n^{\#}$ be $n$ i.i.d. draws from the empirical measure $P_n$.
Let $P_n^{\#}$ be the empirical measure of this bootstrap sample. In the following we define a generalized definition of $Q$ being absolutely continuous w.r.t. $Q_n$:
$Q\ll Q_n$.
\begin{definition}\label{defabscont}
Recall the representation (\ref{Frepresentation}) for a mulivariate real valued cadlag function $F$ in terms of its sections $F_s$.
We will say that $Q_k$ is absolutely continuous w.r.t. $Q_{k,n}$ if for each subset $s\subset\{1,\ldots,m_{1k}\}$, its $s$-specific section $Q_{k,s}$ defined by
$u_s\rightarrow Q_{k}(u_s,0_{-s})$ is absolutely continuous w.r.t. $Q_{n,k,s}$ defined by $u_s\rightarrow Q_{n,k}(u_s,0_{-s})$.
 We use the notation $Q_k\ll Q_{n,k}$. In addition, we use the notation $Q\ll Q_n$ if $Q_k\ll Q_{n,k}$ for each component $k\in \{1,\ldots,K_1\}$.
 Similarly, we use this notation $G\ll G_n$ if $G_k\ll G_{n,k}$ for each component $k\in \{1,\ldots,K_2\}$.
 \end{definition}
 In practice, the HAL-MLE $Q_n=\arg\min_{Q\in Q({\cal M})}P_n {\bf L}_1(Q)$ is attained by a discrete measure  $Q_n$ so that it can be computed by  minimizing the empirical risk over a large linear combination of indicator basis functions (e.g., $2^{m_{1k}} n$  for $Q_{nk}$) under the constraint that the sum of the absolute value of the coefficients is bounded by the specified constant $C_1$ \citep{Benkeser&vanderLaan16}.
 However, $Q_n$ will only have around $n$ non-zero coefficients. In that case, the constraint $Q\ll Q_n$ states that $Q$ is a linear combination of the indicator basis functions that had a non-zero coefficient in $Q_n$.

 Let  $Q_n^{\#}=\arg\min_{Q\in Q({\cal M}),Q\ll Q_n}P_n^{\#}L_1(Q)$ and
$G_n^{\#}=\arg\min_{G\in G({\cal M}),G\ll G_n}P_n^{\#}L_2(G)$ be the corresponding HAL-MLEs of $Q_n=\arg\min_{Q\in Q({\cal M})}P_nL_1(Q)$ and $G_n=\arg\min_{G\in G({\cal M})}P_nL_2(G)$ based on these bootstrap samples. 
Since the empirical measure $P_n^{\#}$ has a support contained in $P_n$, we expect that in many problems $Q_n^{\#}=\arg\min_{Q\in Q({\cal M})}P_n^{\#}L_1(Q)$ satisfies $Q_n^{\#}\ll Q_n$ and similarly for $G_n^{\#}$. Either way, the extra restriction $Q\ll Q_n$ makes the computation of the HAL-MLE on the bootstrap sample much faster than the HAL-MLE $Q_n$ based on the original sample, so that enforcing this extra constraint is only beneficial from a computational point of view.
That is, the computation of $Q_n^{\#}$ only involves minimizing the empirical risk w.r.t. $P_n^{\#}$ over maximally $n$ non-zero coefficients, making the calculation of $Q_n^{\#}$ relatively trivial. 

Let $\epsilon_n^{\#}=\arg\min_{\epsilon}P_n^{\#}L_1(Q_{n,\epsilon}^{\#})$ be the one-step TMLE update of $Q_n^{\#}$ based on the least favorable submodel $\{Q_{n,\epsilon}^{\#}:\epsilon\}$ through $Q_n^{\#}$ at $\epsilon =0$ with score $D^*(Q_n^{\#},G_n^{\#})$ at $\epsilon =0$.  Let $Q_n^{\#*}=Q_{n,\epsilon_n^{\#}}^{\#}$ be the TMLE update which is assumed to solve $r_n^{\#}\equiv \mid P_n^{\#}(Q_n^{\#*},G_n^{\#})\mid =o_{P_n}(n^{-1/2})$, conditional on $(P_n:n\geq 1)$ (just like $r_n=o_P(n^{-1/2})$). Finally, let $\Psi(Q_n^{\#*})$ be the TMLE of $\Psi(Q_n^*)$ based on this nonparametric bootstrap sample. 
We estimate the finite sample distribution of $n^{1/2}(\Psi(Q_n^*)-\Psi(Q_0))$ with the sampling distribution of $Z_n^{1,\#}\equiv n^{1/2}(\Psi(Q_n^{\#*})-\Psi(Q_n^*))$, conditional on $P_n$.
Let $\Phi_n^{\#}(x)=P(n^{1/2}(\Psi(Q_n^{\#*})-\Psi(Q_n^*))\leq x\mid P_n)$ be the cumulative distribution of this bootstrap sampling distribution.
So a bootstrap based 0.95-confidence interval for $\psi_0$ is given by \[
[\psi_n^{*}+q_{0.025,n}^{\#}/n^{1/2},\psi_n^*+q_{0.975,n}^{\#}/n^{1/2} ],\]
 where 
$q_{\alpha,n}^{\#}=\Phi_n^{\#-1}(\alpha)$ is the $\alpha$-quantile of this bootstrap distribution.

One could also apply this  nonparametric bootstrap to $n^{1/2}\mid \Psi(Q_n^*)-\Psi(Q_0)\mid/\sigma_n$, where $\sigma_n^2$ is an estimator of the variance of $D^*(Q_0,G_0)$. It is not clear if this has any advantage, beyond that the confidence interval is now of the form $[\psi_n^{*}+q_{0.025,n}^{\#}\sigma_n/n^{1/2},\psi_n^*+q_{0.975,n}^{\#}\sigma_n/n^{1/2} ]$, where $q_{\alpha,n}^{\#}$ is the $\alpha$-quantile of the cumulative distribution function of $n^{1/2}(\Psi(Q_n^{\#*})-\Psi(Q_n^*))/\sigma_n^{\#}$, conditional on $P_n$, imitating the Wald-type confidence interval.

We now want to prove that $\Phi_n^{\#}$ converges to the cumulative distribution function of limit distribution $N(0,\sigma^2_0)$ so that we are consistently estimating the limit distribution of the TMLE. Importantly, this nonparametric bootstrap confidence interval could potentially dramatically improve the coverage relative to using the first order Wald-type confidence interval since this bootstrap distribution is estimating the variability of the full-expansion of the TMLE, including the exact remainder $\tilde{R}_{2n}$.

In the next subsection we show that the nonparametric bootstrap works for the HAL-MLEs $Q_n$ and $G_n$. Subsequently, not surprisingly,  we can show that this also establishes that the bootstrap works for the one-step TMLE $Q_n^*$ ($K$-th step TMLE for fixed $K$). This provides then the basis for proving that the nonparametric bootstrap is consistent for the HAL-TMLE.

\subsection{Nonparametric bootstrap for HAL-MLE}

The following theorem establishes that the bootstrap HAL-MLE $Q_n^{\#}$  estimates $Q_n$ as well w.r.t. an empirical loss-based dissimilarity $d_{n1}(Q_n^{\#},Q_n)=P_n L_1(Q_n^{\#})-P_nL_1(Q_n)$  as
$Q_n$  estimates  $Q_0$  with respect to $d_{01}(Q_n,Q_0)=P_0L_1(Q_n,Q_0)$. Moreover, it proves that $d_{n1}(Q_n^{\#},Q_n)$ is at minimal equivalent with a  square of an $L^2(P_n)$-norm defined by the exact second-order remainder in a first order Tailor expansion of $P_n L_1(Q)$ at $Q_n$.
The analogue results apply to $G_n^{\#}$. We are stating the theorem for the sum loss function $L_1$, but it can also  be applied to each separate HAL-MLE of $Q_{0,k_1}$ with its loss $L_{1k_1}(Q_{k_1})$ and $G_{0,k_2}$ with its loss $L_{2k_2}(G_{k_2})$ to provide a separate result for each HAL-MLE. In fact, we could simply replace $L_1$ by ${\bf L}_1$ and $L_2$ by ${\bf L}_2$ to obtain the theorem for all components of $Q$.

{\bf Sectional variation norm of HAL-MLE dominates sectional variation norm of bootstrapped HAL-MLE:}
We either assume that $\pl Q_n\pl_v^*=C_1^u$  achieves the maximal allowed value $C_1^u$ or the weaker assumption $\pl Q_n^{\#}\pl_v^*\leq \pl Q_n\pl_v^*$, conditional on $P_n$.
Of course, for the sake of asymptotics we would only need this to hold with probability tending to 1.
The same assumption is used for $G_n$ and $G_n^{\#}$. If $C_1^u$ is chosen so that the sectional variation norm of an MLE $Q_n$ is smaller than $C_1^u$ even though it is a perfect fit of the data in the sense that $P_n L_1(Q_n)=0$ (i.e., smallest possible value), then $\pl Q_n\pl_v^*=C_1^u$ would not be satisfied. Therefore $\pl Q_n\pl_v^*=C_1^u$ requires to make sure that $C_1^u$ is selected small enough relative to sample size so that  the MLE is not a complete overfit of the data. If $C_1^u$ is replaced by a the cross-validation selector $C_{1n}$, our experience is that the HAL-MLE (i.e.,, the Lasso) achieves it maximal allowed value for the sum of the absolute value of its coefficients: i.e,  $\pl Q_n\pl_v^*=C_{1n}$. In fact, all we need is that $\pl Q_n^{\#}\pl_v^*\leq \pl Q_n\pl_v^*$, which is a weaker assumption and could easily be true for all  choices of $C_1^u$: for example, Lasso regression applied to bootstrap sample (i.e., subset of original data but using weights)   might select an $L_1$-norm of its coefficient vector smaller than the $L_1$-norm when applied to the original sample, whatever $C_1^u$  is selected.

\begin{theorem}\label{thnpbootmle}
Recall our assumption (\ref{calFmodel}) or (\ref{calFmodelplus}) on the parameter spaces of $Q$ and $G$.\nl
{\bf Definitions:}
Let $d_{n1}(Q,Q_n)=P_n \{L_1(Q)-L_1(Q_n)\}$ be the loss-based dissimilarity at the empirical measure, where $Q_n=\arg\min_{Q\in Q({\cal M})}P_n L_1(Q)$.
Similarly, let $d_{n2}(G,G_n)=P_n \{L_2(G)-L_2(G_n)\}$ be the loss-based dissimilarity at the empirical measure, where $G_n=\arg\min_{G\in G({\cal M})}P_n L_2(G)$.
Let $P_n R_{2L_1,n}(Q_n^{\#},Q_n)$ be defined as the exact second-order remainder of a first order Tailor expansion of $P_nL_1(Q)$ at $Q_n$:
\[ P_n \{L_1(Q_n^{\#})-L_1(Q_n)\}=P_n \frac{d}{dQ_n}L_1(Q_n)(Q_n^{\#}-Q_n)+P_n R_{2L_1,n}(Q_n^{\#},Q_n),\]
where $\frac{d}{dQ_n}L_1(Q_n)(h)=\left . \frac{d}{d\epsilon} L_1(Q_n+\epsilon h)\right |_{\epsilon =0}$ is the directional derivative in direction $h$.
Similarly,  we define $P_0R_{2L_1,0}(Q_n,Q_0)$ as the exact second-order remainder of a first order Tailor expansion of $P_0L_1(Q)$ at $Q_0$:
\[
P_0\{L_1(Q_n)-L_1(Q_0)\}=P_0\frac{d}{dQ_0}L_1(Q_0)(Q_n-Q_0)+P_0R_{2L_1,0}(Q_n,Q_0).\]
Similarly, we define $P_0R_{2L_2,0}(G_n,G_0)$ and $P_n R_{2L_2,n}(G_n^{\#},G_n)$.
\newline
{\bf Assumption:}
Assume $\pl Q_n\pl_v^*=C_1^u$, $\pl G_n\pl_v^*=C_2^u$ (i.e., they attain the maximal allowed value) or assume that $\pl Q_n^{\#}\pl_v^*\leq \pl Q_n\pl_v^*$  and $\pl G_n^{\#}\pl_v^*\leq \pl G_n\pl_v^*$ with probability 1, conditional on $(P_n:n\geq 1)$.
Suppose that 
\begin{eqnarray}
P_0\{L_1(Q_n)-L_1(Q_0)\}^2&\lesssim& P_0 R_{2L_1,0}(Q_n,Q_0)\nonumber \\
P_n \{L_1(Q_n^{\#})-L_1(Q_n)\}^2&\lesssim& P_nR_{2L_1,n}(Q_n^{\#},Q_n)\nonumber\\
P_0\{L_2(G_n)-L_2(G_0)\}^2&\lesssim& P_0 R_{2L_2,0}(G_n,G_0)\nonumber \\
P_n \{L_2(G_n^{\#})-L_2(G_n)\}^2&\lesssim& P_n R_{2L_2,n}(G_n^{\#},G_n).\label{boundb}
\end{eqnarray}
{\bf Conclusion:}
Then,
\[
d_{n1}(Q_n^{\#},Q_n)=O_P(n^{-1/2-\alpha(d)})\mbox{ and } d_{n2}(G_n^{\#},G_n)=O_P(n^{-1/2-\alpha(d)}).\]
In addition, we have  $P_n\frac{d}{dQ_n}L_1(Q_n)(Q_n^{\#}-Q_n)\geq 0$ so that  $d_{n1}(Q_n^{\#},Q_n)$ is more powerful dissimilarity than the quadratic $P_n R_{2L_1,n}(Q_n^{\#},Q_n)$ (i.e., convergence w.r.t. $d_{n1}$ implies convergence w.r.t. latter):
\[
P_n\{L_1(Q_n^{\#})-L_1(Q_n)\}\geq P_n R_{2L_1,n}(Q_n^{\#},Q_n).\]
Similarly,  we have $P_0\frac{d}{dQ_0}L_1(Q_0)(Q_n-Q_0)\geq 0$ so that $d_{01}(Q_n,Q_0)$ dominates $P_0R_{2L_1,0}(Q_n,Q_0)$:
\[
P_0\{L_1(Q_n)-L_1(Q_0)\}\geq P_0 R_{2L_1,0}(Q_n,Q_0).\]
As a consequence, we also have 
\[
P_n R_{2L_1,n}(Q_n^{\#},Q_n)=O_P(n^{-1/2-\alpha(d)})\mbox{ and } P_n R_{2L_2,n}(G_n^{\#},G_n)=O_P(n^{-1/2-\alpha(d)}).\]
{\bf Bootstrapping HAL-MLE$(C)$ at $C=C_n$:}
This theorem also applies to the case that  $C^u=(C_1^u,C_2^u)$  is replaced by a data adaptive choice $C_n=(C_{1n},C_{2n})$ (i.e., depending on $P_n$) satisfying (\ref{Cn}).
 \end{theorem}
 Note that if $C^u=C_n$, then conditional on $P_n$, $C_n$ is still fixed, so that establishing the latter result only requires checking that the convergence of the bootstrapped HAL-MLE $(Q_{n,C_1}^{\#},G_{n,C_2}^{\#})$ to the HAL-MLE $(Q_{n,C_1},G_{n,C_2})$ at  a fixed $C$ w.r.t. the loss-based dissimilarities $d_{n1}$ and $d_{n2}$ holds uniformly in $C$ between 
 the true sectional variation norms $C_0$ and the model upper bound $C^u$. The validity of this theorem does not rely on  $C_n$  exceeding $C_0$, but the latter is needed for establishing that the HAL-MLE $Q_{n,C_n}$ is consistent for $Q_0$ and thus the efficiency of the HAL-TMLE $\Psi(Q_n^*)$. 
The proof of Theorem \ref{thnpbootmle} is presented in the Appendix \ref{AppendixB}.

Clearly, $P_n R_{2L_1,n}(Q_n^{\#},Q_n)$ will behave as a square of a difference of $Q_n^{\#}$ and $Q_n$.
In our proof below of the validity of the nonparametric bootstrap method for the HAL-TMLE we will need that convergence of $d_{n1}(Q_n^{\#},Q_n)$ and $d_{01}(Q_n,Q_0)$ implies convergence of $d_{01}(Q_n^{\#},Q_0)$ as well. This requires showing that convergence w.r.t. an $L^2(P_n)$-norm implies convergence at the same rate w.r.t. $L^2(P_0)$-norm. For that purpose we note the following lemma, which is also proved in the Appendix \ref{AppendixB}.

\begin{lemma}\label{lemmadndo}
If $P_n f_n^2=O_P(n^{-1/2-\alpha(d)})$ for some $f_n$ with $\pl f_n\pl_v^*<M$ for some $M<\infty$ with probability 1, then we also have
$P_0 f_n^2=O_P(n^{-1/2-\alpha(d)})$. 
\end{lemma}

\subsection{Preservation of rate of convergence for the targeted bootstrap estimator}
It is no surprise that under a weak regularity condition, we have that $\epsilon_n^{\#}=O_P(n^{-1/4-\alpha(d)/2})$ converges at same rate as $Q_n^{\#}$ to $Q_n$.  As a result, the TMLE-update $Q_{n,\epsilon_n^{\#}}^{\#*}$ of $Q_n^{\#}$ converges at  the same rate to $Q_n$ as $Q_n^{\#}$.
A general proof of this result is presented in the  Appendix \ref{AppendixC} under weak regularity conditions on the least favorable submodel. 


\subsection{The nonparametric bootstrap for the  HAL-TMLE}
We can now imitate the efficiency proof for the HAL-TMLE to obtain the desired result for the bootstrapped HAL-TMLE of $\Psi(Q_n^*)$.
In addition to the model assumptions of Theorem \ref{thefftmle} for asymptotic efficiency of the TMLE, we asssume the conditions (\ref{boundb}) of Theorem \ref{thnpbootmle} for validity of the nonparametric bootstrap for the HAL-MLE.
In addition, we assume the very weak condition that convergence of $d_{n1}(Q_n^{\#},Q_n)$ and $d_{01}(Q_n,Q_0$ implies the same convergence of $d_{01}(Q_n^{\#},Q_0)$:
\begin{equation}\label{dn1bounding}
\max(d_{n1}(Q_n^{\#},Q_n),d_{01}(Q_n,Q_0))=O_P(r(n)) \mbox{ implies $d_{01}(Q_n^{\#},Q_0)=O_P(r(n))$,}\end{equation}
 and similarly $\max(d_{n2}(G_n^{\#},G_n),d_{02}(G_n,G_0))=O_P(r(n))$ implies $d_{02}(G_n^{\#},G_0)=O_P(r(n))$.
 To verify this assumption, one  can use that 
$P_n R_{2L_1,n}(Q_n^{\#},Q_n)\leq d_{n1}(Q_n^{\#},Q_n)$, $P_n R_{2L_2,n}(G_n^{\#},G_n)\leq d_{n2}(G_n^{\#},G_n)$,
$P_0R_{2L_1,0}(Q_n,Q_0)\leq d_{01}(Q_n,Q_0)$, and $P_0R_{2L_2,0}(G_n,G_0)\leq d_{02}(G_n,G_0)$, and Lemma \ref{lemmadndo} to translate
$\int f_n^2dP_n=o_P(r(n))$ implies $\int f_n^2 dP_0=O_P(r(n))$.

Finally, we  assume a the empirical analogue of the uniform continuity condition (\ref{contDstar}) on the efficient influence curve:
\begin{equation}\label{contDstaremp}
P_n\{D^*(Q_n^{\#},G_n^{\#})-D^*(Q_n,G_n)\}^2\lesssim d_{n1}(Q_n^{\#},Q_n)+d_{n2}(G_n^{\#},G_n).\end{equation}
Again, to verify this we can use $P_n R_{2L_1,n}(Q_n^{\#},Q_n)\leq d_{n1}(Q_n^{\#},Q_n)$ and $P_n R_{2L_2,n}(G_n^{\#},G_n)\leq d_{n2}(G_n^{\#},G_n)$,
so that it suffices to verify
\[
P_n\{D^*(Q_n^{\#},G_n^{\#})-D^*(Q_n,G_n)\}^2\lesssim P_n R_{2L_1,n}(Q_n^{\#},Q_n)+P_n R_{2L_2,n}(G_n^{\#},G_n).\]



\begin{theorem}\label{thnpboothaltmle}\  \nl
{\bf Assumptions:}
Assume the conditions of Theorem \ref{thefftmle} providing asymptotic efficiency of $\Psi(Q_n^*)$; $\pl Q_n\pl_v^*=C_1^u$, $\pl G_n\pl_v^*=C_2^u$ (i.e., they attain the maximal allowed value) or that $\pl Q_n^{\#}\pl_v^*\leq \pl Q_n\pl_v^*$  and $\pl G_n^{\#}\pl_v^*\leq \pl G_n\pl_v^*$ with probability 1, conditional on $(P_n:n\geq 1)$;
(\ref{boundb}); (\ref{dn1bounding});  (\ref{contDstaremp}) on  loss functions $L_1(Q)$ and $L_2(G)$; $r_n^{\#}=P_n^{\#}D^*(Q_n^{\#*},G_n^{\#})=o_P(n^{-1/2})$, conditional on $(P_n:n\geq 1)$;  and that $Q_n^{\#*}$ preserves rate of convergence of $Q_n^{\#}$ in the sense that the following three conditions hold: 1) $\pl Q_n^{\#*}\pl_v^*<C\pl Q_n^{\#}\pl_v^*$ for some $C<\infty$; 2) $P_n \{D^*(Q_n^{\#*})-D^*(Q_n^{\#})\}^2\rightarrow_p 0$, conditional on $(P_n:n\geq 1)$; 3) $P_0\{D^*(Q_n^{\#*})-D^*(Q_n^{\#})\}^2\rightarrow_p 0$. 
If we use the one-step TMLE $Q_n^*=Q_{n,\epsilon_n}$, then the last three conditions can be replaced by $\epsilon_n^{\#}=O_P(d_{n1}^{1/2}(Q_n^{\#},Q_n))$.

{\bf Conclusion:}
Then, $d_{n1}(Q_n^{\#},Q_n)=O_P(n^{-1/2-\alpha(d)})$, $d_{n2}(G_n^{\#},G_n)=O_P(n^{-1/2-\alpha(d)})$, 
 $P_n R_{2L_1,n}(Q_n^{\#*},Q_n)=O_P(n^{-1/2-\alpha(d)})$ and $P_n R_{2L_2,n}(G_n^{\#},G_n)=O_P(n^{-1/2-\alpha(d)})$.
 
In addition, 
\[
\Psi(Q_n^{\#*})-\Psi(Q_n)=(P_n^{\#}-P_n)D^*(Q_n,G_n)+O_P(n^{-1/2-\alpha(d)}),\]
and thus $Z_n^{1,\#}\equiv n^{1/2}(\Psi(Q_n^{\#*})-\Psi(Q_n^*))\Rightarrow_d N(0,\sigma^2_0)$, conditional on $(P_n:n\geq 1)$.
\newline
{\bf Consistency of the nonparametric bootstrap for HAL-TMLE at cross-validation selector $C_n$:}
This theorem can be applied to $C^u=C_n$ satisfying (\ref{Cn}).
\end{theorem}
{\bf Proof:}
We provide the proof for the one-step TMLE using the condition $\epsilon_n^{\#}=O_P(d_{n1}^{1/2}(Q_n^{\#},Q_n)$. The proof for the general TMLE $Q_n^*$ using conditions 1-3 instead follows immediately from the following proof as well and below we point out how the proof is generalized to this general case.
Firstly, by definition of the remainder $R_{20}()$ we have the following two expansions:
\begin{eqnarray*}
\Psi(Q_n^{\#*})-\Psi(Q_0)&=&(P_n^{\#}-P_0) D^*(Q_n^{\#*},G_n^{\#})+R_{20}(Q_n^{\#*},G_n^{\#},Q_0,G_0)\\
&=& (P_n^{\#}-P_n)D^*(Q_n^{\#*},G_n^{\#})+(P_n-P_0)D^*(Q_n^{\#*},G_n^{\#})\\
&& +R_{20}(Q_n^{\#*},G_n^{\#},Q_0,G_0)\\
\Psi(Q_n^*)-\Psi(Q_0)&=&(P_n-P_0) D^*(Q_n^*,G_n)+R_{20}(Q_n^*,G_n,Q_0,G_0),
\end{eqnarray*}
where we ignored  $r_n=P_nD^*(Q_n^*,G_n)$ and its bootstrap analogue $r_n^{\#}=P_n^{\#}D^*(Q_n^{\#*},G_n^{\#})$ in these two expressions (which were both assumed to be $o_P(n^{-1/2})$).
Subtracting the first equality from the second equality yields:
\begin{eqnarray}
\Psi(Q_n^{\#*})-\Psi(Q_n^*)&=&(P_n^{\#}-P_n)D^*(Q_n^{\#*},G_n^{\#})+(P_n-P_0)\{D^*(Q_n^{\#*},G_n^{\#})-D^*(Q_n^*,G_n)\}\nonumber \\
&&+
R_{20}(Q_n^{\#*},G_n^{\#},Q_0,G_0)-R_{20}(Q_n^*,G_n,Q_0,G_0). \label{helpa}
\end{eqnarray}
Under the conditions of Theorem \ref{thefftmle}, we already established that $R_{20}(Q_n^*,G_n,Q_0,G_0)=O_P(n^{-1/2-\alpha(d)})$.
By  assumption (\ref{boundingR2}), we can bound the first remainder $R_{20}(Q_n^{\#*},G_n^{\#},Q_0,G_0)$ by $f({\bf d}_{01}^{1/2}(Q_n^{\#*},Q_0),{\bf d}_{02}^{1/2}(G_n^{\#},G_0))$.
Theorem \ref{thnpbootmle} established  that  $d_{n1}(Q_n^{\#},Q_n)=O_P(n^{-1/2-\alpha(d)})$ and $d_{n2}(G_n^{\#},G_n)=O_P(n^{-1/2-\alpha(d)})$.
By assumption (\ref{dn1bounding}), this implies also that $d_{01}(Q_n^{\#},Q_0)$ and $d_{02}(G_n^{\#},G_0)$ are $O_P(n^{-1/2-\alpha(d)})$.
Again, by assumption (\ref{boundingR2}) this yields that $R_{20}(Q_n^{\#},G_n^{\#},Q_0,G_0)=O_P(n^{-1/2-\alpha(d)})$.
By assumption, $\epsilon_n^{\#}=O_P(n^{-1/4-\alpha(d)/2})$, using the fact that $f$ is a quadratic polyonomial, this now also establishes that $R_{20}(Q_n^{\#*},G_n^{\#},Q_0,G_0)=O_P(n^{-1/2-\alpha(d)})$. Similarly, if we work with a general TMLE $Q_n^*$ satisfying conditions 1-3, then this result follows as well.

It remains to analyze the two leading empirical process terms in (\ref{helpa}). 
Firstly, replace $Q_n^{\#*}$ and $Q_n^*$ by $Q_n^{\#}$ and $Q_n$, respectively, in these two terms. This generates three additional remainder terms:
\[
\begin{array}{l}
(P_n^{\#}-P_n)\{D^*(Q_n^{\#*},G_n^{\#})-D^*(Q_n^{\#},G_n^{\#})\}\\
(P_n-P_0)\{D^*(Q_n^{\#*},G_n^{\#})-D^*(Q_n^{\#},G_n)\}\\
(P_n-P_0)\{D^*(Q_n^*,G_n)-D^*(Q_n,G_n)\}.
\end{array}
\]
Since  $Q_n^{\#*}=Q_{n,\epsilon_n^{\#}}^{\#}$ and $Q_n^*=Q_{n,\epsilon_n}$,  each of these terms can be written as $f_n(\epsilon_n^{\#})-f_n(0)$ or $f_n(\epsilon_n)-f_n(0)$ for certain specified $f_n$. We  can carry out an exact first order tailor expansion of this $f_n(\epsilon)$ at $\epsilon =0$ to represents  these three terms as $\epsilon_n^{\#}(P_n^{\#}-P_n)f_{1n} $, $\epsilon_n^{\#}(P_n-P_0)f_{2n}$
and $\epsilon_n (P_n-P_0)f_{3n}$, respectively, for certain functions $f_{1n},f_{2n},f_{3n}$. By assumption (\ref{sectionalvarbound}) these functions $f_{1n},f_{2n},f_{3n}$ have a uniformly bounded sectional variation norm. Thus $(P_n-P_0)f_{jn}=O_P(n^{-1/2})$ for $j=1,2,3$, so that  these three terms can be bounded by  $O_P(n^{-1/2})$ times $\max(\epsilon_n,\epsilon_n^{\#})$, which is $O_P(n^{-3/4-\alpha(d)/2})$. The above three terms are also $o_P(n^{-1/2})$ if we work with a general TMLE $Q_n^*$ and the three conditions 1-3 on $Q_n^{\#*}$ apply.

The remainder of the proof now only involves non-targeted $Q_n^{\#}$ and $Q_n$ so that it generally applies to a general TMLE $Q_n^*$.
Let's now return to the two leading terms in (\ref{helpa}) but with $Q_n^{\#*}$ and $Q_n^*$ replaced by  $Q_n^{\#}$ and $Q_n$, respectively.
By our continuity assumption (\ref{contDstaremp}) on the efficient influence curve as function in $(Q,G)$, we have
that  convergence of $d_{n1}(Q_n^{\#},Q_n)+d_{n2}(G_n^{\#},G_n)$ to zero implies convergence of the square of the $L^2(P_n)$-norm of $D^*(Q_n^{\#},G_n^{\#})-D^*(Q_n,G_n)$ at the same rate. 
By empirical process theory \citep{vanderVaart&Wellner11}, this teaches us that $(P_n^{\#}-P_n)D^*(Q_n^{\#},G_n^{\#})=(P_n^{\#}-P_n)D^*(Q_n,G_n)+O_P(n^{-1/2-\alpha(d)}$. This deals with the first leading term in (\ref{helpa}). 

By our continuity condition (\ref{contDstar}) we also have that 
$P_0\{D^*(Q_n^{\#},G_n^{\#})-D^*(Q_n,G_n)\}^2\rightarrow_p 0$ at  this rate. Again, by \citep{vanderVaart&Wellner11} this shows 
$(P_n-P_0)\{D^*(Q_n^{\#},G_n^{\#})-D^*(Q_n,G_n)\}=O_P(n^{-1/2-\alpha(d)})$.
Thus we have shown  that
\[
\begin{array}{l}
(P_n^{\#}-P_n)D^*(Q_n^{\#},G_n^{\#})+(P_n-P_0)\{D^*(Q_n^{\#},G_n^{\#})-D^*(Q_n,G_n)\}\\
=(P_n^{\#}-P_n)D^*(Q_n^{\#},G_n^{\#})
+O_P(n^{-1/2-\alpha(d)}).\end{array}
\]
The latter term can be written as $(P_n^{\#}-P_n)D^*(Q_n,G_n)+(P_n^{\#}-P_n)\{D^*(Q_n^{\#},G_n)-D^*(Q_n,G_n)\}$.
The second term can be analyzed with empirical process theory as above using (\ref{contDstaremp}) to establish that it is $O_P(n^{-1/2-\alpha(d)})$.

Thus, we have now shown
\[
n^{1/2}(\Psi(Q_n^{\#*})-\Psi(Q_n^*))=n^{1/2}(P_n^{\#}-P_n)D^*(Q_n^*,G_n)+o_P(1)\Rightarrow_d N(0,\sigma^2_0).
\]
This completes the proof of the Theorem for the HAL-TMLE at the fixed $C^u$.  As remarked earlier, it follows straightforwardly that this proof applies uniformly to any $C$ in between $C_0$ and $C^u$, and thereby to a selector $C_n$ satisfying (\ref{Cn}).  $\Box$

\paragraph{Remark regarding robustness in underlying data distribution}
Consider the exact second order remainder $R_{20}$  for the HAL-TMLE $\Psi(Q_n^*)-\Psi(Q_0)=(P_n-P_0)D^*(Q_0,G_0)+R_{20}(Q_n^*,G_n,Q_0,G_0)$
and the exact second order remainder $R_{2n}^{\#}$  for the bootstrapped HAL-TMLE $\Psi(Q_n^{\#*})-\Psi(Q_n^*)=(P_n^{\#}-P_n)D^*(Q_n,G_n)+R_{2n}^{\#}$, as specified in the above proof. Under our model assumptions, the bounding of these two remainders only concern empirical processes $(P_n-P_0)$ and $(P_n^{\#}-P_n)$ indexed by a uniform Donsker class. As shown in \citep{vanderVaart&Wellner96}, such empirical processes converge and satisfy exact finite sample bounds that apply uniformly in all possible data distributions. 
Therefore, it follows that we can also establish that the nonparametric bootstrap is consistent for the normal limit distribution of the HAL-TMLE, {\em uniformly in all $P_0\in {\cal M}$}. This would mean that there exist sufficient sample sizes to obtain a particular level of precision in approximating  the normal limit distribution, uniformly in all $P_0\in {\cal M}$. This further demonstrates the robustness of the HAL-MLE, HAL-TMLE, and its bootstrap distribution in statistical models that have uniform model bounds $(M_1,M_2,M_3)$.

\subsection{The nonparametric bootstrap for the exact second-order expansion of the HAL-TMLE}
Recall the exact second-order expansion of the HAL-TMLE:
\begin{equation}\label{displ}
n^{1/2}(\Psi(Q_n^*)-\Psi(Q_0))=n^{1/2}(P_n-P_0)D^*(Q_n^*,G_n)+n^{1/2}R_{20}(Q_n^*,G_n,Q_0,G_0).\end{equation}
Recall that $R_{20}(Q,G,Q_0,G_0)=R_{2P_0}(Q,G,Q_0,G_0)$ potentially depends on $P_0$ beyond $(Q_0,G_0)$. 
Typically, we have 
\begin{equation}
R_{2P_0}(Q,G,Q_0,G_0)=P_0 R_2(Q,G,Q_0,G_0)\mbox{ for some $R_2(Q,G,Q_0,G_0)$.}\label{P0repR2}
\end{equation}
Let $R_{2n}()=R_{2P_n}()$ be obtained by replacing the $P_0$ by the empirical measure $P_n$. Thus, if we have (\ref{P0repR2}), then 
$R_{2n}(Q,G,Q_n,G_n)=P_n R_2(Q,G,Q_n,G_n)$.
 We assume the analogue of the bound (\ref{boundingR2}) on $R_{20}$ for $R_{2n}$:
\begin{equation}\label{boundingR2n}
\mid R_{2n}(Q,G,Q_n,G_n)\mid \leq f({\bf d}_{n1}^{1/2}(Q,Q_n),{\bf d}_{n2}^{1/2}(G,G_n))\end{equation}
for some function $f:\openr^K_{\geq 0}\rightarrow\openr_{\geq 0}$, $K=K_1+K_2$ of the form $f(x)=\sum_{i,j} a_{ij} x_ix_j$, a quadratic polynomial with positive coefficients $a_{ij}\geq 0$. 

Consider the nonparametric bootstrap analogue of the right-hand side of (\ref{displ}):
\[
Z_n^{2,\#}=n^{1/2}(P_n^{\#}-P_n)D^*(Q_n^{\#*},G_n^{\#})+n^{1/2}R_{2n}(Q_n^{\#*},G_n^{\#},Q_n^*,G_n).\]
This bootstrap sampling distribution $Z_n^{\#*}$ provides a very direct estimate of the sampling distribution of $n^{1/2}(\Psi(Q_n^*)-\Psi(Q_0))$.

Let $\Phi_n^{\#}(x)=P(Z_n^{2,\#}\leq x\mid P_n)$  be the cumulative distribution of this bootstrap sampling distribution.
So a bootstrap based 0.95-confidence interval for $\psi_0$ is given by \begin{equation}\label{ciexactexp}
[\psi_n^{*}+q_{0.025,n}^{\#}/n^{1/2},\psi_n^*+q_{0.975,n}^{\#}/n^{1/2} ],\end{equation}
 where 
$q_{\alpha,n}^{\#}=\Phi_n^{\#-1}(\alpha)$ is the $\alpha$-quantile of this bootstrap distribution.

\begin{theorem}\label{thnpbootexactexp}
Under the same conditions as Theorem \ref{thnpboothaltmle} and condition (\ref{boundingR2n}), we have 
$Z_n^{2,\#}=n^{1/2}(P_n^{\#}-P_n)D^*(Q_n,G_n)+O_P(n^{-1/2-\alpha(d)})$, and thereby
\[
Z_n^{2,\#}\Rightarrow_d N(0,\sigma^2_0)\mbox{ conditional on $(P_n:n\geq 1)$}.\]
In particular the above confidence interval (\ref{ciexactexp}) contains $\psi_0$ with probability tending to 0.95 as $n\rightarrow\infty$.
\end{theorem}

One might simplify $Z_n^{2\#}$  by replacing the targeted versions by their initial estimators:
\begin{equation}\label{simplified}
Z_n^{2a,\#}=n^{1/2}(P_n^{\#}-P_n)D^*(Q_n^{\#},G_n^{\#})+n^{1/2}R_{2n}(Q_n^{\#},G_n^{\#},Q_n,G_n).\end{equation}
In this case $Z_n^{2a,\#}$ is the bootstrap sampling distribution of the exact second-order expansion 
\[
n^{1/2}(\psi_n^1-\Psi(Q_0))=n^{1/2}(P_n-P_0)D^*(Q_n,G_n)+R_{20}(Q_n,G_n,Q_0,G_0)\]
of the HAL-one-step estimator $\psi_n^1=\Psi(Q_n)+P_n D^*(Q_n,G_n)$.  The latter bootstrap sampling distribution can also  be used for the HAL-TMLE.

As above, let $\Phi_n^{a\#}(x)=P(Z_n^{2a,\#}\leq x\mid P_n)$  be the cumulative distribution of $Z_n^{2a,\#}$, conditional on $P_n$.
A corresponding bootstrap based 0.95-confidence interval for $\psi_0$ is given by \begin{equation}\label{ciexactexpa}
[\psi_n^1+q_{0.025,n}^{a\#}/n^{1/2},\psi_n^1+q_{0.975,n}^{a\#}/n^{1/2} ],\end{equation}
 where 
$q_{\alpha,n}^{a\#}=\Phi_n^{a\#-1}(\alpha)$. 
We have the analogue of the above theorem for the bootstrap distribution $Z_n^{2a,\#}$ of the one-step estimator, where we can remove the specific conditions needed for the TMLE
$Q_n^*$.
Since the proof is remarkably simple and demonstrates that $Z_n^{2,\#}$ and $Z_n^{2a,\#}$ provide very direct approximations of the sampling distributions of the HAL-TMLE and HAL-one-step estimator,  we show here its proof.

\begin{theorem}\label{thnpbootexactexpa}
Assume the conditions of Theorem \ref{thefftmle}; (\ref{boundb}); (\ref{dn1bounding}); and (\ref{contDstaremp}).

Then, $Z_n^{2a,\#}=n^{1/2}(P_n^{\#}-P_n)D^*(Q_n,G_n)+O_P(n^{-1/2-\alpha(d)})$, and thereby
\[
Z_n^{2a,\#}\Rightarrow_d N(0,\sigma^2_0)\mbox{ conditional on $(P_n:n\geq 1)$}.\]
In particular the above confidence interval (\ref{ciexactexpa}) contains $\psi_0$ with probability tending to 0.95 as $n\rightarrow\infty$.
\end{theorem}
{\bf Proof:}
Consider (\ref{simplified}).
By Theorem \ref{thnpbootmle} we have that $d_{n1}(Q_n^{\#},Q_n)$ and $d_{n2}(G_n^{\#},G_n)$ are $O_P(n^{-1/2-\alpha(d)})$.
Using the bound (\ref{boundingR2n}) implies now that $R_{2n}(Q_n^{\#},G_n^{\#},Q_n,G_n)=O_P(n^{-1/2-\alpha(d)})$. 
Regarding the leading term in $Z_n^{2a,\#}$ we write is as 
\[
n^{1/2}(P_n^{\#}-P_n)\{D^*(Q_n^{\#},G_n^{\#})-D^*(Q_n,G_n)\}+n^{1/2}(P_n^{\#}-P_n)D^*(Q_n,G_n).\]
In the proof of Theorem \ref{thnpboothaltmle} we showed that the first term is $O_P(n^{-1/2-\alpha(d)})$. This completes the proof of Theorem \ref{thnpbootexactexp}.
$\Box$

The above bootstrap distribution  $Z_n^{2,\#}$ is different from the  sampling distribution of $Z_n^{1,\#}=n^{1/2}(\Psi(Q_n^{\#*})-\Psi(Q_n^*))$ used in the previous subsection (and similarly, $Z_n^{2a,\#}$ is different from the bootstrap distribution of the standardized one-step estimator). The advantage of $Z_n^{1\#}$ is that it is an actual sampling distribution of our HAL-TMLE and thereby fully respects that our estimator is a substitution estimator. On the other hand, 
its asymptotic expansion as analyzed in the proof of Theorem \ref{thnpboothaltmle} and the remarkable direct and simple proof of Theorem \ref{thnpbootexactexp} suggests that the sampling distribution $Z_n^{1,\#}$  is more different from the desired sampling distribution 
of $n^{1/2}(\psi_n^*-\psi_0)$ than $Z_n^{2,\#}$. Therefore it will be of interest to compare both bootstrap methods through a simulation study.

As in our previous theorems, the above theorems also apply to the setting in which we replace $C^u$ by the  cross-validation selector $C_n$.

\section{The nonparametric bootstrap for a conservative finite sample bound of exact second-order expansion of HAL-TMLE or HAL-one-step-estimator}\label{sectupperb}

We have the following finite sample upper bound for our HAL-TMLE relative to its target $\Psi(Q_0)$:
\begin{eqnarray*}
 n^{1/2}\mid \Psi(Q_n^*)-\Psi(Q_0)\mid&\leq & \mid  n^{1/2}(P_n-P_0)D^*(Q_n^*,G_n) \mid \\
 &&+
 f({\bf d}_{01}(Q_n^*,Q_0),{\bf d}_{02}(G_n,G_0)) +n^{1/2}r_n\\
 &\equiv&X_n(Q_n^*,G_n)+n^{1/2}r_n,
\end{eqnarray*}
where we defined a process $X_n(Q,G)$ (suppressing its dependence on $P_0$).
Similarly, we have this upper bound for the HAL one-step estimator $\psi_n^1=\Psi(Q_n)+P_n D^*(Q_n,G_n)$:
\[
n^{1/2}\mid \psi_n^1-\Psi(Q_0)\mid \leq  X_n(Q_n,G_n).\]
Let's focus on the latter, which could just as well be used for the sampling distribution of the HAL-TMLE as well.

{\bf How is this upper bound conservative?}
This upper bound is conservative from various points of view. Firstly, the true second-order remainder $R_{20}(Q_n,G_n,Q_0,G_0)$ could have both negative and positive values  that could cancel out a positive or negative value of $(P_n-P_0)D^*(Q_n,G_n)$. For example, in many models $R_{20}$ has a double robust structure
$\int (H_1(Q_n)-H_1(Q_0))(H_2(G_n)-H_2(G_0)) H_3(P_0,P_n) dP_0$.
In these double robust problems $Q_n$ and $G_n$ are based on different factors of the likelihood so that $H_1(Q_n)-H_1(Q_0))$ is generally almost uncorrelated with $H_2(G_n)-H_2(G_0))$, so that such a second order term could be reasonably symmetric  distributed around zero. The above upper bound does not allow any cancelation making it particularly conservative for double robust estimation problems. Secondly, the actual size  of $R_{20}(Q_n,G_n,Q_0,G_0)$ could be significantly smaller than our upper bound. For example, if we have the double robust structure, then the upper bound  bounds a term $\int (H_1(Q_n)-H_1(Q_0))(H_2(G_n)-H_2(G_0)) H_3(P_0,P_n)dP_0$ by Cauchy-Schwarz while bounding $H_3$ by its supremum norm. Since $Q_0$ and $G_0$ are very different functions, the Cauchy-Schwarz bound is very conservative itself, and, the supremum norm bound on $H_3$ will involve replacing a denominator by its smallest value. 
Therefore,  this bound is highly conservative for double robust estimation problems. 

If the second-order remainder has the form $\int (H_1(Q_n)-H_1(Q_0))^2 H_3(P_n,P_0)dP_0$, then this bound is more reasonable by  only being conservative due to the bounding of $H_3$ by its supremum norm and that we do not allow cancelation of  
 a mean zero centered $n^{1/2}(P_n-P_0)D^*(Q_n,G_n)$ with $R_{20}(P_n^*,P_0)$.
 Finally, the sampling distribution of this upper bound is not incorporating the known bounds for $n^{1/2}(\Psi(Q_n^*)-\Psi(Q_0))$ such as, for example, that $\Psi(P_0)$. The first nonparametric bootstrap method of the previous section is a sampling distribution of a substitution estimator thereby respects all the global bounds of the model and target parameter (e.g., $\Psi(P)$ is a probability). Respecting global constraints is particularly important when the target parameter is weakly supported by the data and asymptotics has not kicked in for the given sample size. 
 
 \ 
 \newline
We estimate the distribution of this upper bound with the nonparametric bootstrap. That is, we (conservatively) approximate the sampling distribution of 
$Z_n=\mid n^{1/2}(\Psi(Q_n^*)-\Psi(Q_0)\mid$ or $\mid n^{1/2}(\psi_n^1-\Psi(Q_0))\mid$ with
\[
Z_n^{3,\#}=\mid  n^{1/2}(P_n^{\#}-P_n)D^*(Q_n^{\#},G_n^{\#})\mid +
 f({\bf d}_{n1}(Q_n^{\#},Q_n),{\bf d}_{n2}(G_n^{\#},G_n)) \mid \]
 conditional on $(P_n:n\geq 1)$.
This distribution can now be used to construct an $0.95$-confidence interval. Let $F_n^{\#}(x)=P(Z_n^{3,\#}\leq x\mid (P_n:n\geq 1))$ and 
$q_{n,0.95}^{\#}=F_n^{\#  -1}(0.95)$ be its $0.95$-quantile. Then, $\Psi(Q_n^*)\pm q_{n,0.95}^{\#}/n^{1/2}$ is the resulting $0.95$-confidence interval.
 
{\bf Alternative reasonable upper-bound:}
 It appears to also be reasonable to use as upper-bound of $X_{n1}(Q_n,G_n)\equiv \mid  n^{1/2}(P_n-P_0)D^*(Q_n^*,G_n)  +
 f(d_{01}(Q_n^*,Q_0),d_{02}(G_n,G_0))\mid $. Note that  $X_{n1}(Q_n,G_n)$ is different from $X_n(Q_n,G_n)$, by putting the absolute value outside the sum of the two terms. This is not a deterministic upper bound in the sense that $\mid n^{1/2}(\psi_n^1-\Psi(Q_0))\mid $ is smaller than this bound with probability 1, but we certainly expect this to be a conservative distribution since we are adding a positive bias to a mean zero centered symmetric empirical process at $D^*(Q_n,G_n)$.

In spite of the conservative nature of our upper bound it is still asymptotically sharp.
Again,  not surprisingly, asymptotic consistency of the above conservative sampling distribution  is an immediately corollary of our analysis of the nonparametric bootstrap for the HAL-TMLE. 
\begin{theorem}\label{thupperbound}
Under the same conditions as Theorem \ref{thnpboothaltmle}, $Z_n^{3,\#}=n^{1/2}(P_n^{\#}-P_n)D^*(Q_n,G_n)+O_P(n^{-1/2-\alpha(d)})$, and thereby
\[
Z_n^{3,\#}\Rightarrow_d N(0,\sigma^2_0)\mbox{ conditional on $(P_n:n\geq 1)$}.\]
In particular the above confidence interval (\ref{ciexactexp}) contains $\psi_0$ with probability tending to 0.95 as $n\rightarrow\infty$.
 \end{theorem}
 
\paragraph{More conservative asymptotically sharp sampling distribution and corresponding bootstrap method:}
In the Appendix \ref{sectsupupperb} we propose an even more conservative sampling distribution for $n^{1/2}(\psi_n^1-\psi_0)$ (or $n^{1/2}(\psi_n^*-\psi_0)$) by replacing $X_n(Q_n,G_n)$ by   the supremum of $X_n(Q,G)$ over all $(Q,G)$ in the parameter space for which $d_{01}(Q,Q_0)$ and $d_{02}(G,G_0)$ are smaller than specified constants $x_{1n} $ and $x_{2n}$ chosen so that the probability that $d_{01}(Q_n,Q_0)$ and $d_{02}(G_n,G_0)$ are smaller than these constants are known to be larger than $1-\bar{\alpha}_{n}$ for some small number $\bar{\alpha}_n\rightarrow 0$. We propose a concrete method that expresses $(x_{1n},x_{2n})$  in terms of  a quantile of the supremum norm of a standard empirical cumulative survival function process $n^{1/2}(\bar{P}_n-\bar{P}_0)$. Bootstrapping this latter process provides now an estimator $(x_{1n}^{\#},x_{2n}^{\#})$ of $(x_{1n},x_{2n})$. The distribution of the corresponding supremum of $X_n(Q,G)$  is then estimated with the nonparametric bootstrap as well. The same method could be applied to $X_{n1}(Q,G)$.
Since this supremum might be cumbersome to compute in practice, in Appendix \ref{sectsupupperb} we proceed with proposing a simplified conservative approximation of this supremum in which the second-order remainder is separately maximized by plugging in the values $(x_{1n},x_{2n})$, resulting in an easy to compute sampling distribution. Again, we show that both of these methods are still asymptotically sharp.  
 
\section{Examples}\label{sectexample}

\subsection{Nonparametric estimation of average treatment effect}
Let $O=(W,A,Y)\sim P_0$, where $W\in [0,\tau_1]\subset \openr^{m_1}_{\geq 0}$ is an $m_1$-dimensional vector of baseline covariates, $A\in \{0,1\}$ is a binary treatment, and $Y\in \{0,1\}$ is a binary outcome. For a possible data distribution $P$, let $\bar{Q}(P)=E_P(Y\mid A,W)$, $G(P)=P(A=1\mid W)$, and let $Q_W(P)$ be the cumulative  probability distribution of $W$. Let $Q=(Q_W,\bar{Q})$. Let $g(a\mid W)=P(A=a\mid W)=G(W)^a(1-G(W))^{1-a}$.
Thus $Q_1=Q_W$, $Q_2=\bar{Q}$, $m_{11}=m_1$ and $m_{12}=m_1+1$, in terms of our general notation. Suppose that our model assumes that $G(W)$ depends on a  possible subvector of $W$, and let $m_2 $ be the dimension of this subvector.

{\bf Statistical model:}
 Since $Q_W$ is a cumulative distribution function  it is a monotone $m_1$-variate cadlag function and its sectional variation norm equals its total variation which thus equals 1. 
Let $\delta>0$ be given. 
We assume  $\bar{Q}\in (\delta,1-\delta)$ and that it is  an element of the class of $m_{12}$-dimensional cadlag functions with sectional variation norm bounded by some $C_{12}^u$. (here one can treat $A$ as continuous on $[0,1]$ and assume that $\bar{Q}$ is a step-function in $A$ with single jump at 1, allowing us to embed functions of continuous and discrete covariates in a cadlag function space.)
Similarly, we assume $G\in (\delta,1-\delta)$ and that it is an element of the class of $m_2$-dimensional cadlag functions with sectional variation norm bounded by a $C_2^u$.
Let's denote these parameter spaces for $Q_W,\bar{Q}$ and $G$ with ${\cal F}_{11}$, ${\cal F}_{12}$ and ${\cal F}_2$, respectively.
Let ${\cal F}_1={\cal F}_{11}\times {\cal F}_{12}$ be the parameter space of $Q=(Q_W,\bar{Q})$.
For a given $C_1^u=(C_{11}^u=1,C_{12}^u),C_2^u<\infty$ and $\delta>0$, consider the statistical model \[
{\cal M}=\{P: Q_W\in {\cal F}_{11}, \bar{Q}\in {\cal F}_{12}, G\in {\cal F}_2\}.\]
 Thus, ${\cal M}$ is defined as the set of all possible probability distributions for which the conditional means of $Y$ and $A$ are cadlag functions with sectional variation norm bounded by $C_1$ and $C_2$, respectively, and the conditional density of $A$, given $W$, is bounded away from $0$ and $1$,  $P_W$-a.e, while we make no assumptions on the probability distribution of $W$.
 
{\bf Parameter space of type (\ref{calFmodel}) or (\ref{calFmodelplus})}: As shown in Section 2, we can reparametrize $G=\delta+(1-2\delta)\mbox{expit}(f_2(G)(W))$ and $\bar{Q}=\delta+(1-2\delta)\mbox{expit}(f_1(Q))$, where now $f_1$ and $f_2$ can be any cadlag function that is only restricted by  upper bounds on their sectional variation norm  implied by $C_{12}^u$ and
$C_2^u$, while $C_{12}^l=C_2^l=0$,  so that the parameter space for $f_1$ and $f_2$ is indeed of type (\ref{calFmodel}). Obviously, $Q_W$ is of the type (\ref{calFmodelplus}) with $C_{11}^l=C_{11}^u=1$. This demonstrates that our model ${\cal M}$ can be represented as a model as defined in Section 2.
 
{\bf Target parameter:}
Let $\Psi:{\cal M}\rightarrow\openr$ be defined by $\Psi(P)=\Psi_1(P)-\Psi_0(P)$, where $\Psi_a(P)=E_PE_P(Y\mid A=a,W)$.
Note that $\Psi(P)$ only depends on $P$ through $Q(P)$, so that we will also use the notation $\Psi(Q)$ instead of $\Psi(P)$.
Let's focus on $\Psi_1(P)$ which will also imply the formulas for $\Psi_0(P)$ and thereby $\Psi(P)$.

{\bf Loss functions for $Q$ and $G$:}
Let $L_{11}(Q_W)=\int_x(I(W\leq x)-Q_W(x))^2r(x) dx$ for some weight function $r>0$  be the loss function for $Q_{W,0}$. 
Let $d_{011}(Q_W,Q_{W,0})=P_0L_{11}(Q_W)-P_0L_{11}(Q_{W,0})$ be the corresponding  loss based dissimilarity.
Let $L_{12}(\bar{Q})=-\{Y\log\bar{Q}(A,W)+(1-Y)\log(1-\bar{Q}(A,W))\}$ be the log-likelihood loss function for the conditional mean $\bar{Q}_0$, and let
$d_{012}(\bar{Q},\bar{Q}_0)=P_0 L_{12}(\bar{Q})-P_0L_{12}(\bar{Q}_0)$ be the corresponding Kullback-Leibler dissimilarity.
We can then define the sum-loss $L_1(Q)=L_{11}(Q_W)+L_{12}(\bar{Q})$ for $Q_0$, and its loss-based dissimilarity
 $d_{01}(Q,Q_0)=P_0L_1(Q)-P_0L_1(Q_0)$ which equals the sum of the following two dissimilarities
 \begin{eqnarray*}
d_{012}(Q,Q_0)&=&\int_x (Q_{W}(x)-Q_{W,0}(x))^2 r(x) dx\\
d_{011}(Q_W,Q_{W,0})&=&\int \log \left(\frac{\bar{Q}_0}{\bar{Q}}\right)^y \left( \frac{1-\bar{Q}_0}{1-\bar{Q}}\right)^{1-y}(a,w) 
 dP_0(w,a,y) .\end{eqnarray*}
Let $L_2(G)=-\{A\log G(W)+(1-A)\log(1-G(W))\}$ be the loss function for $G_0=P_0(A=1\mid W)$, and let $d_{02}(G,G_0)=P_0L_2(G)-P_0L_2(G_0)$ be the Kullback-Leibler dissimilarity between $G$ and $G_0$.

{\bf Canonical gradient and corresponding exact second order expansion:}
The canonical gradient of $\Psi_a$ at $P$ is given by:
\[
D^*_a(Q,G)=\frac{I(A=a)}{g(A\mid W)}(Y-\bar{Q}(A,W))+\bar{Q}(1,W)-\Psi_a(Q).\]
The exact second-order remainder $R_{20}^a(P,P_0)\equiv \Psi_a(P)-\Psi_a(P_0)+P_0 D^*_a(P)$ is given by:
\[
R_{20}^a(\bar{Q},G,\bar{Q}_0,G_0)=\int \frac{(g-g_0)(a\mid w)}{g(a\mid w)}(\bar{Q}-\bar{Q}_0)(a,w) dP_0(w).\]

{\bf Bounding the second order remainder:}
By using Cauchy-Schwarz inequality, we obtain the following bound on $R_{20}^a(P,P_0)$:
\[
\mid R_{20}^a(P,P_0)\mid \leq \delta^{-1}\pl \bar{Q}_a-\bar{Q}_{a0}\pl_{P_0}\pl G-G_0\pl_{P_0} ,\]
where $\bar{Q}_a(W)=\bar{Q}(a,W)$, $a\in \{0,1\}$.
Thus, $D^*(P)=D^*_1(P)-D^*_0(P)$, $R_{20}(P,P_0)=R_{20}^1(P,P_0)-R_{20}^0(P,P_0)$, and the upper bound for $R_{20}(P,P_0)$ can be defined as the sum of the two upper bounds for $R_{20}^a(P,P_0)$ in the above inequality, $a\in \{0,1\}$.

By \citep{vanderVaart98} we have $\pl p^{1/2}-p_0^{1/2}\pl_{P_0}^2\leq  P_0 \log p_0/p $. For Bernoulli distributions,  we have
$\pl p-p_0\pl^2_{P_0}\leq 4 \pl p^{1/2}-p_0^{1/2}\pl^2_{P_0}\leq P_0\log p_0/p$. 
From this it follows that $\int (\bar{Q}-\bar{Q}_0)^2(a,w)dP_0(a,w)\leq 4 d_{012}(\bar{Q},\bar{Q}_0)$ and thus 
$\pl \bar{Q}_a-\bar{Q}_{a0}\pl^2_{P_0}\leq 4\delta^{-1}d_{012}(\bar{Q},\bar{Q}_0)$.
Therefore, $\pl \bar{Q}_a-\bar{Q}_{a0}\pl_{P_0}\leq 2\delta^{-1/2} d_{012}^{1/2}(\bar{Q},\bar{Q}_0)$.
Similarly, it follows that $\pl G-G_0\pl_{P_0}\leq 2 d_{02}^{1/2}(G,G_0)$.
This thus shows the following bound on $R_{20}^a(P,P_0)$:
\[
\mid R_{20}^a(P,P_0)\mid \leq 4\delta^{-1.5} d_{012}^{1/2}(\bar{Q},\bar{Q}_0) d_{02}^{1/2}(G,G_0).\]
The right-hand side represents the function $f({\bf d}_{01}^{1/2}(Q,Q_0),{\bf d}_{02}^{1/2}(G,G_0))$ for the parameter $\Psi_a$  in our general notation: $f(x=(x_1,x_2),y)=
4 \delta^{-1.5} x_2 y$.
The sum of these two bounds for $a\in \{0,1\}$ (i.e, $2f()$) provides now a conservative bound for  $R_{20}=R_{20}^1-R_{20}^0$:
\begin{equation}\label{r2upperboundexample1}
\mid R_{20}(P,P_0)\mid \leq f(d_{012}^{1/2}(\bar{Q},\bar{Q}_0),d_{02}(G,G_0))\equiv 8\delta^{-1.5} d_{012}^{1/2}(\bar{Q},\bar{Q}_0) d_{02}^{1/2}(G,G_0).\end{equation}
This verifies (\ref{boundingR2}).
We note that this bound is very conservative due to the arguments we provided in general in the previous section for double robust estimation problems.

{\bf Continuity of canonical gradient:}
Regarding the continuity assumption (\ref{contDstar}), we note that $P_0\{D^*_a(P)-D^*_a(P_0))^2$ can be bounded by $\pl G-G_0\pl_{P_0}^2+\pl \bar{Q}_a-\bar{Q}_{a0}\pl^2_{P_0}$ and $(\Psi_a(Q)-\Psi_a(Q_0))^2$, where the constant depends on $\delta$. The latter square difference can be bounded in terms of $\pl \bar{Q}_a-\bar{Q}_{a0}\pl^2_{P_0}$ and  by applying our integration by parts formula to  $\int \bar{Q}_a(w) d(Q_W-Q_{W0})(w)$
by $d_{011}(Q_W,Q_{W0})$, where the constant depends on $C_1^u$. Thus this proves (\ref{contDstar}) for $D^*=D^*_1-D^*_0$.

{\bf Uniform model bounds on sectional variation norm:} It also follows immediately that the sectional variation norm model bounds $M_1,M_2,M_3$ (\ref{sectionalvarbound})
of $L_1(Q)$, $L_2(G)$ and $D^*(P)$ are all finite, and can be expressed in terms of $(C_1^u,C_2^u,\delta)$. This verifies the model assumptions of Section 2.


{\bf HAL-MLEs:}
 Let $Q_n=\arg\min_{Q\in {\cal F}_1}P_n {\bf L}_1(Q)$ and $G_n=\arg\min_{G\in {\cal F}_2}P_n L_2(G)$ be the HAL-MLEs. Here we can use the above mentioned reparameterizations of $Q$ and $G$ in terms of $f_1$ and $f_2$, respectively, that varies over a parameter space of type $(\ref{calFmodel})$.
 As shown in \citep{vanderLaan15,Benkeser&vanderLaan16}, if one simply sets $\delta=0$, then $\bar{Q}_n$ and $G_n$ can be computed with standard Lasso logisitic regression software using a linear logistic regression model with around $n 2^{m_1}$ indicator basis functions, where $m_1$ is the dimension of $W$.
 The reparameterization would now enforce the bounds $\delta$ and $1-\delta$ for these HAL-MLEs.

 Note that $Q_{W,n}$ is just an unrestricted MLE and thus equals the empirical cumulative distribution function. 
 Therefore, we actually have that $\pl Q_{W,n}-Q_{W,0}\pl_{\infty}=O_P(n^{-1/2})$ in supremum norm, while $d_{012}(\bar{Q}_n,\bar{Q}_0)$ and $d_{02}(G_n,G_0)=O_P(n^{-1/2-\alpha(d)})$
 where $d$ is the dimension of $O$. If $m_2<d-2$, then one  should be able to improve the bound into $n^{-1/2-\alpha(m_2)}$.
 
 {\bf CV-HAL-MLEs:}
 The above HAL-MLEs are determined by $(C_1^u=(1,C_{12}^u),C_2^u)$ and could thus be denoted with $Q_{n,C_1^u}=\hat{Q}_{C_1^u}(P_n)$ and $G_{n,C_2^u}=\hat{G}_{C_2^u}(P_n)$.
 Let $C_{10}=\pl Q_0\pl_v^*=(1,\pl \bar{Q}_0\pl_v^*)$ and $C_{20}=\pl G_0\pl_v^*$, respectively, which are thus smaller than $C_1^u$ and $C_2^u$, respectively.
 We can now define the cross-validation selector that selects the best HAL-MLE over all $C_{1}$ and $C_2 $ smaller than these upper-bounds:
 \begin{eqnarray*}
 C_{1n}&=&\arg\min_{C_{11}=1,C_{12}<C_{12}^u}E_{B_n}P_{n,B_n}^1L_1(\hat{Q}_{C_1}(P_{n,B_n}^0)) \\
C_{2n}&=&\arg\min_{C_2<C_2^u}E_{B_n}P_{n,B_n}^1L_2(\hat{G}_{C_2}(P_{n,B_n}^0)),
\end{eqnarray*}
where $B_n\in \{0,1\}^n$ is a random split in training sample $\{O_i:B_n(i)=0\}$ with empirical measure $P_{n,B_n}^0$ and validation sample
$\{O_i:B_n(i)=1\}$ with empirical measure $P_{n,B_n}^1$.
This defines now the CV-HAL-MLE $Q_n=Q_{n,C_{1n}}$ and $G_n=G_{n,C_{2n}}$ as well. Thus, by setting $C_1^u=C_{1n}$ and $C_2^u=C_{2n}$, our HAL-MLEs equal the CV-HAL-MLE.

 {\bf HAL-TMLE:}
 Let $\mbox{Logit}\bar{Q}_{n,\epsilon}=\mbox{Logit}\bar{Q}_n+\epsilon C(G_n)$, where $C(G_n)(A,W)=(2A-1)/g_n(A\mid W)$.
 Let $\epsilon_n=\arg\min_{\epsilon}P_n L_{11}(\bar{Q}_{n,\epsilon})$. This defines the TMLE $\bar{Q}_n^*=\bar{Q}_{n,\epsilon_n}$ of $\bar{Q}_0$. We can also define a local least favorable submodel $\{Q_{W,n,\epsilon_2}:\epsilon_2\}$ for $Q_{W,n}$ but since $Q_{W,n}$ is an NPMLE one will have that $\epsilon_{2n}=\arg\min_{\epsilon_2}P_n L_{11}(Q_{W,n,\epsilon_2})=0$, and thereby that the TMLE of $Q_0$ for any such 2-dimensional least favorable submodel is given by  
 $Q_n^*=(Q_{W,n},\bar{Q}_n^*)$. It follows that $P_n D^*(Q_n^*,G_n)=0$.
 
 {\bf Preservation of rate for HAL-TMLE:}
The proof in the Appendix \ref{AppendixA} for $\epsilon_n=O_P(d_{012}^{1/2}(\bar{Q}_n,\bar{Q}_0))$ applies to this submodel, so that indeed
$d_{01}(Q_n^*,Q_0)$  converges at same rate as $d_{01}(Q_n,Q_0)$.

{\bf Asymptotic efficiency of HAL-TMLE and CV-HAL-TMLE:}
 Application of Theorem \ref{thefftmle} shows that $\Psi(Q_n^*)$ is asymptotically efficient, where one can either choose $Q_n$ as  a fixed HAL-MLE using $C_1=C_1^u$ or the CV-HAL-MLE using $C_1=C_{1n}$, and similarly, for $G_n$. The preferred estimator would be the CV-HAL-TMLE.
 
 {\bf Finite sample conservative confidence interval:}
 Let's first consider the exact finite sample conservative confidence interval presented in (\ref{consbound1}).
 For this we need bounds $M_3$, $M_{12}$ and $M_2$ on the sectional variation norm of $D^*(Q,G)$, $L_{12}(\bar{Q})$ and $L_2(G)$, respectively.
 These can be expressed in terms of the sectional variation norm bounds $(C_{12}^u,C_2^u)$ on $(\bar{Q}, G)$ and the lower bound $\delta$ of $\min_a g(a\mid W)$ and $\bar{Q}$.
 Here one can use that the sectional variation norm of $1/g(a\mid W)$ can be bounded in terms of $\delta$ and $\pl w\rightarrow g(a\mid w)\pl_v^*$. 
 (\ref{consbound1}) tells us that $\mid n^{1/2}(\psi_n^1-\Psi(Q_0))\mid $ is dominated by the distribution of 
 $Z_n^+=(M_3+f(M_{12}^{1/2},M_2^{1/2}))\pl n^{1/2}(\bar{P}_n-\bar{P}_0)\pl_{\infty}$, where
$f$ is defined by (\ref{r2upperboundexample1}), and 
$\bar{P}(u)=P([u,\tau])=\int_{[u,\tau]} dP(s)$ is the probability that $O\in [u,\tau]$ under $P$. Estimation of the sampling distribution of
$n^{1/2}(\bar{P}_n-\bar{P}_0)$ with $n^{1/2}(\bar{P}_n^{\#}-\bar{P}_n)$ results then in the estimate $Z_n^{+,\#}$ of $Z_n^+$ and corresponding finite sample conservative $0.95$-confidence interval.
(Similarly, we have this bound for the TMLE $\Psi(Q_n^*)$ with $M_3$  and $M_{12}$ replaced by $M_3^*$ and $M_{12}^*$, respectively.)
 However, as we argued in general, this conservative confidence interval will generally be of little practical use by being much too conservative, although these confidence intervals will still have a width or order $n^{-1/2}$.

 {\bf Asymptotic validity of the nonparametric bootstrap for the HAL-MLEs:}
 Firstly, note that the bootstrapped HAL-MLEs \[
 \bar{Q}_n^{\#}=\arg\min_{\pl \bar{Q}\pl_v^*< C_{12}^u,\bar{Q}\ll \bar{Q}_n}P_n^{\#}L_{12}(\bar{Q}),\]
  and 
 $G_n^{\#}=\arg\min_{\pl G\pl_v^*<C_2^u,G\ll G_n}P_n^{\#}L_2(G)$ are easily computed as a standard Lasso regression using $L_1$-penalty $C_{12}^u$ and $C_2^u$ and including the maximally $n$ indicator basis functions with the non-zero coefficients selected by $Q_n$ and $G_n$, respectively. This makes the actual computation of the nonparametric bootstrap distribution a very doable computational problem, even though the single computation of $Q_n$  and $G_n$ is highly demanding for large dimension of $W$.
 
{\bf Verification of conditions for validity of bootstrap for HAL-MLE Theorem \ref{thnpbootmle}:}
  We now want to verify the conditions for our asymptotic consistency of the nonparametric bootstrap of Theorem \ref{thnpbootmle}.
This requires us to establish a second-order expansion for $d_{012}(\bar{Q}_n^*,\bar{Q}_0)$ and $d_{n12}(\bar{Q}_n^{\#*},\bar{Q}_n)$, so that we can specify
the second-order terms $P_0 R_{20,L_{12}}(\bar{Q}_n^*,\bar{Q}_0)$ and $P_n R_{2n,L_{12}}(\bar{Q}_n^{\#*},\bar{Q}_n)$. Subsequently, we have to bound
the square of the $L^2(P_0)$-norm and $L^2(P_n)$ norm of the loss-differences $L_{12}(\bar{Q}_n^*,\bar{Q}_0)$ and $L_{12}(\bar{Q}_n^{\#*},\bar{Q}_n)$ in terms of these second-order terms.
Since this concerns a log-likelihood loss, we consider this problem in general. Let $L(p)=-\log p$ be the log-likelihood loss. 
Firstly, we consider the exact second-order Tailor expansion of $P_0L(p)-P_0L(p_0)$ at $p_0$:
\[
P_0\log p-P_0\log p_0=\int p_0^{-1}(p-p_0) dP_0-P_0 R_{20,L}(P,P_0).\]
Since the first order linear term equals zero, it follows that $P_0R_{20,L}(P,P_0)=-P_0\log p/p_0$. Thus, indeed $P_0 \{L(p)-L(p_0)\}^2$ can be bounded by $P_0 R_{20,L}(P,P_0)$ due to known result that the Kullback-Leibler dissimilarity is equivalent with $\int (p-p_0)^2 d\mu$ if the densities are bounded away from 0.
Similarly, the exact second-order Tailor expansion of $P_n L(p)-P_n L(p_n)$ at $p_n$ is given by:
\[
P_n \log p-P_n \log p_n=P_n p_n^{-1}(p-p_n) -P_n R_{2n,L}(p,p_n) \]
for an exact second-order remainder $P_n R_{2n,L}(p,p_n)$.
By the exact second-order Tailor expansion of the function $\log x$ at $x=p_n(o)$, we obtain
\[
\log p(o)-\log p_n(o)=p_n^{-1}(p-p_n)(o) -\xi(p_n(o),p(o))^{-2} (p-p_n)^2(o),\]
where $\xi(p_n(o),p(o))$ is a value in between $p_n(o)$ and $p(o)$.
Thus, $P_n R_{2n,L}(p,p_n)=P_n \xi(p_n,p)^{-2}(p-p_n)^2$.
If $\min(p_n,p)>\delta$ for some $\delta>0$, then $P_n (p-p_n)^2\lesssim P_n R_{2n,L}(p,p_n)$.
It follows trivially that $P_n (L(p_n^{\#}-L(p_n))^2\lesssim P_n (p_n^{\#}-p_n)^2$, and thus also that 
$P_n (L(p_n^{\#})-L(p_n))^2\lesssim P_n R_{2n,L}(p,p_n)$.
This verifies the conditions on the loss function for the bootstrap theorem \ref{thnpbootmle}.

{\bf Behavior of HAL-MLE under sampling from $P_n$:}
This  shows that $d_{n12}(\bar{Q}_n^{\#},\bar{Q}_n)=O_P(n^{-1/2-\alpha(d)})$, and that this dissimilarity  is equivalent with the square $L^2$-norms $P_n R_{2n,L_{12}}(\bar{Q}_n^{\#},\bar{Q}_n)$, which is equivalent with $\sum_a \int (\bar{Q}_{na}^{\#}-\bar{Q}_{na})^2 dP_n$.
Theorem \ref{thnpbootmle} also shows that $d_{n2}(G_n^{\#},G_n)=O_P(n^{-1/2-\alpha(d)})$ and that this loss based dissimilarity is equivalent with $P_n (G_n^{\#}-G_n)^2$.

{\bf Preservation of rate of TMLE under sampling from $P_n$:}
The proof in the Appendix \ref{AppendixC} for $\epsilon_n^{\#}=O_P(d_{n12}^{1/2}(\bar{Q}_n^{\#},\bar{Q}_n))$ applies to our smooth submodel, so that indeed
$d_{01}(Q_n^{\#*},Q_0)$  converges at same rate as $d_{01}(Q_n,Q_0)$.

{\bf Consistency of nonparametric bootstrap for HAL-TMLE:}
This verifies all conditions of Theorem \ref{thnpboothaltmle} which establishes the asymptotic efficiency and asymptotic consistency of the nonparametric bootstrap.

\begin{theorem}
 We have that $\Psi(Q_n^*)$ is asymptotically efficient, i.e. $n^{1/2}(\Psi(Q_n^*)-\Psi(Q_0))\Rightarrow_d N(0,\sigma^2_0)$, where
$\sigma^2_0=P_0 \{D^*(P_0)\}^2$.
In addition, conditional on $(P_n:n\geq 1)$, $Z_n^{1,\#}=n^{1/2}(\Psi(Q_n^{\#*})-\Psi(Q_n^*))\Rightarrow_d N(0,\sigma^2_0)$.
This can also be applied to the setting in which $C^u$ is replaced by the cross-validation selector $C_n^u$ defined above.
\end{theorem}

{\bf Consistency of nonparametric bootstrap for exact expansion of HAL-TMLE/HAL-one-step:}
Recall the exact second-order remainder of $\Psi(Q_n^*)-\Psi(Q_0)$ defined by $R_{20}(Q_n^*,G_n,Q_0,G_0)=(R_{20}^1-R_{20}^0)(Q_n^*,G_n,Q_0,G_0)$.
Let $R_{2n}(Q_n^{\#*},G_n^{\#},Q_n^*,G_n)$ be the nonparametric bootstrap analogue.
Thus,
\[
R_{2n}^a(Q_n^{\#*},G_n^{\#},Q_n^*,G_n)=\int \frac{(g_n^{\#}-g_n)(a\mid w)}{g_n^{\#}(a\mid w)}(\bar{Q}_n^{\#*}-\bar{Q}_n^*)(a,w) dP_n(w).\] 
Let 
\[
Z_n^{2,\#}=n^{1/2}(P_n^{\#}-P_n)D^*(Q_n^{\#*},G_n^{\#})+n^{1/2}R_{2n}(Q_n^{\#*},G_n^{\#},Q_n^*,G_n).\]

Consider also the conservative version
\[
Z_n^{3,\#}=\mid n^{1/2}(P_n^{\#}-P_n)D^*(Q_n^{\#*},G_n^{\#})\mid + f(d_{n12}(\bar{Q}_n^{\#*},\bar{Q}_n^*),d_{n2}(G_n^{\#},G_n)) \mid .\]
 where the upper bound $f()$ for the remainder is defined in  (\ref{r2upperboundexample1}).
 Application of Theorems \ref{thnpbootexactexp} and \ref{thupperbound} prove the asymptotic consistency of these two nonparametric bootstrap distributions.
 
 \begin{theorem}
 We have $Z_n^{2,\#}\Rightarrow_d N(0,\sigma^2_0)$ and $Z_n^{3,\#}\Rightarrow_d \mid N(0,\sigma^2_0)\mid$, conditional on $(P_n:n\geq 1)$.
 As a consequence, an $0.95$-confidence interval for $\psi_0$ based on $Z_n^{2,\#}$ and $Z_n^{3,\#}$ have asymptotic coverage 0.95 of $\psi_0$.
 This also applies to the setting in which $C^u$ is replaced by the cross-validation selector $C_n$.
 \end{theorem}
 Finally, we remark that our HAL-MLE is really indexed by the  model bounds $(C_{12}^u,C_2^u,\delta)$ and these might all three be unknown to the user. So in that case, we recommend to select
 all three with the cross-validation selector $(C_{12n},C_{2n},\delta_)$ and define the HAL-TMLE and bootstrap of the HAL-TMLE at this fixed  choice $(C_n,\delta_n,M_n)$.






\subsection{Nonparametric estimation of integral of square of density}
{\bf Statistical model,  target parameter, canonical gradient:}
Let $O\in \openr^d$ be a multivariate random variable with probability distribution $P_0$ with support $[0,\tau]$.
Let ${\cal M}$ be a  nonparametric model dominated by Lebesgue measure $\mu$, where we assume that for each $P\in {\cal M}$ its density
$p=dP/d\mu$ is bounded away from below by $\delta>0$ and from above by $M<\infty$. In addition, we assume that all densities are cadlag functions and have sectional variation norm bounded by $C^u<\infty$.
As  shown under the remark in Section 2, we can reparametrize $p=c(f)\{\delta+(M-\delta)\mbox{expit}(f) \}$, where $f$  can be any cadlag function with sectional variation norm bounded from above by some finite constant implied by $C^u$ (while $C^l=0$), in which case our model is of the type  (\ref{calFmodel}).
The target parameter $\Psi:{\cal M}\rightarrow\openr$ is defined as $\Psi(P)=E_Pp(O)=\int p^2(o)d\mu(o)$. 
This target parameter is pathwise differentiable at $P$ with canonical gradient 
\[
D^*(P)(O)=2 (p(O)-\Psi(P)).\]
{\bf Exact second order remainder:}
It implies the following exact second-order expansion:
\[
\Psi(P)-\Psi(P_0)=(P-P_0)D^*(P)+R_{20}(P,P_0),\]
where 
\[
R_{20}(P,P_0)\equiv -\int (p-p_0)^2 d\mu .\]
{\bf Loss function:}
As loss function for $p$ we could consider the log-likelihood loss $L(p)(O)=-\log p(O)$ with $d_0(p,p_0)=P_0\log p_0/p$.
We have $\pl p^{1/2}-p_0^{1/2}\pl_{P_0}^2\leq  P_0 \log p_0/p $ so that 
\begin{eqnarray*}
\mid R_{20}(P,P_0)\mid &=&\int (p-p_0)^2 d\mu \\
&=&\sup_x\frac{(p^{1/2}+p_0^{1/2})^{2}}{p_0} \int (p^{1/2}-p_0^{1/2})^2 dP_0\\
&\leq&  M/\delta P_0\log p_0/p =M/\delta d_0(p,p_0).
\end{eqnarray*}
Alternatively, we could consider the loss function 
\[
L(p)(O)=-2p(O)+\int p^2(o)d\mu(o) .\]
Note that this is  indeed a valid loss function with loss-based dissimilarity given by \begin{eqnarray*}
d_0(p,p_0)&=&P_0 L(p)-P_0L(p_0)\\
&=&-2\int p(o)p_0(o)d\mu(o)+\int p^2d\mu+2\int p_0^2 d\mu-\int p^2_0 d\mu\\
&=& \int (p-p_0)^2 d\mu .\end{eqnarray*}
{\bf Bounding second order remainder:}
Thus, if we select this loss function, then we have
\[
\mid R_{20}(P,P_0)\mid =d_0(p,p_0) .\]
In terms of our general notation, we now have $f(x)=x^2$ for the upper bound on $R_{20}$ so that $\mid R_{20}(P,P_0)\mid =f(d_0^{1/2}(p,p_0))$.
We will proceed with the latter loss function so that our bound on $R_2(P,P_0)$ is sharp.  
In addition, if we use this loss function, we do not need a lower bound $\delta $ for our densities so that we can set $\delta=0$ in our definition of the model ${\cal M}$.
The canonical gradient is indeed continuous in $P$ as stated in (\ref{contDstar}) and the bounds $M_1,M_2,M_3$ (\ref{sectionalvarbound}) are obviously finite and can be expressed in terms of $(C^u,M,\delta)$.
This verifies the assumptions on our model as stated in Section 2.

{\bf HAL-MLE and CV-HAL-MLE:}
Let $p_n=\arg\min_{P\in {\cal M}}P_n L(p)$ be the MLE. Using our reparameterization this can be computed as 
\[
f_n=\arg\min_{ f}P_nL(c(f)\{\delta+(M-\delta)\mbox{expit}f\}),\]
where $f$ can be represented by our general representation (\ref{Frepresentation}), $f(o)=f(0)+\sum_{s\subset \{1,\ldots,d\}}\int_{(0_s,o_s]} df_s(u_s)$,
and constrained to satisfy $\mid f(0)\mid+\sum_{s\subset\{1,\ldots,d\}}\int_{(0_s,\tau_s]} \mid df_s(u_s)\mid \leq C$  for a $C$ implied by $C^u$.
Let's denote this $f_n$ with $f_{n,C}$.
Thus, for a given $C$ computation of $f_{n,C}$ can be done with a Lasso type algorithm. Let $C_n=\arg\min_CE_{B_n}P_{n,B_n}^1L(\hat{p}_C(P_{n,B_n}^0))$ be the cross-validation selector of $C$, as defined in previous example. If we set $C=C_n$, then we obtain the CV-HAL-MLE $f_n=f_{n,C_n}$.
We have $d_0(p_n,p_0)=O_P(n^{-1/2-\alpha(d)})$.

{\bf Universal least favorable submodel:}
We now define the HAL-TMLE. Consider the universal least favorable submodel $\{p_{n,\epsilon}:\epsilon\}$ through the HAL-MLE $p_n$: for $\epsilon\geq 0$
\[
p_{n,\epsilon}=p_n\exp\left(\int_0^{\epsilon} D^*(p_{n,x}) dx\right).\]
This submodel recursively defines $p_{n,\epsilon}$ where one starts calculating $p_{n,dx}$ for an infinitesimal $dx>0$ from $p_n$, and then $p_{n,2dx}$ from $p_{n,dx}$ and $p_n$ etc. This recursive definition generates $\{p_{n,\epsilon}:\epsilon \geq 0\}$. Similarly, one computes $p_{n,-dx}$ from $p_n$, and $p_{n,-2dx}$ from $p_n,p_{n,-dx}$ etc, 
where for $\epsilon<0$  $\int_0^{\epsilon}=-\int_{\epsilon}^0$.  
One can also define  this universal least favorable submodel by recursively applying a local least favorable submodel:
\[
p_{n,\epsilon+d\epsilon}=p_{n,\epsilon,d\epsilon}^{lfm},\]
where $p_{x}^{lfm}$ is a local least favorable submodel through $p$ at parameter value $x$ so that $p_{n,\epsilon,d\epsilon}^{lfm}$ is the local least favorable submodel through $p_{n,\epsilon}$ at parameter value $d\epsilon$. 
A possible local least favorable submodel choice is $p_{x}=(1+x D^*(p)) p$ for $x$ in a small neighborhood around $0$.

{\bf HAL-TMLE:}
Let $\epsilon_n=\arg\min_{\epsilon}P_n L(p_{n,\epsilon})$ be the MLE, and $p_n^*=p_{n,\epsilon_n}$ is the TMLE.
The TMLE of $\Psi(P_0)$ is the plug-in estimator $\psi_n^*=\Psi(P_n^*)=\int p_n^{*2}d\mu$.
It is easily verified that the universal least favorable submodel $p_{\epsilon}=p\exp(\int_0^{\epsilon}D^*(p_x) dx)$ is such that $\log p_{\epsilon}$ is twice differentiable in $\epsilon$. Therefore we can carry out the general proof in the Appendix \ref{AppendixA} establishing that $\epsilon_n=O_P(n^{-1/4-\alpha(d)/2})$.

{\bf Efficiency of HAL-TMLE and CV-HAL-TMLE:}
Application of Theorem \ref{thefftmle} shows that $\Psi(P_n^*)$ is asymptotically efficient, where one can either choose the HAL-MLE with fixed index $C$ implied by $C^u$ or one can set $C=C_n$ equal to cross-validation selector defined above.

{\bf Finite sample conservative confidence interval:}
Let's first consider the exact finite sample conservative confidence interval presented in (\ref{consbound1}).
 For this we need bounds $M_3$, $M_{1}$ on the sectional variation norm of $D^*(P)=2p-\Psi(p)$ and $L(p)=-2p+\Psi(p)$, respectively. These bounds will thus be identical: $M_3=M_1$. We have that $M_1=\sup_{p\in {\cal M}} \pl L(p)\pl_v^*=2 C^u$.
 (\ref{consbound1}) tells us that $\mid n^{1/2}(\psi_n^1-\Psi(P_0))\mid $ is dominated by the distribution of 
 $Z_n^+=4C^u\pl n^{1/2}(\bar{P}_n-\bar{P}_0)\pl_{\infty}$, where 
$\bar{P}(u)=P([u,\tau])=\int_{[u,\tau]} dP(s)$ is the probability that $O\in [u,\tau]$ under $P$. Estimation of the sampling distribution of 
$n^{1/2}(\bar{P}_n-\bar{P}_0)$ with the bootstrap distribution $n^{1/2}(\bar{P}_n^{\#}-\bar{P}_n)$ results then in $Z_n^{+,\#}$ and  corresponding finite sample conservative $0.95$-confidence interval, which can also be used for $\Psi(P_n^*)$.
 In the previous example, this confidence interval appeared to be much too conservative to be practically useful, but in this example, since our bound on $R_{20}(P,P_0)$ is sharp and the sectional variation norm bounds $M_1=M_3=2C^u$ are easily determined in terms of the sectional variation norm bound $C^u$ of the model, this appears to be an interesting finite sample conservative confidence interval. We propose to apply this confidence interval to $C^u=C_n$ and the corresponding bounds $M_{1n}=M_{3n}=2C_n$.
 
 

{\bf Asymptotic validity of the nonparametric bootstrap for the HAL-MLE Theorem \ref{thnpbootmle}:}
As remarked in the previous example, computation of the HAL-MLE  $p_n^{\#}=\arg\min_{\pl p\pl_v^*\leq C,p\ll p_n}P_n^{\#}L(p)$ is much faster than the computation of $p_n=\arg\min_{\pl p\pl_v^*\leq C^u}P_nL(p)$, due to only having to minimize the empirical risk over the bootstrap sample over the linear combinations of indicator functions that had non-zero coefficients in $p_n$. 
The conditions for our asymptotic consistency of the nonparametric bootstrap of Theorem \ref{thnpbootmle} hold, as we show now. 
We have $P_0 L(P)-P_0L(P_0)=\int (p-p_0)^2 d\mu$ and thus $P_0 R_{2,L}(P,P_0)=\int (p-p_0)^2 d\mu$.
 It easily follows that  $P_0(L(P)-L(P_0))^2$ can be bounded by $P_0 R_{2,L}(P,P_0)$.
We now have to establish the second-order exact expansion for  $d_n(P_n^{\#*},P_n^*)=P_n L(P_n^{\#*})-P_n L(P_n^*)$.
We have
\begin{eqnarray*}
d_n(P_n^{\#*},P_n^*)&=&P_n\{ L(P_n^{\#*})-L(P_n^*)\}\\
&=&P_n\{ -2(p_n^{\#*}-p_n^*)+\int p_n^{\#*2} d\mu-\int p_n^{*2} d\mu\\
&=&P_n\{ -2(p_n^{\#*}-p_n^*)+\int (p_n^{\#*}-p_n^*)(p_n^{\#*}+p_n^*) d\mu  \\
&=&P_n\{ (-2+2p_n^*) (p_n^{\#*}-p_n^*)d\mu  \\
&=&\int (p_n^{\#*}-p_n^*)^2 d\mu  \\
\end{eqnarray*}
Thus, $P_n R_{2,L}(P_n^{\#*},P_n^*)=\int (p_n^{\#*}-p_n^*)^2 d\mu$. Clearly, $P_n (L(P_n^{\#*})-L(P_n^*))^2$ can be bounded by $P_n R_{2,L}(P_n^{\#*},P_n^*)$.
Application of Theorem \ref{thnpbootmle} now shows that $\int (p_n^{\#}-p_n)^2 d\mu=O_P(n^{-1/2-\alpha(d)})$.

{\bf Preservation of rate for HAL-TMLE under sampling from $P_n$:}
We can carry out the general proof in the Appendix \ref{AppendixC} establishing that $\epsilon_n^{\#}=O_P(n^{-1/4-\alpha(d)/2})$.

{\bf Asymptotic consistency of the bootstrap for the HAL-TMLE:}
This verifies all conditions of Theorem \ref{thnpboothaltmle} which establishes the asymptotic efficiency and asymptotic consistency of the nonparametric bootstrap.
\begin{theorem}
Consider the model ${\cal M}$ defined by upper and lower bound $M<\infty$, $\delta\geq 0$ on the densities on support $[0,\tau]$, and the assumption that the sectional variation norm of the densities over $[0,\tau]$ is bounded by $C^u<\infty$. 

We have that $\Psi(P_n^*)$ is asymptotically efficient, i.e. $n^{1/2}(\Psi(P_n^*)-\Psi(P_0))\Rightarrow_d N(0,\sigma^2_0)$, where
$\sigma^2_0=P_0 \{D^*(P_0)\}^2$.

In addition, conditional on $(P_n:n\geq 1)$, $Z_n^{1,\#}=n^{1/2}(\Psi(P_n^{\#*})-\Psi(P_n^*))\Rightarrow_d N(0,\sigma^2_0)$.

This theorem can also be applied to the setting in which $C^u=C_n$.
\end{theorem}

{\bf Asymptotic consistency of the bootstrap for the exact second-order expansion of HAL-TMLE:}
We have $n^{1/2}(\Psi(P_n^*)-\Psi(P_0))=n^{1/2}(P_n-P_0)D^*(P_n^*)-n^{1/2}\int (p_n^*-p_0)^2 d\mu$.
Let  $Z_n^{2,\#}$ be the nonparametric bootstrap estimator of this exact second-order expansion:
\[
Z_n^{2,\#}=n^{1/2}(P_n^{\#}-P_n)D^*(P_n^{\#*})-n^{1/2}\int (p_n^{\#*}-p_n^*)^2 d\mu.\]
In addition, consider the upper bound:
\[
\mid n^{1/2}(\Psi(P_n^*)-\Psi(P_0))\mid\leq \mid n^{1/2}(P_n-P_0)D^*(P_n^*)\mid +\int (p_n^*-p_0)^2 d\mu .\]
Let $Z_n^{3,\#}$ be the nonparametric bootstrap estimator of this conservative sampling distribution:
\[
Z_n^{3,\#}\equiv \mid n^{1/2}(P_n^{\#}-P_n)D^*(P_n^{\#*})\mid +n^{1/2}\int (p_n^{\#*}-p_n^*)^2 d\mu .\]
Application of Theorems \ref{thnpbootexactexp} and \ref{thupperbound} prove the asymptotic consistency of these two nonparametric bootstrap distributions.

\begin{theorem}
We have that conditional on $(P_n:n\geq 1)$, $Z_n^{j,\#}\Rightarrow_d N(0,\sigma^2_0)$ as $n\rightarrow\infty$, $j=2,3$. 
 As a consequence, an $0.95$-confidence interval for $\psi_0$ based on $Z_n^{2,\#}$ and $Z_n^{3,\#}$ have asymptotic coverage 0.95 of $\psi_0$.
 This theorem can also be applied to setting in which $C^u=C_n$.
 \end{theorem}
 
 Finally, we remark that our HAL-MLE is really indexed by the hypothesized model bounds $(C^u,\delta,M)$ and these might all three be unknown to the user. So in that case, we recommend to select
 all three with the cross-validation selector $(C_n,\delta_n,M_n)$ and define the HAL TMLE and bootstrap of the HAL-TMLE at this fixed choice $(C_n,\delta_n,M_n)$.




%

\section{Discussion}\label{sectdisc}
In parametric models and, more generally, in  models small enough so that  the MLE is still well behaved, one can use the nonparametric bootstrap to estimate the sampling distribution of the MLE. It is generally understood that in these small models the nonparametric bootstrap  outperforms estimating the sampling distribution with a normal distribution (e.g., with variance estimated as the sample variance of the influence curve of the MLE), by picking up the higher order behavior of the MLE, {\em if asymptotics has not set in yet}. In such small models, reasonable sample sizes already achieve the normal approximation in which case the Wald type confidence intervals will perform well. 
Generally speaking, the nonparametric bootstrap is a valid method when the estimator is a compactly differentiable function of the empirical measure, such as the Kaplan-Meier estimator (i.e., one can apply the functional delta-method to analyze such estimators) \citep{Gill89,vanderVaart&Wellner96}. These are estimators that essentially do not use  smoothing of any sort.

On the other hand, efficient estimation of a pathwise differentiable target parameter in  large realistic models generally requires estimation of the data density, and thereby    machine learning such as super-learning to estimate the relevant parts of the data distribution. Therefore, efficient one-step estimators or TMLEs are not compactly differentiable functions of the data distribution.
Due to this reason, we moved away from using the nonparametric bootstrap to estimate its sampling distribution, since it represents a generally inconsistent method (e.g., a cross-validation selector behaves very differently under sampling from the empirical distribution than under sampling from the true data distribution). Instead we estimated the normal limit distribution by estimating the variance of the influence curve of the estimator.  

Such an influence curve based method is asymptotically consistent and therefore results in asymptotically valid $0.95$-confidence intervals. However,
in such large models the nuisance parameter estimators will converge at low rates (like $n^{-1/4}$ or lower) with large constants depending on the size of the model, so that for normal sample sizes the exact second-order remainder could be easily larger than the leading empirical process term with its normal limit distribution. 
So one has to pay a significant  price for using the computationally attractive influence curve based confidence intervals, by generally reporting overly optimistic confidence intervals. That is, for small models the bootstrap is available but not that important since estimators will quickly achieve asymptotics, while in large models it appears to 
not be available even though it is  crucial since estimators generally achieve asymptotics for very large sample sizes.  One might argue that one should use a smooth bootstrap instead by sampling from an  estimator of the density of the data distribution. General results show that such a smooth bootstrap method will be asymptotically valid as long as the density estimator is consistent. This is like carrying out a simulation study for the estimator in question using an estimator of the true data distribution as sampling distribution. However, estimation of the actual density of the data distribution is itself a very hard problem, with bias heavily affected by the curse of dimensionality, and, in addition, it can be immensely burdensome to construct such a density estimator and sample from it when the data is complex and high dimensional.

As demonstrated in this article, the HAL-MLE provided a solution to this bottleneck. The HAL-MLE($C^u$)  of the nuisance parameter is an actual MLE minimizing the empirical risk over a highly nonparametric parameter space (depending on the model ${\cal M}$) in which it is assumed that the sectional variation norm of the nuisance parameter is bounded by universal constant $C^u$. This MLE is still well behaved by being consistent at a rate that is in the worst case still faster than $n^{-1/4}$. However, this MLE is not an interior MLE, but will be on the edge of its parameter space: the MLE will itself have sectional variation norm equal to the maximal allowed value $C^u$. Nonetheless, our analysis shows that it is still  a smooth enough function of the data (while not being compactly differentiable at all) that it is equally well behaved under sampling from the empirical distribution.  

As a consequence of this robust behavior of the HAL-MLE, for  models in which the nuisance parameters of interest are cadlag functions with a universally bounded sectional variation norm (beyond possible other assumptions), we presented asymptotically consistent estimators of the sampling distribution of the HAL-TMLE and HAL-one step estimator of the target parameter of interest using the nonparametric bootstrap. Our proposals range from a  bootstrap estimator of the HAL-TMLE itself, a  bootstrap estimator of the exact second-order expansion of the HAL-TMLE, and  two bootstrap estimators of  conservative upper bounds on the exact second-order expansion of the HAL-TMLE. In addition, we presented a highly conservative finite sample sampling distribution based on applying general  integration by parts formulas to the leading empirical process term and a conservative bound of the exact second-order remainder. We also provided slight variations of these proposals, corresponding with the HAL-one step estimator.

Our estimators of the sampling distribution are highly sensitive to the curse of dimensionality, just as the sampling distribution of the HAL-TMLE itself: specifically, the HAL-MLE on a bootstrap sample will  converge just as slowly to its truth as under sampling from the true distribution. Therefore, in high dimensional estimation problems, we  expect highly significant gains in valid inference relative to Wald type confidence intervals  that are purely  based on the normal limit distribution of the HAL-TMLE.

In general, the  user will typically not know how to select the upper bound $C^u$ on the sectional variation norm of the nuisance parameters (except if the nuisance  parameters are cumulative distribution functions). Therefore, we recommend to select this bound with cross-validation just as we use cross-validation to select the sectional variation norm bound in the HAL-MLE. Due to the oracle inequality for the cross-validation selector $C_n$ (which only relies on a bound on the supremum norm of the loss function), the data adaptively selected upper bound will be selected larger than the true sectional variation norm  $C_0$ of the nuisance parameters $(Q_0,G_0$, as sample size increases. 
Therefore, our bootstrap estimators will still be guaranteed to be consistent for its normal limit distribution while incorporating its higher order behavior.

For small sample sizes, one would most likely select a bound smaller than the sectional variation norm of the true nuisance parameter (optimally trading off bias and variance of the HAL-MLE), and as sample size increases it will get larger and larger till at some large enough sample size it will plateau, having reached the sectional variation norm of the true nuisance parameter. The advantage of this data adaptive choice of the bound is that our resulting bootstrap inference will  have adapted to the true underlying sectional variation norm once it has reached that plateaux. The disadvantage is that for small sample size the selected  model (implied by the selected bound) might be smaller than a model containing the true data distribution, so that our bootstrap methods might still be optimistic for such small sample sizes. Nonetheless, it will evaluate the finite sample sampling distribution of the HAL-TMLE (relative to its truth under $P_n$ satisfying this same bound)  in a correctly specified model (just too small model). In particular, it still has the same first order behavior as the sampling distribution of the actual HAL-TMLE (using cross-validation to select the bound), but it may underestimate  its higher order behavior.  In these settings there is still use for our conservative sampling distributions whose conservative nature might outweigh the potential underestimation of uncertainty due to a selecting a  bound smaller than the true sectional variation norm of the nuisance parameter. Simulations will likely shed light on this.

The sectional variation norm plays a fundamental role in this work. The sectional variation norm of a function can be interpreted as a measure of complexity or smoothness of the function: it represents the sum of the absolute value of the coefficients in our integral presentation (\ref{Frepresentation}) of the function as an infinite linear combination of indicator basis functions.
The sectional variation norm of the true nuisance parameter such as a regression function can be viewed as a  general measure of degree of sparsity.
For example, a regression function of $d$ variables that is a sum of  functions of maximally three variables  has a sectional variation norm that behaves as $d^3$ instead of the worst case behavior $2^d$.  
This definition of degree of sparsity (i.e., the true function has a certain sectional variation norm) is not dependent on a choice of a main term regression  model as in the typical Lasso literature (far from a saturated model). The HAL-MLE and HAL-TMLE using cross-validation to select the sectional variation norm bounds will adapt to this underlying sparsity and so will our inference for the target parameter using this selected bound as fixed in the bootstrap of the corresponding HAL-TMLE. This demonstrates the enormous importance of this measure of sparsity for the behavior of the cross-validated HAL-TMLE (and HAL-MLE).

Presumably, by selecting another basis and corresponding function representation, one could also define sparsity as the sum of the absolute value of the coefficients of these basis functions in its representation. 
Possible advantages of the indicator basis and its representation (\ref{Frepresentation}) is that it allows approximation of discontinuous functions; it it easy to determine the subset of indicator basis functions that are relevant for the given sample, making the implementation of HAL-MLE doable; and the HAL-MLE has good convergence properties, and is highly robust (as shown by our bootstrap results), which might not be available for many other basis choices. The latter appears to be due to the feature of the collection of indicator basis function in that it represents a Donsker class, while many other choices of basis functions cannot be embedded in a Donsker class (e.g. Fourrier series). Therefore, we wonder if the indicator basis and its function representation (\ref{Frepresentation})  is a particular powerful (and possibly unique) choice for defining a measure of complexity of a true parameter, bounding the model accordingly, defining an MLE of the nuisance parameters for such a model, selecting the bound with cross-validation, and using the nonparametric bootstrap to estimate the sampling distribution of its corresponding TMLE treating the selected bound as fixed, for the sake of inference as carried out in this article.

This article focused on a HAL-TMLE that represents the statistical  target parameter  $\Psi(P)$ as a   function  $\Psi(Q_1(P),\ldots,Q_{K_1}(P))$ of variation independent nuisance parameters $(Q_1,\ldots,Q_{K_1})$. In some examples  it has important advantages to represent  $\Psi(P)$ in terms of recursively defined nuisance parameters. For example, the longitudinal  one-step TMLE  of causal effects of multiple time point interventions  in \citep{Gruber&vanderLaan12,Petersen&Schwab&vanderLaan13} relies on a sequential   regression representation of the target parameter \citep{Bang&Robins05}.  In this case,  the next regression is defined as the regression of the previous regression on a shrinking history, across a number of regressions, one for each time point at which an intervention takes place. By fitting each of these sequential regressions with an HAL-regression we obtain the analogue of the HAL-TMLE for this sequential regression type TMLE. Our convergence results for the HAL-MLE and the bootstrapped HAL-MLE can be applied to these HAL-MLEs of each regression, in which case the outcome is the HAL-MLE fit of the previous regression. Some additional work will be needed to deal with the dependence on the previous regression to analyze this type of sequential HAL-TMLE, but we conjecture that the nonparametric bootstrap will be valid for this type of non-variation independent HAL-TMLE as well.

\subsection*{Acknowledgement.}
This research is funded by NIH-grant 
5R01AI074345-07.


\bibliography{combined,TLB2}


\appendix 
\section*{Appendix.}


\section{Proof that the one-step TMLE $Q_n^*$ preserves rate of convergence of $Q_n$}\label{AppendixA}
The following lemma establishes that the one-step TMLE $Q_n^*=Q_{n,\epsilon_n}$ preserves the rate of convergence of $Q_n$, where $Q_{\epsilon}$ is a univariate local least favorable submodel through $Q$ at $\epsilon =0$. 
Recall the notation $L_1(Q_1,Q_2)=L_1(Q_1)-L_1(Q_2)$.
\begin{lemma}\label{lemmathalmle}
Assume $\epsilon\rightarrow L_1(Q_{n,\epsilon})$ be  differentiable with a uniformly bounded derivative on an interval $\epsilon \in (-\delta,\delta)$ for some $\delta>0$.
Define $\epsilon_{0n}=\arg\min_{\epsilon}P_0 L_1(Q_{n,\epsilon})$, $\epsilon_n=\arg\min_{\epsilon}P_n L_1(Q_{n,\epsilon})$, and assume the weak regularity condition  
\begin{equation}\label{Qnstarconv}
\epsilon_{0n}=O_P(d_{01}^{1/2}(Q_n,Q_0)).\end{equation}
Then,
\begin{equation}\label{thalmle}
d_{01}(Q_n^*,Q_0)\leq d_{01}(Q_n,Q_0)+O_P(n^{-1/2-\alpha(d)}).\end{equation}
Specifically,
\begin{eqnarray*}
d_{01}(Q_n^*,Q_0)&\leq& d_{01}(Q_n,Q_0)
-(P_n-P_0)L_1(Q_{n,\epsilon_{0n}},Q_0) .
\end{eqnarray*}
\end{lemma}
This also proves that the $K$-th step TMLE using a {\em finite} $K$ (uniform in $n$) number of iterations satisfies $d_{01}(Q_n^*,Q_0)\leq d_{01}(Q_n,Q_0)+O_P(n^{-1/2-\alpha(d)})$. So  if $r_n=P_n D^*(Q_{n,\epsilon_n},G_n)$ is not yet $o_P(n^{-1/2})$, then  one should consider  a $K$-th step  TMLE to guarantee that $r_n$ is small enough to be neglected (we know that the fully iterated TMLE will solve $P_n D^*(Q_n^*,G_n)=0$, but this one is harder to analyze).
\newline
{\bf Proof of Lemma \ref{lemmathalmle}:}
Using the MLE properties of $\epsilon_{0n}$ and $Q_n$, subsequently, we obtain
\begin{eqnarray*}
P_0L_1(Q_n^*)-P_0L_1(Q_0)&=&\{ P_0L_1(Q_{n,\epsilon_n})-P_0L_1(Q_{n,\epsilon_{0n}})\}
+\{P_0L_1(Q_{n,\epsilon_{0n}})-P_0L_1(Q_n)\}\\
&&+\{P_0L_1(Q_n)-P_0L_1(Q_0)\}\\
&\leq & P_0L_1(Q_{n,\epsilon_{0n}},Q_n) +P_0L_1(Q_n,Q_0)\\
&=&(P_0-P_n)L_1(Q_{n,\epsilon_{0n}},Q_n) +P_n L_1(Q_{n,\epsilon_{0n}},Q_n)+d_{01}(Q_n,Q_0)\\
&\leq &-(P_n-P_0)L_1(Q_{n,\epsilon_{0n}},Q_n)+d_{01}(Q_n,Q_0)\\
\end{eqnarray*}
An exact first order Tailor expansion $f(\epsilon_{0n})-f(0)=\frac{d}{d\xi}f(\xi)\epsilon_{0n}$ for a $\xi$ between $0$ and $\epsilon_{0n}$ applied to $f(\epsilon)=L_1(Q_{n,\epsilon})$ yields:
\[
L_1(Q_{n,\epsilon_{0n} })-L_1(Q_{n})=\epsilon_{0n}\frac{d}{d\xi }L_1(Q_{n,\xi} ).\]
Using that $\epsilon_{0n}$ converges as fast to zero as $d_{01}(Q_n,Q_0)$ and that $(P_n-P_0)L_1(Q_{n,\epsilon_{0n}},Q_0)$ is an evaluation of an  empirical process indexed by  the Donsker class of cadlag functions with sectional variation norm bounded by universal $M_1$, we obtain $(P_n-P_0)L_1(Q_{n,\epsilon_{0n}},Q_n)=O_P(n^{-1/2-\alpha(d)})$, analogue to the proof of this rate of convergence for the  HAL-MLE. So we have
\[
d_{01}(Q_n^*,Q_0)\leq d_{01}(Q_n,Q_0)+O_P(n^{-1/2-\alpha(d)}) .\]
$\Box$ 

\section{Asymptotic convergence of bootstrapped HAL-MLE: Proof of Theorem \ref{thnpbootmle}.}\label{AppendixB}

Consider  a given set $A$ of possible values for $(s,u_s)$, where $s$ is a subset of $\{1,\ldots,m\}$ and $u_s$ is a value in $[0_s,\tau_s]$.  
Let $D_{m,C^l,C^u}[0,\tau]=\{f\in D[0,\tau]: C_l\leq \pl f\pl_v^*\leq C^u\}$ be the set of $m$-variate real valued cadlag functions on $[0,\tau]$ with sectional variation norm between $C^l$ and $C^u$. Let $C^l<C^u$. Consider the case that the parameter space $Q({\cal M})$ of $Q(P)=\arg\min_{Q\in Q({\cal M})}PL_1(Q)$ is given by
 \begin{equation}\label{B1}
{\cal F}^{np}_A\equiv \{F\in D_{m,C^l,C^u}[0,\tau]: dF_s(u_s)=I_{(s,u_s)\in A}dF_s(u_s)\}.\end{equation}
Here we use the short-hand notation $g(x)=I_{x\in A}g(x)$ to state that $g(x)=0$ for $x\not \in A$.
Recall that each function $F\in D_{m,C^l,C^u}[0,\tau]$ can be represented as $F(x)=F(0)+\sum_{s\subset \{1,\ldots,m\}}\int_{(0_s,x_s]} dF_s(u_s)$, while 
$\pl F\pl_v^*=\mid F(0)\mid +\sum_{s\subset\{1,\ldots,m\}}\int_{(0_s,\tau_s]}\mid dF_s(u_s)\mid$. 
We also consider the case that the parameter space $Q({\cal M})$ equals 
\begin{equation}\label{B2}
{\cal F}^{np+}_A\equiv \{F\in D_{m,C^l,C^u}[0,\tau]: dF_s(u_s)=I_{(s,u_s)\in A}dF_s(u_s),dF_s(u_s)\geq 0,F(0)\geq 0\}.\end{equation}
For this case, we allow that $C^l=C^u$.
If we select $A$ is unrestricted, then ${\cal F}^{np}_A$ is the set of all $m$-variate cadlag functions with sectional variation norm bounded by $(C^l,C^u)$, and 
${\cal F}^{np+}_A$ is the set of non-negative monotone functions  with positive total mass between $C^l$ and $C^u$  (i.e.,we obtain a set of  cumulative distributions if $C^l=C^u=1$).
So the set $A$ indicates that $dF_s(u_s)$ is only non-zero for $(s,u_s)\in A$ and thereby restricts the set of possible functions in the parameter space.

In our model defined in Section 2 we assumed that each component of $Q$ and $G$ has its own loss and such a parameter space, and the HAL-MLE can be computed separately for each of these components. Thus, by applying our results below to each component of $Q$ and $G$ separately with $L_1(Q)$ replaced by its loss function, $Q_n$ and $Q_n^{\#}$ replaced by its HAL-MLE and  HAL-MLE applied to bootstrap sample, respectively, one obtains the desired result for each component.
 
 The next Theorem \ref{ThB} shows that $d_{n1}(Q_n^{\#},Q_n)=P_n \{L_1(Q_n^{\#})-L_1(Q_n)\}$ converges at rate $n^{-1/2-\alpha(d)}$ and that this empirical loss-based dissimilarity $d_{n1}(Q_n^{\#},Q_n)$ dominates a quadratic dissimilarity (making it equally powerful as $d_{01}(Q_n,Q_0)$).
 
\begin{theorem}\label{ThB}
Recall Definition \ref{defabscont} of $Q\ll Q_n$. 
Let $Q_n=\arg\min_{Q\in Q({\cal M}}P_nL_1(Q)$, $Q_n^{\#}=\arg\min_{Q\in Q({\cal M}),Q\ll Q_n}P_n^{\#}L_1(Q)$ be the HAL-MLE and HAL-MLE on the bootstrap sample, respectively, where $Q({\cal M})$ either equals ${\cal F}_A^{np}$ (\ref{B1}) or ${\cal F}_A^{np+}$ (\ref{B2}). 
Either assume $\pl Q_n\pl_v^*=C^u$ or assume that $\pl Q_n^{\#}\pl_v^*\leq \pl Q_n\pl_v^*$ (and $\pl Q_n\pl_v^*>C^l$ if $C^l<C^u$), conditional on $(P_n:n\geq 1)$.

Let $P_n R_{2L_1,n}(Q_n^{\#},Q_n)$ be defined by
\[
P_n \{L_1(Q_n^{\#})-L_1(Q_n)\}=P_n \frac{d}{dQ_n}L_1(Q_n)(Q_n^{\#}-Q_n)+P_n R_{2L_1,n}(Q_n^{\#},Q_n).\]

We have $P_n\frac{d}{dQ_n}L_1(Q_n)(Q_n^{\#}-Q_n)\geq 0$ so that  
\[
d_{n1}(Q_n^{\#},Q_n) \geq P_n R_{2L_1,n}(Q_n^{\#},Q_n).\]
Thus, if  
\[
P_n \{L_1(Q_n^{\#})-L_1(Q_n)\}^2\lesssim P_nR_{2L_1,n}(Q_n^{\#},Q_n),\]
then 
\begin{equation}\label{b3}
P_n \{L_1(Q_n^{\#})-L_1(Q_n)\}^2\lesssim d_{n1}(Q_n^{\#},Q_n).\end{equation}
Suppose the latter (\ref{b3}) holds. 
By Lemma \ref{generalhalmle} below, then
\[ d_{n1}(Q_n^{\#},Q_n)=O_P(n^{-1/2-\alpha(d)}).\]
This, on its turn then  implies
 $P_nR_{2L_1,n}(Q_n^{\#},Q_n)=O_P(n^{-1/2-\alpha(d)})$.
\end{theorem}

In order to provide the reader a concrete example, we provide here the corollary for the squared error loss.
\begin{corollary}
Consider definition and assumptions of Theorem \ref{ThB}.
Suppose that $L_1(Q)(O)=(Y-Q(X))^2$ is the squared error loss. 
Then,
we have
\[
d_{n1}(Q_n^{\#},Q_n)\geq P_n (Q_n^{\#}-Q_n)^2 .\]
Since  $P_n  \{L_1(Q_n^{\#})-L_1(Q_n)\}^2\lesssim P_n (Q_n^{\#}-Q_n)^2$, this implies
$P_n  \{L_1(Q_n^{\#})-L_1(Q_n)\}^2\lesssim d_{n1}(Q_n^{\#},Q_n)$.
By Lemma \ref{generalhalmle}, this shows 
\[
d_{n1}(Q_n^{\#},Q_n)=O_P(n^{-1/2-\alpha(d)}).
\]
This, on its turn then implies  $P_n (Q_n^{\#}-Q_n)^2=O_P(n^{-1/2-\alpha(d)})$.
\end{corollary}

 
To prove Theorem \ref{ThB}, we first present the following straightforward lemma, which follows immediately by just imitating the proof of the convergence of the HAL-MLE $Q_n$ itself but now under sampling from $P_n$. Here one uses that $P_n^{\#}L_1(Q_n^{\#})-P_n^{\#}L_1(Q_n)\leq 0$, since $Q_n$ is an element of the parameter space over which $Q_n^{\#}$ minimizes $P_n^{\#}L_1(Q)$.

\begin{lemma}\label{generalhalmle}
Consider the above setting.
We have
\begin{eqnarray}
0&\leq & d_{n1}(Q_n^{\#},Q_n)\equiv P_n\{ L_1(Q_n^{\#})-L_1(Q_n)\}\nonumber \\
&=& -(P_n^{\#}-P_n)\{L_1(Q_n^{\#})-L_1(Q_n)\}+P_n^{\#}\{L_1(Q_n^{\#})-L_1(Q_n)\}\nonumber \\
&\leq& -(P_n^{\#}-P_n)\{L_1(Q_n^{\#})-L_1(Q_n)\}.\label{boota}
\end{eqnarray}
As a consequence, by empirical process theory \citep{vanderVaart&Wellner11}, we have $d_{n1}(Q_n^{\#},Q_n)=P_n L_1(Q_n^{\#})-P_nL_1(Q_n)$ is $O_P(n^{-1/2})$, and
if $\pl L_1(Q_n^{\#})-L_1(Q_n)\pl_{P_n}^2\lesssim d_{n1}(Q_n^{\#},Q_n)$, then we have 
$d_{n1}(Q_n^{\#},Q_n)=O_P(n^{-1/2-\alpha(d)})$.
\end{lemma}

{\bf Proof of Theorem \ref{ThB}:}
We first prove the results for the general loss function, and subsequently, we will consider the special case that $L_1(Q)$ is the squared error loss.
Consider the  $h$-specific path  \[
Q_{n,\epsilon}^h(x)=(1+\epsilon h(0))Q_n(0)+\sum_s \int_{(0_s,x_s]}(1+\epsilon h_s(u_s)) dQ_{n,s}(u_s))\]
 for $\epsilon \in [0,\delta)$ for some $\delta>0$, where $h$ is uniformly bounded, and,  if $C^l<C^u$,  \[
 r(h,Q_n)\equiv  h(0)\mid Q_n(0)\mid+\sum_s \int_{(0_s,\tau_s]} h_s(u_s)\mid dQ_{n,s}(u_s)\mid \leq 0,\] while if $C^l=C^u$, then $r(h,Q_n)=0$.  Let ${\cal H}=\{h: r(h,Q_n)\leq 0,\pl h\pl_{\infty}<\infty\}$ be the set of possible functions $h$(i.e., functions of $s,u_s$), which defines a collection of paths $\{Q_{n,\epsilon}^h:\epsilon\}$ indexed by $h\in {\cal H}$. Consider a given $h\in {\cal H}$ and let's denote this path with $Q_{n,\epsilon}$, suppressing the dependence on $h$ in the notation.
For $\epsilon\geq 0$ small enough we have $(1+\epsilon h(0))>0$ and $1+\epsilon h_s(u_s)>0$. Thus, for $\epsilon\geq $ small enough we have
\begin{eqnarray*}
\pl Q_{n,\epsilon}\pl_v^*&=&(1+\epsilon h(0))\mid Q_n(0)\mid+\sum_s \int_{(0_s,\tau_s]}(1+\epsilon h_s(u_s))\mid dQ_{n,s}(u_s)\mid \\
&=& \pl Q_n\pl_v^*+\epsilon\left\{ h(0)\mid Q_n(0)\mid+\sum_s \int_{(0_s,\tau_s]}h_s(u_s) \mid dQ_{n,s}(u_s)\mid \right\} \\
&=&\pl Q_n\pl_v^*+\epsilon r(h,Q_n)\\
&\leq & \pl Q_n\pl_v^*,
\end{eqnarray*}
by assumption that $r(h,Q_n)\leq 0$. If $C^l=C^u$ and thus $r(h,Q_n)=0$, then the above shows $\pl Q_{n,\epsilon}\pl_v^*=\pl Q_n\pl_v^*$.
Thus, for a small enough $\delta>0$ $\{Q_{n,\epsilon}:0\leq \epsilon<\delta\}$ represents a path of cadlag functions with sectional variation norm  bounded from below by $C^l$ and smaller than or equal to $\pl Q_n\pl_v^*\leq C^u$. In addition, we have that $dQ_{n,s}(u_s)=0$ implies $(1+\epsilon h_s(u_s))dQ_{n,s}(u_s)=0$ so that the support of $Q_{n,\epsilon}$ is included in the support $A$ of $Q_n$ as defined by ${\cal F}_A^{np}$. Thus, this proves that for $\delta>0$ small enough this path $\{Q_{n,\epsilon}:0\leq \epsilon\leq \delta\}$ is indeed a  submodel of the parameter space of $Q$, defined as ${\cal F}_{A}^{np}$ or ${\cal F}_A^{np+}$.

We also have that
 \[
 Q_{n,\epsilon}-Q_n=\epsilon\left\{ Q_n(0)h(0)+ \sum_s \int_{(0_s,x_s]} h_s(u_s)dQ_{n,s}(u_s)\right\} .\]
 Thus, this path generates a direction $f(h,Q_n)$  at $\epsilon=0$ given by:
 \[
 \frac{d}{d\epsilon}Q_{n,\epsilon}=f(h,Q_n)\equiv Q_n(0)h(0)+ \sum_s \int_{(0_s,x_s]} h_s(u_s)dQ_{n,s}(u_s) .\]
 Let ${\cal S}\equiv \{f(h,Q_n): h\in {\cal H}\}$ be the collection of directions generated by our family of paths.
 By definition of the MLE $Q_n$, we also have that $\epsilon \rightarrow P_n L_1(Q_{n,\epsilon})$ is minimal over $[0,\delta)$ at $\epsilon =0$. 
This shows that the derivative of $P_n L_1(Q_{n,\epsilon})$ from the right at $\epsilon =0$ is  non-negative:
 \[
 \frac{d}{d\epsilon^+ }P_n L_1(Q_{n,\epsilon })\geq 0\mbox{ at $\epsilon =0$}.\]
 This derivative is given by $P_n \frac{d}{dQ_n}L_1(Q_n)(f(h,Q_n))$, where $d/dQ_nL_1(Q_n)(f(h,Q_n))$ is the directional (Gateaux) derivative of $Q\rightarrow L_1(Q)$ at $Q_n$ in  in direction $f(h,Q_n)$.
  Thus for each $h\in {\cal H}$, we have
\[
P_n \frac{d}{dQ_n}L_1(Q_n)(f(h,Q_n)) \geq 0 .\]
Suppose that
\begin{equation}\label{keya1}
Q_n^{\#}-Q_n\in {\cal S}=\{f(h,Q_n):h\in {\cal H}\}.\end{equation}
Then, we have
\[
P_n \frac{d}{dQ_n}L_1(Q_n)(Q_n^{\#}-Q_n)\geq 0.\]
Combined with the stated second-order Tailor expansion of $P_n L_1(Q)$ at $Q=Q_n$ with exact second-order remainder  $P_nR_{2L_1,n}(Q_n^{\#},Q_n)$, this proves
\[
P_n\{L_1(Q_n^{\#})-L_1(Q_n)\}\geq P_nR_{2L_1,n}(Q_n^{\#},Q_n).\]
Thus it remains to show (\ref{keya1}).

In order to prove (\ref{keya1}), let's solve explicitly for $h$ so that $Q_n^{\#}-Q_n=f(h,Q_n)$ and then verify that $h\in {\cal H}$ satisfies its assumed constraints (i.e., $r(h,Q_n)\leq 0$ if $C^l<C^u$ or $r(h,Q_n)=0$ if $C^l=C^u$,  and $h$ is uniformly bounded).
We have
\begin{eqnarray*}
Q_n^{\#}-Q_n&=&Q_n^{\#}(0)-Q_n(0)+\sum_s\int_{(0_s,x_s]} d(Q_{n,s}^{\#}-dQ_{n,s})(u_s)\\
&=& Q_n^{\#}(0)-Q_n(0)+\sum_s \int_{(0_s,x_s]} \frac{d(Q_{n,s}^{\#}-dQ_{n,s})}{dQ_{n,s}} dQ_{n,s}(u_s),
\end{eqnarray*}
where we used that $Q_{n,s}^{\#}\ll Q_{n,s}$ for each subset $s$.
Let $h(Q_n^{\#},Q_n)$ be defined by
\begin{eqnarray*}
h(Q_n^{\#},Q_n)(0)&=&(Q_n^{\#}(0)-Q_n(0))/Q_n(0)\\
h_s(Q_n^{\#},Q_n)&=&\frac{d(Q_{n,s}^{\#}-dQ_{n,s})}{dQ_{n,s}}\mbox{ for all subsets $s$} .
\end{eqnarray*}
For this choice $h(Q_n^{\#},Q_n)$, we have $f(h,Q_n)=Q_n^{\#}-Q_n$.
First, consider the case $Q({\cal M})={\cal F}_A^{np}$ or $Q({\cal M})={\cal F}_A^{np+}$, but$C^l<C^u$.
We now need to verify if  $r(h,Q_n)\leq 0$ for this choice $h=h(Q_n^{\#},Q_n)$.
We have
\begin{eqnarray*}
r(h,Q_n)&=&\frac{Q_n^{\#}(0)-Q_n(0)}{Q_n(0)}\mid Q_n(0)\mid+\sum_s\int_{(0_s,\tau_s]}\frac{dQ_{n,s}^{\#}-dQ_{n,s}}{dQ_{n,s}}\mid dQ_{n,s}\mid\\
&=&I(Q_n(0)>0)\{Q_n^{\#}(0)-Q_n(0)\}+I(Q_n(0)\leq 0)\{Q_n(0)-Q_n^{\#}(0)\}\\
&&+\sum_s\int_{(0_s,\tau_s]}I(dQ_{n,s}\geq 0) d(Q_{n,s}^{\#}-dQ_{n,s})\\
&&+\sum_s\int_{(0_s,\tau_s]} I(dQ_{n,s}<0)d(Q_{n,s}-Q_{n,s}^{\#})\\
&=&-\pl Q_n\pl_v^*+Q_n^{\#}(0)\{ I(Q_n(0)>0)-I(Q_n(0)\leq 0)\} \\
&&+\sum_s \int_{(0_s,\tau_s]}\{I(dQ_{n,s}\geq 0) -I(dQ_{n,s}\leq 0\} dQ_{n,s}^{\#}\\
&\leq&-\pl Q_n\pl_v^*+\mid Q_n^{\#}(0)\mid +\sum_s\int_{(0_s,\tau_s]} \mid dQ_{n,s}^{\#}(u_s)\mid \\
&=&-\pl Q_n\pl_v^*+\pl Q_n^{\#}\pl_v^*\\
&\leq &0,
\end{eqnarray*}
since $\pl Q_n^{\#}\pl_v^*\leq \pl Q_n\pl_v^*$, by assumption.
Thus, this proves that indeed $r(h,Q_n)\leq 0$ and thus that $Q_n^{\#}-Q_n\in {\cal S}$.
Consider now the case that $Q({\cal M})={\cal F}_A^{np+}$ {\em and } $C^l=C^u$. Then $\pl Q_n\pl_v^*=\pl Q_n^{\#}\pl_v^*=C^u$. We now need to show that $r(h,Q_n)=0$ for this choice $h=h(Q_n^{\#},Q_n)$. We now use the same three equalities as above, but now use that $dQ_{n,s}(u_s)\geq 0$ and $Q_n(0)\geq 0$, by definition of ${\cal F}_A^{np+}$, which then shows $r(h,Q_n)=0$.
This proves (\ref{keya1}) and thereby Theorem \ref{ThB}.

Consider  now the squared error loss $L_1(Q)=(Y-Q(X))^2$. Then,
\begin{eqnarray*}
d_{n1}(Q_n^{\#},Q_n)&=&\frac{1}{n}\sum_i \{2Y_iQ_n(X_i)-2Y_i Q_n^{\#}(X_i)+Q_n^{\#2}(X_i)-Q_n^2(X_i)\}\\
&=& \frac{1}{n}\sum_i\{2  (Q_n-Q_n^{\#})(X_i)  Y_i+Q_n^{\#2}(X_i)-Q_n^2(X_i)\}\\
&=&\frac{1}{n}\sum_i \{2(Q_n-Q_n^{\#})(X_i)(Y_i-Q_n(X_i))\\
&&\hfill +2(Q_n-Q_n^{\#})Q_n(X_i)+Q_n^{\#2}(X_i)-Q_n^2(X_i)\}\\
&=&\frac{1}{n}\sum_i 2(Q_n-Q_n^{\#})(X_i)(Y_i-Q_n(X_i))+\frac{1}{n}\sum_i (Q_n-Q_n^{\#})^2(X_i).\end{eqnarray*}
Note that the first term corresponds with $P_n\frac{d}{dQ_n}L_1(Q_n)(Q_n^{\#}-Q_n)$ and the second-order term with $P_nR_{2L_1,n}(Q_n^{\#},Q_n)$, where
$R_{2L_1,n}(Q_n^{\#},Q_n)=(Q_n^{\#}-Q_n)^2$.
We want to show that $ \frac{1}{n}\sum_i 2(Q_n-Q_n^{\#})(X_i)(Y_i-Q_n(X_i))\geq 0$.
The general equation above $P_n \frac{d}{dQ_n}L_1(Q_n)(f(h,Q_n)) \geq 0$ corresponds for the squared error loss with:
\[
-2 \frac{1}{n}\sum_i f(h,Q_n)(Y_i-Q_n(X_i)) \geq 0.\]
As we showed in general above  we have that $f(h,Q_n)$ can be chosen to be equal to $Q_n^{\#}-Q_n$.
So this proves that $n^{-1}\sum_i(Q_n^{\#}-Q_n)(X_i)(Y_i-Q_n(X_i))\leq 0$. This proves the desired result for the squared error loss.
Specifically,
\[
d_{n1}(Q_n^{\#},Q_n)\geq P_n (Q_n-Q_n^{\#})^2.\]
Clearly, $P_n  \{L_1(Q_n^{\#})-L_1(Q_n)\}^2\leq C P_n (Q_n^{\#}-Q_n)^2\leq Cd_{n1}(Q_n^{\#},Q_n)$ for some $C<\infty$. 
Lemma \ref{generalhalmle} now shows \[
d_{n1}(Q_n^{\#},Q_n)=O_P(n^{-1/2-\alpha(d)}).\]
This completes the proof of Theorem \ref{ThB} and its corollary for the squared error loss. $\Box$

We now prove Lemma \ref{lemmadndo}.
\begin{lemma}\label{dndo}
Suppose that $\int f^2_n dP_n=O_P(n^{-1/2-\alpha(d)})$ and we know that $\pl f_n\pl_v^*<M$ for some $M<\infty$. 
Then $\int f_n^2dP_0=O_P(n^{-1/2-\alpha(d)})$.
\end{lemma}
{\bf Proof:}
We have 
\begin{eqnarray*}
\int f_n^2 dP_0&=&-\int f_n^2d(P_n-P_0)+\int f_n^2 dP_n\\
&=&-\int f_n^2 d(P_n-P_0)+O_P(n^{-1/2-\alpha(d)}).
\end{eqnarray*}
We have $\int f_n^2 d(P_n-P_0)=O_P(n^{-1/2})$.
This proves that $\int f_n^2 dP_0=O_P(n^{-1/2})$. By asymptotic equicontinuity of the empirical process indexed by cadlag functions with uniformly bounded sectional variation norm, it follows now also that $\int f_n^2 d(P_n-P_0)=O_P(n^{-1/2-\alpha(d)})$. Thus, this proves that indeed
that $\int f_n^2 dP_0=O_P(n^{-1/2-\alpha(d)})$ follows from $\int f_n^2dP_n=O_P(n^{-1/2-\alpha(d)})$.
$\Box$

\section{Preservation of rate of one-step TMLE $Q_{n,\epsilon_n}$ under sampling from $P_n$}\label{AppendixC}
As the following proof demonstrates, if $\epsilon\rightarrow  L_1(Q_{\epsilon})$ is twice differentiable, and the minima $\tilde{\epsilon}_n^{\#}$ of $P_n L_1(Q_{n,\epsilon}^{\#})$ and $\epsilon_{0,n}^{\#}$ of $P_0L_1(Q_{n,\epsilon}^{\#})$ are interior minima with derivative equal to zero, then
$d_{01}(Q_{n,\epsilon_n^{\#}}^{\#},Q_0)=O_P(n^{-1/2-\alpha(d)})$, showing that indeed the one-step TMLE $Q_{n,\epsilon_n}$ also preserves the rate of convergence of the HAL-MLE under sampling from $P_n$.

Recall that $\sup_{\epsilon}\pl Q_{\epsilon}\pl_v^*< C \pl Q\pl_v^*$ for some $C<\infty$ so that the least favorable submodel preserves the bound on the sectional variation norm.

We define
\begin{eqnarray*}
\epsilon_n^{\#}&=&\arg\min_{\epsilon} P_n^{\#}L_1(Q_{n,\epsilon}^{\#})\\
\tilde{\epsilon}_n^{\#}&=&\arg\min_{\epsilon}P_n L_1(Q_{n,\epsilon}^{\#}).
\end{eqnarray*}
Under a weak regularity condition, we have $\mid \epsilon_n^{\#}-\tilde{\epsilon}_n^{\#}\mid =O_P(n^{-1/2})$.
Specifically, this can be shown as follows.
\begin{eqnarray*}
0&\leq & P_n L_1(Q_{n,\epsilon_n^{\#}}^{\#})-P_nL_1(Q_{n,\tilde{\epsilon}_n^{\#}}^{\#})\\
&=& (P_n-P_n^{\#})L_1(Q_{n,\epsilon_n^{\#}}^{\#},Q_{n,\tilde{\epsilon}_n^{\#}}^{\#})\\
&&+ P_n^{\#}L_1(Q_{n,\epsilon_n^{\#}}^{\#},Q_{n,\tilde{\epsilon}_n^{\#}}^{\#})\\
&\leq &-(P_n^{\#}-P_n)L_1(Q_{n,\epsilon_n^{\#}}^{\#},Q_{n,\tilde{\epsilon}_n^{\#}}^{\#}).\\
\end{eqnarray*}
The last term is a bootstrapped empirical process which is thus $O_P(n^{-1/2})$. 
An exact first order Tailor expansion at $\tilde{\epsilon}_n^{\#}$ allows us to write the last term as a $(P_n^{\#}-P_n)f_n^{\#}(\epsilon_n^{\#}-\tilde{\epsilon}_n^{\#})$ for a specified function $f_n$.
A second-order Tailor expansion at $\tilde{\epsilon}_n^{\#}$ of the left-hand side of this inequality and using that the first derivative at $\tilde{\epsilon}_n^{\#}$ is zero (since it is a minimum) shows that the left-hand side is a quadratic term behaving as $(\epsilon_n^{\#}-\tilde{\epsilon}_n^{\#})^2$.
Since $(P_n^{\#}-P_n)f_n^{\#}=O_P(n^{-1/2})$, this proves that  $\mid \epsilon_n^{\#}-\tilde{\epsilon}_n^{\#}\mid^2 =O_P(n^{-1/2}\mid \epsilon_n^{\#}-\tilde{\epsilon}_n^{\#}\mid)$, and thus $\mid \epsilon_n^{\#}-\tilde{\epsilon}_n^{\#}\mid =O_P(n^{-1/2})$.

Let $\epsilon_{0,n}^{\#}=\arg\min_{\epsilon}P_0 L_1(Q_{n,\epsilon}^{\#})$.
Under a weak regularity condition, we also have $\tilde{\epsilon}_n^{\#}-\epsilon_{0,n}^{\#}=O_P(n^{-1/2})$.
This is shown similarly:
\begin{eqnarray*}
0&\leq& P_0 L_1(Q_{n,\tilde{\epsilon}_n^{\#}}^{\#})-P_0 L_1(Q_{n,\epsilon_{0,n}^{\#}}^{\#})\\
&=&(P_0-P_n)L_1(Q_{n,\tilde{\epsilon}_n^{\#}}^{\#},Q_{n,\epsilon_{0,n}^{\#}}^{\#})+P_n L_1(Q_{n,\tilde{\epsilon}_n^{\#}}^{\#},Q_{n,\epsilon_{0,n}^{\#}}^{\#})\\
&\leq&- (P_n-P_0)L_1(Q_{n,\tilde{\epsilon}_n^{\#}}^{\#},Q_{n,\epsilon_{0,n}^{\#}}^{\#}).
\end{eqnarray*}
The last term is an empirical process which is thus $O_P(n^{-1/2})$.
An exact first order Tailor expansion allows us to write the last term as $(P_n-P_0)f_n^{\#}(\tilde{\epsilon}_n^{\#}-\epsilon_{0,n}^{\#})$.
A second-order Tailor expansion at ${\epsilon}_{0,n}^{\#}$ of the left-hand side of this inequality and using that the first derivative at the minimum $\epsilon_{0,n}^{\#}$ equals zero shows that the left-hand side is a quadratic term behaving as $(\tilde{\epsilon}_n^{\#}-\epsilon_{0,n}^{\#})^2$.
Since $(P_n-P_0)f_n^{\#}=O_P(n^{-1/2})$, this proves that indeed $\mid \tilde{\epsilon}_n^{\#}-\epsilon_{0,n}^{\#}\mid =O_P(n^{-1/2})$.


We will now bound $d_{01}(Q_{n,\tilde{\epsilon}_n^{\#}}^{\#},Q_0)$, where we use the latter  $\mid \tilde{\epsilon}_n^{\#}-\epsilon_{0,n}^{\#}\mid =O_P(n^{-1/2})$, and, combined with $\epsilon_n^{\#}-\tilde{\epsilon}_n^{\#}=O_P(n^{-1/2})$ this will give the desired result for $d_{01}(Q_{n,\epsilon_n^{\#}}^{\#},Q_0)$.
Using that $\epsilon_{0,n}^{\#}$ minimizes $P_0L_1(Q_{n,\epsilon}^{\#})$ and $\tilde{\epsilon}_n^{\#}$ minimizes $P_n L_1(Q_{n,\epsilon}^{\#})$ provides us with the following two  subsequent inequalities:
\begin{eqnarray*}
0&\leq &P_0 L_1(Q_{n,\tilde{\epsilon}_n^{\#}}^{\#})-P_0L_1(Q_0)\\
&=&P_0 L_1(Q_{n,\tilde{\epsilon}_n^{\#}}^{\#},Q_{n,\epsilon_{0,n}^{\#}}^{\#})\\
&&+P_0 L_1(Q_{n,\epsilon_{0,n}^{\#}}^{\#},Q_n^{\#})+P_0 L_1(Q_n^{\#},Q_0)\\
&\leq&P_0 L_1(Q_{n,\tilde{\epsilon}_n^{\#}}^{\#},Q_{n,\epsilon_{0,n}^{\#}}^{\#})+P_0 L_1(Q_n^{\#},Q_0)\\
&=&
-(P_n-P_0) L_1(Q_{n,\tilde{\epsilon}_n^{\#}}^{\#},Q_{n,\epsilon_{0,n}^{\#}}^{\#})+P_n L_1(Q_{n,\tilde{\epsilon}_n^{\#}}^{\#},Q_{n,\epsilon_{0,n}^{\#}}^{\#})\\
&&+P_0 L_1(Q_n^{\#},Q_0)\\
&\leq & -(P_n-P_0) L_1(Q_{n,\tilde{\epsilon}_n^{\#}}^{\#},Q_{n,\epsilon_{0,n}^{\#}}^{\#})+P_0 L_1(Q_n^{\#},Q_0).\\
\end{eqnarray*}
Using an exact first order Tailor expansion at $\epsilon_{0,n}^{\#}$, the first  empirical process term can be represented as  $(P_n-P_0)f_n^{\#}(\tilde{\epsilon}_n^{\#}-\epsilon_{0,n}^{\#})$, which is thus $O_P(n^{-1})$. The second term can be written as $P_0L_1(Q_n^{\#},Q_n)+P_0L_1(Q_n,Q_0)$. We have $P_0 L_1(Q_n,Q_0)=O_P(n^{-1/2-\alpha(d)})$.
We also have 
\begin{eqnarray*}
P_0L_1(Q_n^{\#},Q_n)&=&-(P_n-P_0)L_1(Q_n^{\#},Q_n)+P_n L_1(Q_n^{\#},Q_n)\\
&=&-(P_n-P_0)L_1(Q_n^{\#},Q_n)+O_P(n^{-1/2-\alpha(d)}),\end{eqnarray*}
by Theorem \ref{thnpbootmle}. The first term $(P_n-P_0)L_1(Q_n^{\#},Q_n)=O_P(n^{-1/2-\alpha(d)})$ as well, as shown earlier.
Thus, we have shown
\[P_0 L_1(Q_{n,\tilde{\epsilon}_n^{\#}}^{\#})-P_0L_1(Q_0)=O_P(n^{-1/2-\alpha(d)}).\]
Thus, we have shown $d_{01}(Q_{n,\tilde{\epsilon}_n^{\#}}^{\#},Q_0)=O_P(n^{-1/2-\alpha(d)})$ and $\epsilon_n^{\#}-\tilde{\epsilon}_n^{\#}=O_P(n^{-1/2})$.
By using that $d_{01}(Q,Q_0)$ behaves as a  square of an $L^2(P_0)$-norm, it follows trivially that this implies $d_{01}(Q_{n,\epsilon_n^{\#}}^{\#},Q_0)=O_P(n^{-1/2-\alpha(d)})$.

\section{Nonparametric bootstrap to estimate a supremum norm of the upper-bound of exact second order expansion of HAL-one-step estimator}\label{sectsupupperb}



 By Lemma \ref{lemma2} for the HAL-MLEs we have ${\bf d}_{01}(Q_n,Q_0)\leq -(P_n-P_0){\bf L}_1(Q_n,Q_0)$ and ${\bf d}_{02}(G_n,G_0)\leq -(P_n-P_0){\bf L}_2(G_n,G_0)$.
  Applying these two upper bounds for ${\bf d}_{01}(Q_n,Q_0)$ and ${\bf d}_{02}(G_n,G_0)$ to our upper bound for the exact second-order remainder 
  in the expansion $\psi_n^1-\psi_0=(P_n-P_0)D^*(Q_n,G_n)+R_{20}(Q_n,G_n,Q_0,G_0)$ yields the following conservative bound for the HAL-one-step estimator $\psi_n^1=\Psi(Q_n)+P_n D^*(Q_n,G_n)$:
\begin{eqnarray*}
 n^{1/2}\mid \psi_n^1-\Psi(Q_0)\mid& \leq & \mid  n^{1/2}(P_n-P_0)D^*(Q_n^*,G_n) \mid \\
 &&\hspace*{-5cm}
 +
 f(\sqrt{\mid n^{1/2}(P_n-P_0){\bf L}_1(Q_n,Q_0)\mid},\sqrt{ \mid n^{1/2}(P_n-P_0){\bf L}_2(G_n,G_0)\mid})\mid .
\end{eqnarray*}
Suppose that $\sup_{P,P_1\in {\cal M}}P \{{\bf L}_1(Q(P_1))-{\bf L}_1(Q(P))\}\leq c_1\in \openr^{K_1}_{\geq 0}$ and $\sup_{P,P_1\in {\cal M}}P\{{\bf L}_2(G(P_1))-{\bf L}_2(G(P))\}<c_2\in \openr^{K_2}_{\geq 0}$, then we also know that $R_2(Q_n^*,G_n,Q_0,G_0)\leq f(c)$, where $c=(c_1,c_2)$. Therefore, we assume that $f$ is chosen so that these global bounds on the loss-based dissimilarities ${\bf d}_{01}$ and ${\bf d}_{02}$ are respected even when $f$ is applied to an $x\in\openr^{K}_{\geq 0}$ with $x\geq c$ that is outside these bounds: i.e., 
$f(x)=f(\min(x,c))$.

Define the following process $X_n=(X_n(Q,G): (Q,G)\in {\cal F})$:
\begin{eqnarray*}
X_n(Q,G)&=&\mid  n^{1/2}(P_n-P_0)D^*(Q,G) \mid \\
&&+
 f(\sqrt{\mid n^{1/2}(P_n-P_0){\bf L}_1(Q,Q_0)\mid},\sqrt{ \mid n^{1/2}(P_n-P_0){\bf L}_2(G,G_0)\mid})\mid .\end{eqnarray*}
{\bf Finite sample bound:}
 Thus, we can state that \[
 \mid n^{1/2}(\psi_n^1-\psi_0)\mid \leq \mid X_n(Q_n,G_n)\mid .\]
 It remains to upper bound $X_n(Q_n,G_n)$.\nl
{\bf Proposed method of inference:}
 This finite sample bound in terms of a stochastic process evaluated at our HAL-MLEs suggests the following method for inference:
 \begin{itemize}
 \item For user supplied $\alpha_n=(\alpha_{1n},\alpha_{2n})$, determine an $x_n=(x_{1n},x_{2n})$  so that $P(d_{01}(Q_n,Q_0)>x_{1n})\leq \alpha_{1n}$ and $P(d_{02}(G_n,G_0)>x_{2n})\leq \alpha_{2n}$. For example, $x_n$ could be defined so that  $\mid (P_n-P_0)L_1(Q_n,Q_0)\mid >x_{1n}$ with probability $\alpha_{1n}$ and 
 $\mid (P_n-P_0)L_2(G_n,G_0)\mid >x_{2n}$ with probability $\alpha_{2n}$.
  \item Define
 \[ {\cal F}(x_n)=\{(Q,G)\in {\cal F}: d_{01}(Q,Q_0)<x_{1n},d_{02}(G,G_0)<x_{2n}\},\]
 and note that $P((Q_n,G_n)\in {\cal F}(x_n))\geq 1-\bar{\alpha}_n$, where $\bar{\alpha}_n\equiv \alpha_{1n}+\alpha_{2n}$.
 \item We have 
 \[
P(\mid X_n(Q_n,G_n)\mid > x)\leq  P(\pl X_n\pl_{{\cal F}(x_n)}>x)+\bar{\alpha}_n ,\]
where $\pl X\pl_{{\cal F}}=\sup_{f\in {\cal F}}\mid X(f)\mid$.
Let 
\[
F_{n,x_1}(x)=P(\pl  X_n \pl_{{\cal F}(x_1)} \leq x),\]
and
\[
q_{n,0.95}=F_{n,x_n}^{-1}(0.95+\bar{\alpha}_n) .\]
We have  $\psi_n^1\pm q_{n,0.95}/n^{1/2}$ contains $\psi_0$ with probability at least $0.95$.
\item Let $X_n^{\#}$ be the nonparametric bootstrap estimate of $X_n$:
\begin{eqnarray*}
X_n^{\#}(Q,G)&=&\mid  n^{1/2}(P_n^{\#}-P_n)D^*(Q,G) \mid \\
&&\hspace*{-2cm} +
 f(\sqrt{\mid n^{1/2}(P_n^{\#}-P_n){\bf L}_1(Q,Q_n)\mid}, \sqrt{\mid n^{1/2}(P_n^{\#}-P_n){\bf L}_2(G,G_n)\mid}) .\end{eqnarray*}
 \item Let $x_n^{\#}$ be an estimator of $x_n$ using the nonparametric bootstrap. A highly conservative method is presented below. Alternatively, one could determine $x_n^{\#}$ based on the  conservative sampling distributions $n^{1/2}(P_n^{\#}-P_n)L_1(Q_n^{\#},Q_n)$ for $d_{01}(Q_n,Q_0)$ and $n^{1/2}(P_n^{\#}-P_n)L_2(G_n^{\#},G_n)$ for $d_{02}(G_n,G_0)$. 
\item Let  \[ {\cal F}_n(x_n^{\#})=\{(Q,G)\in {\cal F}: d_{n1}(Q,Q_n)<x_{1n}^{\#},d_{n2}(G,G_n)<x_{2n}^{\#}\}\] be the nonparametric bootstrap version of ${\cal F}(x_n)$.
Let 
\[
F_{n,x_1}^{\#}(x)=P(\pl  X_n^{\#} \pl_{{\cal F}_n(x_1)} \leq x\mid (P_n:n\geq 1)).\]
Let 
\[
q_{n,0.95}^{\#}=F_{n,x_n^{\#}}^{\#-1}(0.95+\bar{\alpha}_n) .\]
\item The proposed 0.95-confidence interval is given by:
\[
\psi_n^1\pm q_{n,0.95}^{\#}/n^{1/2} .\]
\end{itemize}
{\bf Method for determining cut-offs $x_n$ for loss-based dissimilarities of HAL-MLEs:}
In order to implement the above confidence interval we  need to derive a method for determining $x_n$ and its bootstrap version $x_n^{\#}$.
We will now present a conservative definition of $x_n$ and its bootstrap estimate $x_n^{\#}$.
By our integration by parts lemma, we have
\begin{eqnarray*}
d_{01}(Q_n,Q_0)&\leq & \pl \bar{P}_n-\bar{P}_0\pl_{\infty}\pl L_1(Q_n,Q_0)\pl_v^*  \\
&\leq & 2M_1\pl \bar{P}_n-\bar{P}_0\pl_{\infty},\\
d_{02}(G_n,G_0)&\leq &\pl \bar{P}_n-\bar{P}_0\pl_{\infty}\pl L_2(G_n,G_0)\pl_v^* \\
&\leq& 2M_2\pl \bar{P}_n-\bar{P}_0\pl_{\infty}.
\end{eqnarray*}
Let  $T_{n}(x)\equiv P(\pl n^{1/2}(\bar{P}_n-\bar{P}_0)\pl_{\infty}>x)$ so that
\[
P(d_{01}(Q_n,Q_0)>x)\leq T_n( 1/2 n^{1/2}M_1^{-1} x).\]
Given $\alpha_{1n},\alpha_{2n}$ we can select $x_{1n}=2M_1n^{-1/2}T_n^{-1}(\alpha_{1n})$ and $x_{2n}=2M_2n^{-1/2}T_n^{-1}(\alpha_{2n})$.
Then, we have that $P(d_{01}(Q_n,Q_0)>x_{1n})\leq \alpha_{1n}$ and $P(d_{02}(G_n,G_0)>x_{2n})\leq \alpha_{2n}$.
We also know from empirical process theory that the supremum norm of $n^{1/2}(\bar{P}_n-\bar{P}_0)$ has an exponential tail $\exp(-C x)$ so that
$T_n^{-1}(\alpha_{1n})$ behaves as $\log \alpha_{1n}^{-1}$, and the same applies to $T_n^{-1}(\alpha_{2n})$. This shows that one can select $\alpha_{1n}$ and $\alpha_{2n}$ as numbers that converge to zero at a polynomial rate (e.g. $\alpha_{1n}=n^{-1}$), while still preserving that $\max(x_{1n},x_{2n})\rightarrow 0$ at rate $\log n/n^{1/2}$. 
Since $n^{1/2}(\bar{P}_n-\bar{P}_0)$ is an empirical process indexed by a class of indicators, we can consistently and robustly estimate the distribution of $\pl n^{1/2}(\bar{P}_n-\bar{P}_0)\pl_{\infty}$  with the nonparametric bootstrap \citep{vanderVaart&Wellner96}, which, for completeness,  is  stated in the following lemma.
\begin{lemma}
We have that (uniformly in $x$) $T_{n}(x)\rightarrow T_0(x)$, where $T_0(x)=P(\pl \bar{X}_0\pl_{\infty}>x)$, and 
$\bar{X}_0$ is the limit Gaussian process of $\bar{X}_n=n^{1/2}(\bar{P}_n-\bar{P}_0)$.
We have
\begin{eqnarray*}
P(n^{1/2}d_{01}(Q_n,Q_0)>x)&\leq& T_{n}(1/2n^{1/2}x/M_1)\\
P(n^{1/2}d_{02}(G_n,G_0)>x)&\leq & T_{n}(1/2n^{1/2}x/M_2).
\end{eqnarray*}
Let $T_n^{\#}(x)=P(\pl n^{1/2}(\bar{P}_n^{\#}-\bar{P}_n)\pl_{\infty}>x\mid (P_n:n\geq 1))$. Then, uniformly in $x$,
$T_n^{\#}(x)\rightarrow T_0(x)$ as $n\rightarrow\infty$.
\end{lemma}
We estimate $x_n$ with its nonparametric bootstrap analogue: 
\begin{eqnarray*}
x_{1n}^{\#}&=&2M_1 n^{-1/2}T_n^{\#-1}(\alpha_{1n})\\
x_{2n}^{\#}&=&2M_2n^{-1/2}T_n^{\#-1}(\alpha_{2n}).
\end{eqnarray*}

We note that our nonparametric bootstrap estimator of the sampling distribution of $X_n$ and $x_n$ only relies on how well the nonparametric bootstrapped empirical process indexed by Donsker class (functions with uniformly bounded sectional variation norm) and class of indicators approximates the sampling distribution of the empirical process. Therefore, the asymptotic consistency of this nonparametric bootstrap method  follows straightforwardly from the asymptotic consistency of the nonparametric bootstrap for these empirical processes. This is presented in the following theorem. 

 \begin{theorem}
 Consider definitions $X_n$,
 $x_{1n}=2M_1n^{-1/2}T_n^{-1}(\alpha_{1n})$ and $x_{2n}=2M_2n^{-1/2}T_n^{-1}(\alpha_{2n})$. 
 Let $\alpha_n$ converge to zero at a rate $n^{-p}$ for some finite $p$. Then, $x_n$ converges to zero at rate $\log n/n^{1/2}$.
 
{\bf Finite sample oracle confidence interval:} We have
 \[
 \mid n^{1/2}(\psi_n^1-\psi_0)\mid \leq \mid X_n(Q_n,G_n)\mid,\]
 and
 \[
P(\mid X_n(Q_n,G_n)\mid > x)\leq  P(\pl X_n\pl_{{\cal F}(x_n)}>x)+\bar{\alpha}_n .\]
Let 
\[
F_{n,x_1}(x)=P(\pl  X_n \pl_{{\cal F}(x_1)} \leq x),\]
and
\[
q_{n,0.95}=F_{n,x_n}^{-1}(0.95+\bar{\alpha}_n) .\]
Then,  $\psi_n^1\pm q_{n,0.95}/n^{1/2}$ contains $\psi_0$ with probability at least $0.95$. 

{\bf Weak convergence of process $X_n$:}
  Let $Z_n$ be the empirical process $n^{1/2}(P_n-P_0)$ indexed by the class of functions 
 \begin{equation}\label{calf1}
 {\cal F}_1\equiv \{D^*(Q,G),L_1(Q,Q_0),L_2(G,G_0): (Q,G)\in {\cal F}\}.
 \end{equation}
  $X_n$ is a continuous function of $Z_n$:
\[
X_n(Q,G)=g(Z_n)\equiv \mid Z_n(D^*(Q,G))\mid +f(\mid Z_n(L_1(Q,Q_0))\mid,\mid Z_n(L_2(G,G_0))\mid) .\]
We know that $Z_n\Rightarrow_d Z_0$ for a Gaussian process in $\ell^{\infty}({\cal F}_1)$.
The continuous mapping theorem shows that $X_n\Rightarrow_d X_0=g(Z_0)$ , where 
$X_0$ is a simple function of $Z_0$ defined by
\[
X_0(Q,G)= \mid Z_0(D^*(Q,G))\mid +f( \mid Z_0(L_1(Q,Q_0))\mid, \mid Z_0(L_2(G,G_0))\mid) .\] 
Therefore, $\pl X_n\pl_{{\cal F}(x)}\Rightarrow \pl X_0\pl_{{\cal F}(x)}$, uniformly in $x$.
In particular, uniformly in $x$
\[
P(\pl X_n\pl_{{\cal F}(x_n)}>x)\rightarrow P(N(0,\sigma^2_0)>x).\]

{\bf Weak convergence of nonparametric bootstrap process $X_n^{\#}$:}
Let $X_n^{\#}$, $x_n^{\#}$, ${\cal F}_n(x_n^{\#})$ be the nonparametric bootstrap version of $X_n$, $x_n$ and ${\cal F}(x_n)$ defined above.
 We have that uniformly in $x$,  conditional on $(P_n:n\geq 1)$, $\pl X_n^{\#}\pl_{{\cal F}(x)}\Rightarrow_d \pl X_0\pl_{{\cal F}(x)}$.
 In particular,
 \[
 P(\pl X_n^{\#}\pl_{{\cal F}(x_n^{\#})}>x\mid (P_n:n\geq 1))+\bar{\alpha}_n\rightarrow P(\mid N(0,\sigma^2_0)\mid >x).\]
 
 {\bf Nonparametric bootstrap estimate of finite sample confidence interval and its asymptotic consistency:}
 Let 
\[
F_{n,x_1}^{\#}(x)=P(\pl  X_n^{\#} \pl_{{\cal F}_n(x_1)} \leq x\mid (P_n:n\geq 1)).\,\]
and
\[
q_{n,0.95}^{\#}=F_{n,x_n^{\#}}^{\#-1}(0.95+\bar{\alpha}_n) .\]

 Then,  $q_{0.95,n}^{\#}(x_n)\rightarrow 1.96$ so that the confidence interval $\psi_n^*\pm q_{0.95,n}^{\#}/n^{1/2}$ is an asymptotic $0.95$-confidence interval.
 \end{theorem}
 
 The width of this proposed confidence interval is a function of the user supplied $\alpha_n=(\alpha_{1n},\alpha_{2n})$. To remove this choice from consideration, one might determine the $\alpha_n$ that minimizes $q_{n,0.95}^{\#}=F_{n,x_n^{\#}}^{\#-1}(0.95+\bar{\alpha}_n) $, where we note that $x_n^{\#}$ depends on $(\alpha_{1n},\alpha_{2n})$.

 \subsection{Simplified conservative version}
 The implementation of the computation of the supremum over $(Q,G)$ of $X_n(Q,G)$ can be computationally challenging.
 Therefore, here we pursue an easy to compute and robust conservative version of this approach. Consider a sequence of cut-offs  $x_n=(x_{1n},x_{2n})$ for which  $P(d_{01}(Q_n,Q_0)\leq x_{1n}, d_{02}(G_n,G_0)\leq x_{2n})\geq 1-\bar{\alpha}_{n}$, where $\bar{\alpha}_n\rightarrow 0$.
 Here one could conservatively replace $d_{01}(Q_n,Q_0)$ and $d_{02}(G_n,G_0$ by their empirical process bounds
 $\mid (P_n-P_0)L_1(Q_n,Q_0)\mid $ and $\mid (P_n-P_0)L_2(G_n,G_0)\mid$, respectively.
 
  The basic idea is to replace the second-order remainder
 by the upper bound $n^{1/2}(f(x_n^{1/2})$ (which holds with probability at least $1-\bar{\alpha}_n$), while keeping the leading term $\mid n^{1/2}(P_n-P_0)D^*(Q_n,G_n)\mid$ for what it is. We think the latter makes sense since we suspect that the supremum of $\mid n^{1/2}(P_n-P_0)D^*(Q,G)\mid $ over all $(Q,G)\in {\cal F}_{x_n}$ will be close to its evaluation at $Q_n,G_n$.
 
 Let $I_n$ be the indicator that $(Q_n,G_n)\in {\cal F}_{x_n}$ so that $P(I_n=1)\geq 1-\bar{\alpha}_n$, and let $\tilde{Z}_n\equiv \mid n^{1/2}(P_n-P_0)D^*(Q_n,G_n)\mid$.
We have the following bounding for the tail probability of $\mid n^{1/2}(\psi_n^1-\psi_0)\mid$:
\begin{eqnarray*}
P(\mid n^{1/2}(\psi_n^1-\psi_0)\mid >x)&\leq&P( \tilde{Z}_n +n^{1/2} f(x_{1n}^{1/2},x_{2n}^{1/2})>x,I_n=1)\\
&&\hspace*{-4cm}+
P(\tilde{Z}_n +n^{1/2} f(d_{01}(Q_n,Q_0)^{1/2},d_{02}(G_n,G_0)^{1/2})>x,I_n=0)\\
&\leq& P( \tilde{Z}_n +n^{1/2} f(x_{1n}^{1/2},x_{2n}^{1/2})>x)+\bar{\alpha}_n\\
&=&
P( \tilde{Z}_n>x-n^{1/2} f(x_{1n}^{1/2},x_{2n}^{1/2}))+\bar{\alpha}_n.
\end{eqnarray*}
 Let $\Phi_n(x)=P(\tilde{Z}_n>x)$.
Let $\tilde{q}_{n,0.05} $ be the solution in $x$ of 
\[
\Phi_n(x-n^{1/2}f(x_{1n}^{1/2},x_{2n}^{1/2}))+\bar{\alpha}_n=0.05.\]
Thus,
\[
\tilde{q}_{n,0.05}=n^{1/2}f(x_n^{1/2})+\Phi_n^{-1}(0.05-\bar{\alpha}_n).\]
For this choice we have $P(\mid n^{1/2}(\psi_n^1-\psi_0)\mid > \tilde{q}_{n,0.05})\leq 0.05$ so that 
\[
\psi_n^1\pm \tilde{q}_{n,0.05}/n^{1/2}\]
is a finite sample $>0.95$-confidence interval.

Let $x_n^{\#}$ be the above presented bootstrap estimator of $x_n$, $\tilde{Z}_n^{\#}=\mid n^{1/2}(P_n^{\#}-P_n)D^*(Q_n^{\#},G_n^{\#})\mid$, and 
$\Phi_n^{\#}(x)=P(\tilde{Z}_n^{\#}>x\mid (P_n:n\geq 1))$ be the bootstrap estimator of $\Phi_n$.
Then, our bootstrap estimator of this confidence interval is given by:
\[
\psi_n^1\pm \tilde{q}_{n,0.05}^{\#}/n^{1/2},\]
where 
\[
\tilde{q}_{n,0.05}^{\#}=n^{1/2}f(x_n^{\#1/2})+\Phi_n^{\#-1}(0.05-\bar{\alpha}_n).\]

 In order for this confidence interval to be asymptotically sharp one will need to make sure that $x_n$ converges to zero faster than $n^{-1/2}$. Our conservative bound for $x_n^{\#}$  in the previous subsection is too large. Instead, we might estimate $x_n^{\#}$ so that 
 \[
 P( A_{1}(P_n^{\#},P_n)<x_{1n}^{\#},A_2(P_n^{\#},P_n)<x_{2n}^{\#}\mid (P_n:n\geq 1))\geq 1-\bar{\alpha}_n,\]
 where $A_1(P_n^{\#},P_n)\equiv \mid (P_n^{\#}-P_n)L_1(Q_n^{\#},Q_n)\mid$ and $A_2(P_n^{\#},P_n)\equiv \mid (P_n^{\#}-P_n)L_2(G_n^{\#},G_n)\mid$.
 Then, $x_n^{\#}$ converges to zero at a rate faster than $n^{-1/2}$ for slowly converging $\bar{\alpha}_n$.

\end{document}